\documentclass[review]{elsarticle}

\usepackage[a4paper,left=20mm,right=20mm,top=20mm,bottom=20mm]{geometry}
\usepackage{lineno,hyperref}
\modulolinenumbers[5]

\journal{Computers \& Chemical Engineering}

\usepackage{subcaption}
\usepackage{amsmath}
\usepackage{float}
\usepackage{booktabs}
\usepackage{amsthm}
\usepackage{todonotes}
\newtheorem{theorem}{Theorem}[section]
\theoremstyle{definition}
\newtheorem{definition}{Definition}
\newtheorem{assumption}{Assumption}
\newtheorem{corollary}{Corollary}[theorem]
\newtheorem{lemma}{Lemma}
\newtheorem{remark}{Remark}
\usepackage{verbatim}
\usepackage{rotating}
\usepackage{amssymb}
\usepackage{tabularx}

\usepackage{algorithm}
\usepackage{algpseudocode}
\algnewcommand\NonumProcedure[1]{%
  \item[\textbf{procedure}] \textsc{#1}
}
\algnewcommand\NonumEndProcedure[1]{%
  \item[\textbf{end procedure}] \textsc{#1}
}
\algnewcommand{\Inputs}[1]{%
  \State \textbf{Inputs:}
  \Statex \hspace*{\algorithmicindent}\parbox[t]{.8\linewidth}{\raggedright #1}
}
\algnewcommand{\Initialize}[1]{%
  \State \textbf{Initialize:}
  \Statex \hspace*{\algorithmicindent}\parbox[t]{1.2\linewidth}{\raggedright #1}
}

\newcommand\CONDITION[2]%
  {\begin{tabular}[t]{@{}l@{}l@{}}
     #1&#2
   \end{tabular}
  }
\algdef{SE}[WHILE]{While}{EndWhile}[1]%
  {\algorithmicwhile\ \CONDITION{#1}{\ \algorithmicdo}}%
  {\algorithmicend\ \algorithmicwhile}
\algdef{SE}[FOR]{For}{EndFor}[1]%
  {\algorithmicfor\ \CONDITION{#1}{\ \algorithmicdo}}%
  {\algorithmicend\ \algorithmicfor}
\algdef{S}[FOR]{ForAll}[1]%
  {\algorithmicforall\ \CONDITION{#1}{\ \algorithmicdo}}
\algdef{SE}[REPEAT]{Repeat}{Until}{\algorithmicrepeat}[1]%
  {\algorithmicuntil\ \CONDITION{#1}{}}
\algdef{SE}[IF]{If}{EndIf}[1]%
  {\algorithmicif\ \CONDITION{#1}{\ \algorithmicthen}}%
  {\algorithmicend\ \algorithmicif}%
\algdef{C}[IF]{IF}{ElsIf}[1]%
  {\algorithmicelse\ \algorithmicif\ \CONDITION{#1}{\ \algorithmicthen}}

\DeclareMathOperator*{\argmax}{argmax}
\DeclareMathOperator*{\argmin}{argmin}
\usepackage{color}
\usepackage{setspace}
\usepackage{caption}
\definecolor{frenchblue}{rgb}{0.0, 0.45, 0.73}
\definecolor{skyblue}{rgb}{0.3010, 0.7450, 0.9330}
\definecolor{cyan1}{rgb}{0, 0.96, 1}
\definecolor{gray}{rgb}{0.5, 0.5, 0.5}
\definecolor{bleudefrance}{rgb}{0.19, 0.55, 0.91}

\hypersetup{colorlinks = true, allcolors = bleudefrance}

\def\subspace{region }
\def\subspaces{regions }
\def\Subspace{Region }

\biboptions{authoryear}

\makeatletter
\def\BState{\State\hskip-\ALG@thistlm}
\makeatother

\begin{document}

\begin{frontmatter}

\title{A deterministic optimization algorithm for mixed-integer nonlinear bi-objective programming with applications to molecular and process design}
\tnotetext[mytitlenote]{Corresponding author. E-mail addresses: c.adjiman@imperial.ac.uk}

\author[1]{Lauren Ye Seol Lee}
\author[2]{George Jackson}
\author[2]{Amparo Galindo}
\author[2]{Claire S. Adjiman$^{*,}$ }


\address[1]{Department of Chemical Engineering, University College London, Torrington Pl, London WC1E 7JE, United Kingdom}
\address[2]{Department of Chemical Engineering, Sargent Centre for Process Systems Engineering, Imperial College London, South Kensington Campus, London SW7 2AZ, United Kingdom}

\begin{abstract}
Many practical multiobjective optimization (MOO) problems include discrete decision variables and/or nonlinear model equations and exhibit disconnected or smooth but nonconvex Pareto surfaces. Scalarization methods, such as the weighted-sum and sandwich (SD) algorithms, are common approaches to solving MOO problems but may fail on nonconvex or discontinuous Pareto fronts. In the current work, motivated by the normal boundary intersection (NBI) method and the SD algorithm, we present SDNBI, a new algorithm for bi-objective optimization (BOO) designed to address the theoretical and numerical challenges associated with the reliable solution of general nonconvex and discrete BOO problems. The main improvements in the algorithm are the efficient exploration of the nonconvex regions of the Pareto front and, uniquely, the early identification of regions where no additional Pareto solutions exist. The performance of the SDNBI algorithm is assessed based on the accuracy of the approximation of the Pareto front constructed over the disconnected nonconvex objective domains. The new algorithm is compared with two MOO approaches, the modified NBI method and the SD algorithm, using benchmark problems and two molecular and process design case studies. The SDNBI algorithm outperforms the modified NBI and SD algorithms in tackling this collection of convex, nonconvex-continuous, and integer / mixed-integer nonlinear problems, both in terms of computational cost and of the overall quality of the Pareto-optimal set, suggesting that it is a promising alternative for solving challenging BOO problems.
\end{abstract}


\end{frontmatter}


\section{Introduction}
Many real-world problems in operational research, engineering, biology, and chemistry can be characterized as multiobjective optimization (MOO) problems, due to the presence of multiple objective functions that cannot be easily merged into a single metric. The objective functions are typically conflicting, so that MOO problems do not have a single solution that is optimal for all objectives simultaneously. Instead, a set of points, each of which corresponds to a trade-off between the objectives, is sought, such that an improvement in one objective can only be achieved through a sacrifice in another \citep{Miettinen1998multi}. These points are commonly known as Pareto-optimal solutions (or Pareto points, or nondominated points). In most practical applications, it is not possible to derive an analytical expression that describes the loci of these points in the space of the objective functions (the Pareto frontier) \citep{deb2001multi}. As a result, the development of MOO algorithms that can construct an approximation of the Pareto frontier efficiently is of special interest. This remains relevant in process systems engineering, where multiobjective formulations are increasingly used to quantify trade-offs among economic, environmental, safety, and operational criteria in process design and operation, as well as in experimental process development \citep{isella2025integrating,rahbari2025multi,zhang2025multi}. The main criteria in the development of most MOO algorithms are (1) to find Pareto points within reasonable computational time and (2) to generate these points so that they are as diverse as possible and distributed evenly along the Pareto frontier. 

Among the several classes of MOO solution approaches, such as stochastic methods \citep{zimmermann1978fuzzy,huang2006interactive,serafini1994simulated}, evolutionary algorithms \citep{li2019comparison,deb2000fast,deb2013evolutionary} and exact solution strategies \citep{belotti2013branch,przybylski2017multi}, one of the most widely used concepts is scalarization, in which the multiple objective functions are combined into a single objective function through weight parameters and/or additional constraints, making it possible to identify solutions of the MOO problem by solving a series of single objective optimization (SOO) problems using standard optimization methods. A popular scalarization method is the weighted-sum method \citep{marler2004}, in which a set of Pareto-optimal solutions is generated by varying the weights assigned to the multiple objective functions. The weighted-sum method is easy to implement and the problem is of the same degree of computational complexity as the original MOO problem since there are no additional constraints involved. However, this approach suffers from the theoretical limitations that it is not possible to capture any Pareto-optimal solutions on nonconvex regions of the Pareto frontier \citep{Miettinen1998multi} and that an even distribution of the weights does not guarantee an even distribution of solutions on the Pareto front. The $\epsilon$-constraint method, on the other hand, does not require any convexity assumptions on the Pareto frontier, but it has the drawback of being sensitive to the range of values of the objective functions in the feasible region, possibly resulting in infeasible SOOs or repeated Pareto-optimal solutions \citep{ehrogott2005moo}. Even when their range is known, the choice of $\epsilon$ values is not straightforward and it may be necessary to include additional constraints to remove regions of the feasible set. As a result, the computational cost of this approach is strongly dependent on the number of objective functions. Extensions such as the augmented $\epsilon$-constraint method have been developed to improve the systematic generation of Pareto-optimal solutions \citep{mavrotas2009effective} and continue to be used for mixed-integer process and energy-system design problems \citep{schricker2026multi}.

To deal with these shortcomings, several scalarization methods have been proposed with a special focus on the development of a systematic way to determine scalarization parameters and to capture Pareto points from nonconvex regions of the Pareto surfaces. \cite{kim2005adaptive} proposed the adaptive weighted-sum algorithm. To derive an initial representation of the Pareto frontier, the algorithm starts by generating some solution points using the weighted-sum method with a (large) uniform step size in the weighting coefficients. Thereafter, the algorithm continues to solve the weighted-sum scalarization with additional upper bounds on the objective functions for the regions that need further refinement. These regions are determined by computing the distance between neighboring solutions, to control the distribution and uniformity of the Pareto points. With this method, it is possible to investigate nonconvex regions of the Pareto front. Furthermore, the approach provides a systematic way of obtaining more uniformly distributed Pareto points. However, the algorithm can generate dominated solutions 
and it cannot be guaranteed to cover the entire set of nondominated points.

\cite{DASDENNIS1998} introduced another scalarization approach, the normal boundary intersection (NBI) method, for identifying the solutions of nonconvex MOO problems. The method generates a well-spaced set of points on the boundary of the space of objective function values by iteratively solving a scalarized subproblem with respect to reference points placed on the convex hull of the individual minima (CHIM) of the objectives. In contrast to weighted-sum methods where problem size is preserved in the subproblems, a new vector of equality constraints is introduced, resulting in a size increase. The performance of the method was tested for large-scale and highly constrained optimization problems in the area of optimal control \citep{logist2010fast}. While the method was found to produce evenly distributed solutions, its focus is on identifying the boundary of the feasible objective space, which may include some points dominated by others. Furthermore, it may be difficult to achieve convergence when there is no feasible solution that satisfies the equality constraints for the given choice of scalarization parameter values. Following similar concepts, \cite{ismail2002effective} proposed the normal constraint (NC) method as an alternative to the NBI method to avoid getting dominated solutions and improve convergence behavior. Instead of introducing equality constraints in the subproblems, as in the NBI method, a set of inequality constraints is used to reduce the size of the feasible region. \cite{messac2003normalized} further modified the NC method by incorporating a so-called Pareto filter, which is used to eliminate non-Pareto-optimal solutions. However, these normal constraint methods share the same drawbacks as the NBI method because they rely on the CHIM as a reference line or plane. Further adjustments were made by \cite{messac2004normal} and \cite{shao2007finding}, with the aim of reducing the likelihood of missing regions of the nondominated set. \cite{messac2004normal} used an extension of the reference plane defined by the CHIM, so that the entire Pareto set can be enclosed by a hypercube. These studies indicate that non-connectivity of the Pareto surface and the presence of isolated Pareto-optimal points can cause severe problems for the NBI method, whereas the WS and NC methods are less likely to be affected by the disconnected nature of the Pareto set. 
In further work on the NBI method, \cite{shao2016discrete} suggested an extension of an approach initially designed for the solution of multiobjective linear programming (MOLP) problems \citep{shao2007finding}, in which they combined the global shooting method \citep{benson1997towards} with the NBI method to guarantee coverage of the Pareto frontier as well as an even distribution of the Pareto points. More recently, adaptive NBI-based strategies have been developed to reduce the computational effort associated with generating Pareto front approximations. \cite{mores2023multi}, for example, combined NBI with a divide and conquer strategy in which a significance criterion determines whether individual regions of the Pareto frontier require further refinement.

While the methods discussed so far are designed to identify evenly distributed points on the Pareto frontier, the development of Benson’s outer approximation algorithm \citep{benson1998outer} has motivated the emergence of another class of algorithms designed to generate outer and/or inner approximations of the Pareto frontier. The majority of such methods have been focused on obtaining polyhedral approximations to ensure an even spread of points over the Pareto surface. \cite{shao2008approximating} and \cite{ehrgott2011approximation} presented improvements to Benson’s algorithm for the solution of MOLP problems by introducing the sandwich (SD) algorithm so that the Pareto frontier is located between the inner and outer approximations. The SD algorithm and its variants have several useful properties. As the real Pareto surface is sandwiched between the inner and outer approximation, an upper bound on the approximation error can be calculated; this provides valuable information on the accuracy of the Pareto-optimal solutions generated during the course of the algorithm. Furthermore, the approximation quality can be improved efficiently by adding more solutions to the region where the error bound is largest. However, the method is limited to MOLP problems and cannot easily be applied to general convex MOO problems. 
Modified versions of the SD algorithm were proposed by \cite{Solanki1993, craft2006approx, Rennen2011, klamroth2003unbiased, bokrantz2013vertex} to address general convex MOO problems. 
\cite{Solanki1993} extended the SD algorithm to make it applicable to higher-dimensional MOLP problems
. At each iteration, a weighted-sum subproblem is solved for the facet that exhibits the largest distance to a point constrained to the outer approximation. Within the algorithm, the accuracy of the Pareto front approximation is improved by recursively constructing inner and outer approximations through the generation of a convex hull and the identification of supporting hyperplanes. Although this method is directed towards the solution of MOLP problems, it can also be applied to general convex MOO problems. \cite{klamroth2003unbiased} used two separate algorithms for generating the inner and outer approximations of the Pareto set and combined them by alternatively solving one iteration of each algorithm without using the weighted-sum method. A major feature of their algorithm is the use of vertices as representative points to construct an inner approximation that maximizes the distance from the current inner approximation toward the outer approximation. Other authors have proposed further developments of this approach, refining the method for calculating the error between inner and outer approximations \citep{craft2006approx} or modifying the way in which parameters are chosen \citep{Rennen2011}. \cite{craft2006approx} calculated the error by considering the hyperplanes that pass through the corner points of a facet rather than by solving a linear programming problem. The next weight vector is taken from a linear combination of the weight vectors used to obtain the corner points of the facet. \cite{Rennen2011} further improved the algorithm by incorporating dummy points to ensure that all facets have a non-negative normal, which is crucial for the weighted-sum method. To reduce the computational expense, which increases exponentially with the number of objectives, \cite{bokrantz2013vertex} suggested employing a vertex enumeration method. 

While sandwich-type algorithms have provided a promising route to solving MOO problems, their applicability is often limited to convex MOO problems as they rely on weighted-sum subproblems. Accordingly, SD algorithms cannot achieve complete coverage of the regions of the Pareto frontier that are nonconvex, disconnected, or consist of discrete points.

In this study, we address these limitations by introducing a new scalarization scheme and a corresponding solution algorithm, the SDNBI algorithm, that generates a high-quality approximation of nonconvex and/or discontinuous Pareto fronts in bi-objective optimization (BOO) problems. Such problems may contain nonconvex functions and involve a mixed set of continuous and integer decision variables. By high quality, we mean that the algorithm generates a set of well-distributed points across the Pareto front, including in nonconvex regions, and that we provide, within an arbitrary tolerance, an enclosure for the Pareto front which consists of piecewise linear inner and outer approximations as well as regions that are identified as being devoid of Pareto points. 

The development of the approach is focused on the interplay between the modified NBI (mNBI) method suggested by \cite{shukla2007normal} and the SD algorithm, exploiting the strengths of both. We present three main modifications and we investigate their theoretical properties: 1) the validity of the inner and outer approximations derived from the solution of the scalarized subproblem and convex hull generation, 2) the completeness of decomposing of an objective search space based on the convexity of the Pareto front, and 3) the effectiveness of modifying the subproblems in avoiding unnecessary search steps around the disconnected or isolated portions of the Pareto front. To assess the performance of the proposed SDNBI algorithm, numerical tests are conducted for five bi-objective benchmark problems (MOP1, SCH2, TNK, ZDT3, and ZDT5) and compared with the results produced by the original SD algorithm and the mNBI method. The performance of the different algorithms is compared based on reliability and efficiency criteria. While not presented in this paper, it is worth noting that an early version of the SDNBI algorithm has also been successfully applied to challenging molecular-design problems, demonstrating its broader applicability and potential \citep{Lee2020}.

The remainder of this article is organized as follows. In section \ref{sec:Preliminary}, all relevant preliminaries, including the notation, definitions and the formulation of the problems are introduced. In section \ref{sec:Motivation}, we describe the NBI, mNBI and SD algorithm and their properties. In section \ref{sec:SDNBI}, we propose the new algorithm and highlight its characteristics. We perform numerical experiments in Sections \ref{sec:benchmark} and \ref{sec:results} to investigate algorithmic efficiency. Finally, we state the main conclusions, also including future perspectives for this research in Section \ref{sec:conclusion}.


\section{Preliminaries}\label{sec:Preliminary}
The generic mathematical formulation of the MOO problem is defined as:

\begin{equation}\tag{MOP}\label{eq:ch6:2-MOP}
\begin{array}{ccl}
& \text { minimize } & f_{1}(\boldsymbol{x}), \ldots, f_{m}(\boldsymbol{x}) \\
& \text { subject to } & \boldsymbol{x} \in X:=\{\boldsymbol{x} \in \mathbb{R}^{n_1} \times \mathbb{N}^{n_2} \ | \ \boldsymbol{g(x)} \leq \boldsymbol{0},  \boldsymbol{h(x)}=\boldsymbol{0}\},
\end{array}
\end{equation} 
where $f_j: \mathbb{R}^{n} \rightarrow \mathbb{R}, \ j=1,2, ...,m$, are the $m$ objective functions, $\boldsymbol{x}$ is an $n$-dimensional vector of variables ($n=n_1+n_2$), which consists of $n_1$ continuous variables and $n_2$ integer variables, $X\neq\emptyset$ is an nonempty feasible set, $\boldsymbol{g(x)}$ is a $p$-dimensional vector of inequality constraints and $\boldsymbol{h(x)}$ is a $q$-dimensional vector of equality constraints, $q\leq n$.
Note that we use $\mathbb{R}$ to denote the set of real numbers, $\mathbb{R}_{+}$ to denote the set of non-negative real numbers, and $\mathbb{R}_{++}$ to denote the set of positive real numbers. 

We make two assumptions about the nature of the problem. 
\begin{assumption} 
A constraint qualification holds for the weighted sum scalarization of problem \eqref{eq:ch6:2-MOP} whenever the discrete variables are relaxed or fixed, so that the Karush-Kuhn-Tucker (KKT) optimality conditions apply.
\end{assumption}
\begin{assumption} \label{as:boundedset}
There exist $U\in \mathbb{R}$ and $L \in \mathbb{R}$ such that  $L\le f_j(\boldsymbol{x}) \le U, \forall j\in \{1,2, ...,m\}, \forall \boldsymbol{x}\in X$.
\end{assumption}

When the objective functions conflict with each other, no single solution can simultaneously minimize all scalar objective functions. Therefore, it is necessary to introduce a new notion of optimality or Pareto efficiency \citep{ehrgott2005book}.

\begin{definition}\label{ObjectiveSpace}
The \textit{objective space} is defined as the set of all feasible combinations of the objective function values, i.e., $Z=\{ \boldsymbol{f}(\boldsymbol{x}) \mid \boldsymbol{x}\in X \}$
\end{definition}
By Assumption \ref{as:boundedset}, the set $Z$ is a closed bounded set. 

\begin{definition}\label{ParetoPoint}
A point $\boldsymbol{x}^{*} \in X$ is called \textit{weakly efficient} (or \textit{weakly Pareto-optimal}) if there exists no feasible point $\boldsymbol{x} \in X$ such that $f_j(\boldsymbol{x}) < f_j(\boldsymbol{x}^{*})$ for all $j=1, ..., m$. A point $\boldsymbol{x}^{*} \in X$ is called an \textit{efficient solution} (or \textit{Pareto-optimal}) if there exists no feasible point $\boldsymbol{x} \in X$ such that $f_j(\boldsymbol{x}) \leq f_j(\boldsymbol{x}^{*})$ for all $j = 1, ..., m$, and $f_{j'}(\boldsymbol{x}) < f_{j'}(\boldsymbol{x}^{*})$ for at least one $j' \in \{ 1,...,m \}$. If $\boldsymbol{x}^{*}$ is an efﬁcient solution, the point in objective space $\boldsymbol{z} = \boldsymbol{f}(\boldsymbol{x}^{*})$ is referred to as a \textit{nondominated solution} or \textit{Pareto point} in the objective space $Z$. 
\end{definition}

\begin{definition} \citep{marler2004} 
The set of efficient solutions $\Omega^*$ is the \textit{complete set of efficient solutions} if it contains all possible efficient solutions $\boldsymbol{x}^{*}$. The set $\Psi^*$ is the \textit{efficient frontier} (or \textit{Pareto front}) formed by all nondominated solutions, $\Psi^*=\{\boldsymbol{f}(\boldsymbol{x}^{*})\ | \ \boldsymbol{x}^*\in \Omega^*\}$. 
\end{definition}

\begin{definition} \citep{marler2004} 
Any subset of the set of efficient solutions $\boldsymbol{X_E}\subset \Omega^*$ is an approximation of the complete set $\Omega^*$ and has a corresponding approximation of the efficient frontier $\boldsymbol{Z_E}=\{\boldsymbol{f}(\boldsymbol{x}^{*}) \ | \ \boldsymbol{x}^{*}\in \boldsymbol{X_E} \}$. Note that $\boldsymbol{Z_E}\subset \Psi^*$.
\end{definition}

\begin{definition} 
The ideal objective vector or utopia point $\boldsymbol{f}^{id}$ is defined as the objective vector whose components are the optimal objective function values of each single-objective problem, i.e., $f^{id}_j=\underset{\boldsymbol{x}\in \Omega^*}{\operatorname{min}} \ f_{j}(\boldsymbol{x}),\ j=1,...,m$. The objective vector $\boldsymbol{f}^{nd}$ is defined as the nadir objective vector and is such that $f^{nd}_{j}=\underset{\boldsymbol{x}\in \Omega^*}{\operatorname{max}} \ f_{j}(\boldsymbol{x}),\ j=1,...,m$.
\end{definition}

\begin{definition}\label{df:supported} 
Let $\boldsymbol{x}^{*}\in \boldsymbol{X_{E}}$. If there is a vector $\boldsymbol{\lambda}\in\mathbb{R}^{m}_{++}$ such that $\boldsymbol{x}^{*}$ is an optimal solution to $\underset{\boldsymbol{x}\in X}{\operatorname{min}} \ \boldsymbol{{\lambda}}^{\top}\boldsymbol{f}(\boldsymbol{x})$, then $\boldsymbol{x}^{*}$ is called a \textit{supported efficient solution} and  $\boldsymbol{z}=\boldsymbol{f}(\boldsymbol{x}^*)$ is called a \textit{supported nondominated objective vector}.
\end{definition}


\begin{definition}\label{df:convexPS}
A set of Pareto points is said to be \textit{convex}, if all nondominated points $\boldsymbol{z}\in \Psi^{*}$ are on the boundary of the convex envelope of $\boldsymbol{Z}$, which is the smallest convex set that contains $\boldsymbol{Z}$. 
\end{definition}
Note this does not imply that the Pareto front is continuous, but rather it is a relaxed interpretation of \textit{convexity} introduced to distinguish between a supported and an unsupported nondominated objective vector. Only supported nondominated points are attainable by means of a weighted-sum scalarization approach.

\section{Background and Motivation}\label{sec:Motivation}
We start with a short recapitulation of the two MOO algorithms that motivate the new algorithm proposed in the current work.

\subsection{Sandwich Algorithm}\label{section2-sandwich}
The SD algorithm proposed by \cite{Solanki1993} is a scalarization method developed aimed at approximating a (convex) Pareto front while solving as few optimization subproblems as possible. While there are variations of the approach, SD algorithms are based on the successive solution of weighted-sum subproblems in which the multiple objective functions are aggregated into a single objective function using a weight vector, $\boldsymbol{w}$, as follows:
\begin{equation}\tag{WSP$w$}\label{eq:ch6:weightedsum}
\begin{aligned} 
&\underset{\boldsymbol{x}\in X}{\operatorname{min}} \quad \boldsymbol{w}^{\top} \boldsymbol{f}(\boldsymbol{x})=\sum_{i=1}^{m} w_{i} f_{i}(\boldsymbol{x}). 
\end{aligned}
\end{equation}
The key concepts necessary to describe an SD algorithm are set out in this section by adopting the definitions and notation of \cite{Rennen2011}.  A detailed description of the algorithm, including figures and a pseudo-algorithm, is provided in the Supplementary Materials (see Section S1 therein).

\begin{definition}  
\textit{Anchor point} -- Anchor point $\boldsymbol{z}^{Ai}$ is a point in objective space equal to the optimal objective function vector obtained by solving problem \eqref{eq:ch6:weightedsum} with weight vector $\boldsymbol{w}^i$ such that $w^i_i = 1-\delta$ and $w^i_j=\delta$, $j=1, \ldots, m, \   j \neq i$, where $\delta$ is a small positive infinitesimal scalar value. Thus $\boldsymbol{z}^{Ai} = \boldsymbol{f}(\boldsymbol{x}^*)$, where $\boldsymbol{x}^* = \argmin\limits_{\boldsymbol{x}\in X} \boldsymbol{w}^{i\top}\boldsymbol{f}(\boldsymbol{x})$.  
\end{definition}

\begin{definition}  
\textit{Hyperplane} -- A \textit{hyperplane in the objective space} is given by $H(\boldsymbol{w}, b)=\left\{\boldsymbol{z} \in Z \mid {\boldsymbol{w^{\top}} \boldsymbol{z}=b}\right\}$ with $\boldsymbol{w} \in \mathbb{R}^{m} \backslash\{\boldsymbol{0}\}$ and $b \in \mathbb{R}$. 
\end{definition}\label{df:hyperplane}

\begin{remark}The vector $\boldsymbol{w}$ is a normal to the hyperplane. In this paper, the vector $\boldsymbol{w}$ is always taken to be a unit normal vector, such that $\|\boldsymbol{w}\|=1$. 
\end{remark}
A half-space and the inner normal are defined as follows.

\begin{definition} 
\textit{Half-space} --
The set $HS(\boldsymbol{w}, b)=\left\{\boldsymbol{z} \in \mathbb{R}^m \mid \boldsymbol{w}^{\top} \boldsymbol{z} \geq b\right\}$ is the \textit{half-space}. 
\end{definition}

\begin{remark}The vector $\boldsymbol{w}$, $\|\boldsymbol{w}\|=1$, is an \textit{inner unit normal} of the half-space. 
\end{remark}

\begin{definition} 
\textit{Supporting hyperplane of a convex set} -- Suppose $C$ is a convex set and $\boldsymbol{z}_0$ is a point that lies on the boundary of $C$, $\boldsymbol{bd}\ C$, i.e., $\boldsymbol{z}_0 \in \boldsymbol{bd} \ C$. If $\boldsymbol{w}\neq \boldsymbol{0}$ satisfies $\boldsymbol{w^{\top}z} \geq \boldsymbol{w^{\top}z}_0$ for all $\boldsymbol{z}\in C$, then a hyperplane $H(\boldsymbol{w},b)$ defined by $\boldsymbol{w^{\top}z}=b$, where $b=\boldsymbol{w^{\top}z}_0$, is called a \textit{supporting hyperplane} to $C$ at point $\boldsymbol{z}_0$. This is equivalent to saying that the hyperplane $H(\boldsymbol{w},b)$ supports $C$ at $\boldsymbol{z}_0$. 
\end{definition}

\begin{definition} \textit{$k$-face, facets and extreme points} -- 
A set of points $F$ of dimensionality $k$ is a $k$-face of $C$ if there exists a hyperplane $H(\boldsymbol{w},b)$ that supports $C$ at $\boldsymbol{z}_0$ and for which it holds that $H(\boldsymbol{w}, b) \cap C=F .$ If $C \subseteq \mathbb{R}^{m}$, its $(m-1)$-faces are the facets $F_{S}$ and its $0$-faces are the extreme points (vertices).
\end{definition}

For $m=2$, we denote a $p$-th facet which is defined with two extreme points $\boldsymbol{z}^{p1}$ and $\boldsymbol{z}^{p2}$ as $F^{p}_{S}(\boldsymbol{z}^{p1},\boldsymbol{z}^{p2})$, where it holds $z^{p1}_1<z^{p2}_1$. A graphical interpretation of supporting hyperplanes, half-spaces, $(m-1)$-faces, and 0-faces in two-dimensional space is given in Figure S1 in the Supplementary Materials.

\begin{definition} 
\textit{Inner and outer approximations} -- A set $IPS \subseteq Z$  is an inner approximation of $\Psi^*$ if it satisfies:
\[
\forall \boldsymbol{z}' \in IPS, \exists \boldsymbol{z} \in \Psi^*: \boldsymbol{z}'\geq \boldsymbol{z}, 
\]
where the inequality is understood component-wise.
Similarly, a set $OPS \subseteq \mathbb{R}^m$ is an outer approximation of $\Psi^*$ if it satisﬁes:
\[
\forall \boldsymbol{z}' \in OPS, \exists \boldsymbol{z} \in \Psi^*: \boldsymbol{z}'\leq \boldsymbol{z}.
\]
\end{definition}

The above definition does not specify how the inner and outer approximation can be constructed. The SD algorithm is based on iteratively refining inner and outer approximations of the Pareto front. An inner approximation is formed via a convex hull of the nondominated points identified up to the current iteration,

\begin{equation}
IPS = \operatorname{convhull}(\boldsymbol{Z_E})
= \left\{ 
\sum_{k=1}^{N_E}\lambda_k \boldsymbol{z}^k \;\middle|\;
\lambda_k \ge 0,\; \sum_{k=1}^{N_E}\lambda_k = 1
\right\},
\end{equation}
where \(Z_E=\{z^1,\ldots,z^{N_E}\} \subseteq \Psi^{*}\) denotes the current set of nondominated points and $N_E=|Z_E|$. An outer approximation is constructed from the supporting half-spaces obtained through successive solutions  of weighted-sum subproblems \eqref{eq:ch6:weightedsum} and it is defined by multiple Pareto points and their associated weight coefficients:
\begin{equation}\label{eq:si:OPS}
\begin{aligned}
OPS= \{ \boldsymbol{z} \ | \ \boldsymbol{W}^{\top}\boldsymbol{z} \geq \boldsymbol{B}\},\\
\end{aligned}
\end{equation}
where $\boldsymbol{W}$ is a matrix with column $k$ corresponding to a particular weight vector $\boldsymbol{w}^k$ and $\boldsymbol{B}$ is a column vector where the $k$-th element is obtained by the vector product
between the transpose of the $k$-th column of $\boldsymbol{W}$ and the corresponding nondominated solution $\boldsymbol{z}^{k}\in \boldsymbol{Z_E}$, i.e., $\boldsymbol{B}=\boldsymbol{W}^{\top}\boldsymbol{Z_E}$.

These approximations progressively sandwich the Pareto front, improving with each iteration. The process continues until the maximum approximation error, $d_{max}$, which measures the distance between the inner and outer approximations, is less than a predefined tolerance $\epsilon$. The error $d_{max}$ is calculated as the largest facet-specific error, $d_{error,p}$, as defined by Equations (S2) and (S3). The major advantage of this method is that the weight vectors are systematically chosen to be normal to the facets generated by the inner approximation and the accuracy of the approximation is measured during the course of the algorithm, providing a natural stopping criterion.

\subsection{Normal Boundary Intersection}
The original NBI method was proposed by \cite{DASDENNIS1998} to generate uniformly distributed boundary points for a nonlinear multiobjective optimization problem, as illustrated in Figure S3(a). To further ensure Pareto-optimality of all solutions, \cite{shukla2007normal} modified the method using a goal-attainment approach. In the modified NBI (mNBI) method, the individual minima or anchor points are found as a first step. The CHIM is then generated as the set of all convex combinations of these individual minima. The CHIM can be expressed as $\{\boldsymbol{\Phi\beta}: {\boldsymbol{\beta}} \in \mathbb{R}^{m}_{+}, \sum_{j=1}^{m} \beta_{j}=1 \}$ where $\boldsymbol{\Phi} \in \mathbb{R}^{m \times m}$ is a matrix with $i^{\mathrm{th}}$ column $\boldsymbol{z}^{Ai}-\boldsymbol{f}^{id}$. 
The mNBI subproblem (\ref{eq:ch6:mNBI}) is formulated to search for the maximum distance $t$ along a normal vector $\boldsymbol{\bar{n}}$ to the CHIM at a point defined by a specific choice of $\boldsymbol{\beta}$. The resulting subproblem is described as:
\begin{equation}\tag{mNBI$\beta$}\label{eq:ch6:mNBI}
\begin{aligned}
&\begin{array}{c}
\max\limits_{\boldsymbol{x}\in X,t} \ t \\
\text { s.t. } \boldsymbol{\Phi \beta}+t\boldsymbol{\bar{n}}\geq\boldsymbol{f}(\boldsymbol{x})-
\boldsymbol{f}^{id}\\
t \in \mathbb{R}, \ \boldsymbol{\bar{n}} \in \mathbb{R}^{m}
\end{array}
\end{aligned}
\end{equation}

A comprehensive description of both the NBI and mNBI methods is provided in Section S2 of the Supplementary Materials, including the basic scheme of the mNBI algorithm (Algorithm S2) and a schematic illustration of the bi-objective optimization procedure (Figure S3).

\subsection{Discussion}
The SD algorithm has been successfully applied to several practical problems, such as the intensity-modulated radiation therapy optimization \citep{craft2006approx,Rennen2011}, the design of distillation processes for separating a binary mixture of chloroform and acetone  \citep{bortz2014multi}, the integrated design of solvents and processes \citep{burger2015hierarchical}, the design of a plate-fin heat-sink to improve economic and thermal performance \citep{andersson2018moo} and the solution of parameter estimation problems for an equation-of-state model \citep{graham2022multi}, yet there are remaining difficulties to overcome. Since the SD algorithm is closely related to weighted-sum methods, unsupported solutions that lie in the nonconvex regions of the objective space cannot be attained. Thus, there is no weight vector $\boldsymbol{w}$ suitable for finding an efficient solution $\boldsymbol{x}^{*}$ in such a region. 

On the other hand, the mNBI method can yield solutions that lie on the nonconvex (concave) parts of the Pareto front as well as points on the convex regions. These nondominated points can be generated with the appropriate choice of $\boldsymbol{\beta}$. The efficiency of the NBI method and its variants has been investigated in the context of the optimal design of a methyl ethyl ketone production process by carrying out bi-objective optimization with an economic value and potential environmental impact as the two objectives \citep{lim2001efficient}, the development of optimal bidding strategies for the participants of oligopolistic energy markets \citep{vahidinasab2010normal}, the design of an optimal heat exchanger that maximizes the amount of heat transfer and minimizes tube length \citep{siddiqui2012improving}, and the solution of a scheduling problem for pumped hydrothermal power in which a cost function and an emission metric were considered as bi-objective criteria \citep{simab2018multi}. However, unlike the SD algorithm, the mNBI method does not provide information on approximation accuracy, nor the maximum achievable approximation quality level in the presence of disconnected regions, thus making it difficult to determine when to stop the algorithm and how to choose the number of iterations $N_{\beta}$. These drawbacks also preclude the user from choosing the set of reference points $\boldsymbol{\Phi\beta}^{k}\in C^{ref}$ systematically. 
Given the advantages derived from each algorithm and their similarity in terms of the use of convex hulls, an extension of the SD algorithm through its combination with the mNBI method is discussed in the subsequent sections.

\section{Proposed Algorithm}\label{sec:SDNBI}
We present here the SDNBI algorithm, a novel adaptive MOO approach for the solution of BOO problems. It is more general than the SD algorithm in that it is applicable to nonconvex problems, in particular to those with a disconnected Pareto front and feasible region. A new feature is the use of the mNBI subproblem, such that the mNBI parameters are systematically determined by the SD algorithm, in an adaptive fashion. Although the SD algorithm and the mNBI method are applicable to higher-dimensional MOO cases ($m>2$) \citep{bokrantz2013vertex}, we focus here on the solution of two-dimensional (BOO) problems and adapt the terminology appropriately, as described in the following definition.

\begin{definition}
\textit{Line, half-plane, tangent} -- A \textit{line} $H_2(\boldsymbol{w},b)$ and a \textit{half-plane} $HS_2(\boldsymbol{w},b)$ are special cases of the hyperplane $H(\boldsymbol{w},b)$ and half-space $HS(\boldsymbol{w},b)$, respectively, defined for the two-dimensional objective space, $m=2$. If a line is tangential to a curve in the bi-objective space at the point $\boldsymbol{z}_0$, it is called a \textit{tangent} at $\boldsymbol{z}_0$.
\end{definition}

\subsection{Scalarization scheme: Setting the parameters of the mNBI subproblem}
The main idea behind the proposed SDNBI algorithm is to replace the weighted-sum subproblem that appears in the SD algorithm with the mNBI subproblem \eqref{eq:ch6:mNBI} so as to obtain solutions on the nonconvex part(s) of the Pareto front.
Consider a facet $F^{p^k}_{S}(\boldsymbol{z}^{p1},  \boldsymbol{z}^{p2})$ connecting two nondominated points $\boldsymbol{z}^{p1}$ and $\boldsymbol{z}^{p2}$ at some iteration $k$. The normal vector of the facet can be used as the normal vector $\boldsymbol{\bar{n}}^{k}$ in the mNBI subproblem. Then, $\boldsymbol{\Phi}^{k}$ is defined by the two extreme points of the facet, with its $i^{\rm th}$ column given by $\boldsymbol{z}^{pi}-\boldsymbol{z}^{o}$, $i=1,2$, where $\boldsymbol{z}^{o}=[\operatorname{min}(z^{p1}_1,z^{p2}_1), \operatorname{min}(z^{p1}_2,z^{p2}_2)]^{\top}=[z^{p1}_1,z^{p2}_2]^{\top}$. Subsequently, the reference point $\boldsymbol{\Phi}^k \boldsymbol{\beta}^k$ on the facet can be determined by choosing an appropriate $\boldsymbol{\beta}^k$. In this study,  $\boldsymbol{\beta}^k$ is chosen such that the reference point $\boldsymbol{\Phi}^k \boldsymbol{\beta}^k$ is located at the midpoint of the facet, in order to generate a Pareto point that may be evenly placed between $\boldsymbol{z}^{p1}$ and $\boldsymbol{z}^{p2}$. 
While in the mNBI algorithm, $\boldsymbol{\bar{n}}$ remains fixed throughout and the $\boldsymbol{\beta}$ values are chosen \textit{a priori}, in the SDNBI algorithm, they are updated at each iteration.  
%
\subsection{Constructing valid inner and outer approximations}
The solution of the mNBI subproblem at some iteration $k$ with the proposed parameter values results in point $\boldsymbol{z}^k$ that is used to update the approximations of the Pareto front. This requires careful analysis due to the potential nonconvexity of the efficient frontier. In the original SD algorithm, the outer approximation $OPS$ is improved at iteration $k$ by adding a line tangential to the recently identified nondominated point $\boldsymbol{z}^k$, $H_2(\boldsymbol{w}^{k},\boldsymbol{b}^{k})$, and its positive half-plane $HS_2(\boldsymbol{w}^{k},b^{k})$ where $b^{k}=\boldsymbol{w}^{k \top}\boldsymbol{z}^k$. The weight vector $\boldsymbol{w}^k$ is used in the weighted-sum scalarization problem as described in Section \ref{section2-sandwich}. The inner approximation $IPS$ of the Pareto front is obtained as a polyhedral approximation of the set of points $\boldsymbol{z}\in \boldsymbol{Z_E} \subseteq \Psi^*$. This polyhedral approximation is defined as the set of all convex combinations of points in $\boldsymbol{Z_E}$ in the objective space and it is denoted by convexhull($\boldsymbol{Z_E}$). 

\begin{figure}[t]
	\centering
	\includegraphics[height=0.30\textheight]{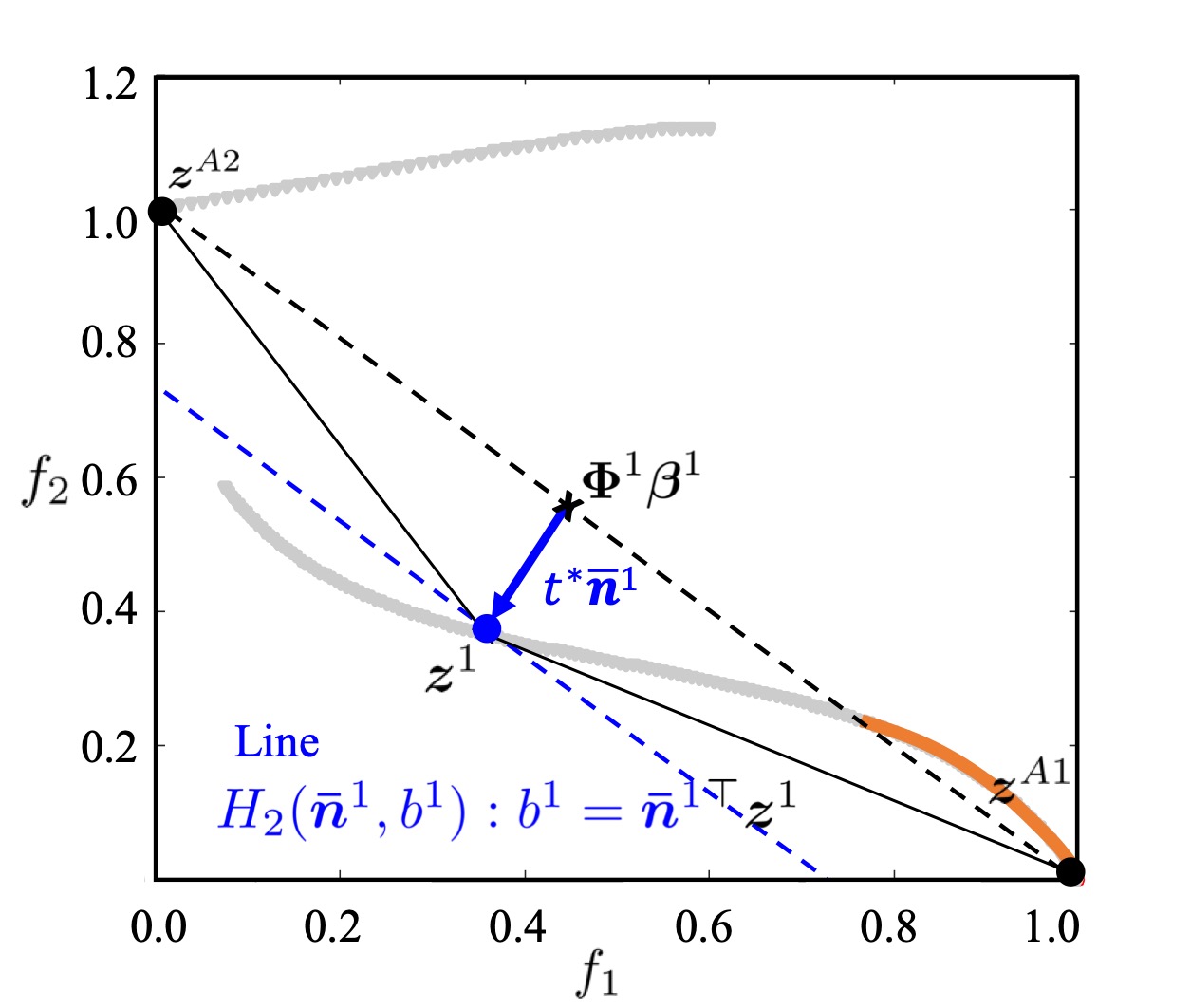}
	\vspace{-0.4cm}
	\captionsetup{font=small}\caption{Illustration of a nonconvex Pareto front and the generation of an inner approximation that cuts off part of the front. The grey curve represents the boundary of the feasible region in bi-objective space for problem MOP1. The black dashed line is a facet obtained by generating a convex hull using $\boldsymbol{z}^{A1}$ and $\boldsymbol{z}^{A2}$. The nondominated point is $\boldsymbol{z}^1$ generated by the solution of \eqref{eq:ch6:mNBI} with the normal vector $\boldsymbol{\bar{n}}^1$ and the reference point $\boldsymbol{\Phi}^{1}\boldsymbol{\beta}^{1}$. The corresponding line $H_2(\boldsymbol{\bar{n}}^1, b^1)$ (blue dashed line) is not tangent to the Pareto front. If the next upper bounds, i.e., inner approximations, are obtained by constructing the convex hull based on the updated set of the Pareto points, $\boldsymbol{Z_E}=\{\boldsymbol{z}^{A2}, \boldsymbol{z}^{1}, \boldsymbol{z}^{A1}\}$, the solid black lines are obtained, cutting off the part of the Pareto front shown in orange.}\label{fig:ch6:w-mnbi}
\end{figure}

Thus, in the original SD algorithm, the Pareto front must be convex to ensure that valid upper and lower bounds are obtained that do not cut off any part of the front. If the Pareto front is nonconvex, the line $H_2(\boldsymbol{\bar{n}}^{k},b^{k})$ generated as a solution of problem \eqref{eq:ch6:mNBI} does not necessarily support the Pareto front at point $\boldsymbol{z}^k$. A simple example where the use of the mNBI subproblem to generate a new Pareto point leads to erroneous inner and outer approximations is shown in Figure \ref{fig:ch6:w-mnbi}. It is clear that the Pareto frontier is not supported by the line $H_2(\boldsymbol{\bar{n}}^1,b^1)$, where ${b^1}=\boldsymbol{\bar{n}}^1 \top
\boldsymbol{z}^1$, at nondominated point $\boldsymbol{z}^1$. The inner approximation generated by facets $F^{p^1}_{S}(\boldsymbol{z}^{A2},\boldsymbol{z}^1)$ and $F^{p^1}_{S}(\boldsymbol{z}^1,\boldsymbol{z}^{A1})$ is valid only for Pareto points $\boldsymbol{z}^{\top}=(z_1,z_2)$ in objective space such that $z_1 \leq z^1_1$. We address these issues in the remainder of this subsection.  

\subsubsection{Tangent at the solution of mNBI}
To analyze further the properties of the mNBI subproblem, it is useful to partition the vector of variables into a continuous variable vector $\boldsymbol{x} \in \mathbb{R}^{n_1}$ and a binary variable vector, denoted as $\boldsymbol{y} \in Y= \{0,1\}^{n_2}$. Without loss of generality, any function $q(\boldsymbol{x}, \boldsymbol{y})$ in the MOO problem, where $q$ can refer to the objective function, an equality or an inequality constraint, can be written as \citep{floudas1995} 
\begin{equation}
    q(\boldsymbol{x}, \boldsymbol{y}) = c_q^{\top}\boldsymbol{y} + q_x(\boldsymbol{x}).
\end{equation}
Then the mNBI subproblem is given by
\begin{equation}\label{eq:ch6:mNBI-MINLP}
\begin{aligned}
&\begin{array}{c}
\max\limits_{\boldsymbol{x}\in X,\boldsymbol{y}\in Y,t} \ t \\
\text { s.t. } \boldsymbol{\Phi \beta}+t\boldsymbol{\bar{n}}\geq\boldsymbol{f}_x(\boldsymbol{x}) + \boldsymbol{c}^{\top}_f \boldsymbol{y}-
\boldsymbol{f}^{id}\\
t \in \mathbb{R}, \ \boldsymbol{\bar{n}} \in \mathbb{R}^{m}
\end{array}
\end{aligned}
\end{equation}

At the solution ($\hat{\boldsymbol{x}}$, $\hat{\boldsymbol{y}}$, $t^*$) of the MINLP problem, one can obtain Lagrange multipliers at fixed  $\boldsymbol{y}=\hat{\boldsymbol{y}}$, provided that the corresponding NLP meets an appropriate constraint satisfaction. Therefore, given an efficient solution ($\hat{\boldsymbol{x}}$,$\hat{\boldsymbol{y}}$) obtained by solving an mNBI subproblem and at which $t=t^*$, there exist $\boldsymbol{\mu}^{*}$ and $\boldsymbol{\nu}^{*}$ such that the KKT optimality conditions for optimization problem \eqref{eq:ch6:mNBI-MINLP} at fixed $\hat{\boldsymbol{y}}$ are satisfied so that:
\begin{eqnarray}
\nabla_{\boldsymbol{x}}\mathcal{L}(\hat{\boldsymbol{x}},t^*,\boldsymbol{\mu}^{*},\boldsymbol{\nu}^{*};\hat{\boldsymbol{y}}) &=& \boldsymbol{\mu}^{*\top} {\nabla_{\boldsymbol{x}} \boldsymbol{f}_x\left(\hat{\boldsymbol{x}}\right)}+ \boldsymbol{\nu}^{*\top} {\nabla_{\boldsymbol{x}} \hat{\boldsymbol{h}}_x\left(\hat{\boldsymbol{x}}\right)} =\boldsymbol{0} ,\label{eq:ch6:NBIKKT1}\\
{\nabla_{t} \mathcal{L}}(\hat{\boldsymbol{x}},t^*,\boldsymbol{\mu}^{*},\boldsymbol{\nu}^{*};\hat{\boldsymbol{y}}) &=& -1 + \boldsymbol{\mu}^{*\top} \boldsymbol{\bar{n}} = 0, \label{eq:ch6:NBIKKT2}
\end{eqnarray}
where $\mathcal{L}(\boldsymbol{x}, t, \boldsymbol{\mu}, \boldsymbol{\nu};\hat{\boldsymbol{y}}) 
=\ -t+ \boldsymbol{\mu}^{\top}\left(\boldsymbol{f}_x(\boldsymbol{x})+\boldsymbol{c}_f^{\top} \hat{\boldsymbol{y}}-\boldsymbol{f}^{id}-\boldsymbol{\Phi\beta}-t \boldsymbol{\bar{n}}\right)+ \boldsymbol{\nu}^{\top}(\hat{\boldsymbol{h}}_x(\boldsymbol{x})+\boldsymbol{c}_h^{\top} \hat{\boldsymbol{y}})$, $\boldsymbol{\mu}\in \mathbb{R}^m$ represents the vector of Lagrange multipliers corresponding to the augmented objective constraints $\boldsymbol{f}_x(\boldsymbol{x})+\boldsymbol{c}_f^{\top} \hat{\boldsymbol{y}}-\boldsymbol{f}^{id}-\boldsymbol{\Phi\beta}-t \boldsymbol{\bar{n}} \leq \boldsymbol{0}$, and $\boldsymbol{\nu} \in \mathbb{R}^s$ is the vector of Lagrange multipliers for the $s$ active constraints in the set $\{\boldsymbol{g}_x({\boldsymbol{x}})+\boldsymbol{c}_g^{\top} \hat{\boldsymbol{y}}\leq \boldsymbol{0},\boldsymbol{h}_x({\boldsymbol{x}}) + \boldsymbol{c}_h^{\top} \hat{\boldsymbol{y}}= \boldsymbol{0}\}$, represented by the vector $\boldsymbol{\hat{h}}_x(\boldsymbol{x})\in \mathbb{R}^{s}$, $q \leq s\leq p+q$. Furthermore, at the KKT point, the following complementarity condition must hold:
\begin{equation}\label{eq:ch6:NBIKKT-compl}
\begin{array}{l}
\boldsymbol{\mu}^{*\top}\left(\boldsymbol{f}_x(\hat{\boldsymbol{x}})+ \boldsymbol{c}_f^{\top} \hat{\boldsymbol{y}}-\boldsymbol{f}^{id}-\boldsymbol{\Phi\beta}-t^* \boldsymbol{\bar{n}}\right)=0.
\end{array}
\end{equation}
From Equations \eqref{eq:ch6:NBIKKT2} and \eqref{eq:ch6:NBIKKT-compl}, it can be deduced that at least one of the augmented objective constraints must be active. This can be shown by contradiction. Let us assume that all constraints in the augmented objective constraints are inactive, i.e., $\boldsymbol{f}_x(\hat{\boldsymbol{x}})+ \boldsymbol{c}_f^{\top} \hat{\boldsymbol{y}}-\boldsymbol{f}^{id}-\boldsymbol{\Phi\beta}-t \boldsymbol{\bar{n}} < \boldsymbol{0}$ and $\boldsymbol{\mu}^{*}=\boldsymbol{0}$. Then, Equation \eqref{eq:ch6:NBIKKT2} is reduced to ${\nabla_{t} \mathcal{L}} = -1 \neq 0$, which is a violation of the KKT necessary conditions. Therefore, if we solve problem \eqref{eq:ch6:mNBI} for any choice of parameters $(\boldsymbol{\bar{n},\beta})$, there exists a corresponding normal vector $\boldsymbol{w}'\geq \boldsymbol{0}$, where the inequality is understood component-wise,  that defines a tangent to the Pareto front at the nondominated solution, $\boldsymbol{z}$, and is given by:
\begin{equation}\label{eq:ch6:normal-mNBI}
\begin{array}{l}
 \boldsymbol{w}'=\dfrac{ 1 }{ \sum_{i=1}^{m} \mu_{i}^{*} }\boldsymbol{\mu}^{*}, \  \sum_{i=1}^{m} w'_{i}=1. 
\end{array}
\end{equation}

As \citet{DASDENNIS1998} and \citet{shukla2007normal} have shown, an even stronger relationship between the NBI/mNBI and weighted-sum subproblems exists when the solution of the NBI/mNBI subproblem is located on a convex part of the Pareto front. A detailed discussion of this relationship is provided in Section S3 of the Supplementary Materials.

\subsubsection{Decomposition of the search space}\label{sec-decomposition}
\begin{figure}[t]
	\centering
	\includegraphics[height=0.30\textheight]{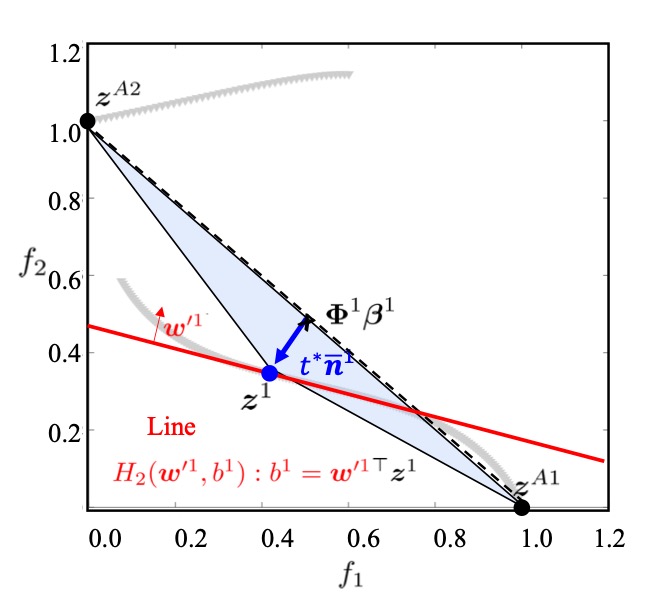}
	\vspace{-0.4cm}
	\captionsetup{font=small}\caption{Representation of the line $H_2(\boldsymbol{w}'^1,b^1)$ passing through $\boldsymbol{z}^1$ with normal vector $\boldsymbol{w}'^1$ (red solid line), calculated based on Equation \eqref{eq:ch6:normal-mNBI}. The grey curve represents the boundary of the feasible region in bi-objective space and the black solid lines are the facets obtained by constructing a convex hull based on the current Pareto points $\boldsymbol{z}^{A1}$, $\boldsymbol{z}^{A2}$ and $\boldsymbol{z}^{1}$. Facet  $(\boldsymbol{z}^{1},\boldsymbol{z}^{A1})$ and line $H_2(\boldsymbol{w}'^1,b^1)$ cannot be used to construct approximations of the Pareto front due to  nonconvexity.}\label{fig:ch6:w-mnbi-tangent}
\end{figure}

Given the line $H_2(\boldsymbol{w}'^1,b^1)$ tangential to the Pareto front at $\boldsymbol{z}^1$, the next step is to construct the outer and inner approximations. However, one can observe from Figure \ref{fig:ch6:w-mnbi-tangent} that the tangent at $\boldsymbol{z}^1$, $H_2(\boldsymbol{w}'^{1},b^1)$, intersects the nonconvex part of the Pareto front and cannot be used to generate a valid outer approximation. Similarly, the facet $(\boldsymbol{z}^1,\boldsymbol{z}^{A1})$ is below the Pareto front and cannot be used as an inner approximation. To provide valid inner and outer approximations of the Pareto front, we introduce a systematic decomposition of the search region based on the convexity of the Pareto front at the Pareto points identified, such that each Pareto-optimal point $\boldsymbol{z}^k$ is supported by a line $H_2(\boldsymbol{w}'^k,b^k)$ in a given space and the collection of supporting lines represents the upper or lower bounds of the Pareto front. 

The idea underlying the decomposition of the search space is to divide the objective space into a number of sub-regions whenever there exist Pareto points that are exterior to the outer approximation as defined in the SD algorithm, i.e., the polyhedral set given by the supporting lines that have been derived so far. 
\begin{definition} \textit{\Subspace} --
An $l^{\rm th}$ \subspace $C^{l}(\boldsymbol{z}^{i},\boldsymbol{z}^{j})$ is a subset of bi-objective space defined by the two nondominated points $\boldsymbol{z}^{i}$ and $\boldsymbol{z}^{j}$, where  $z^i_1 < z^j_1$, without loss of generality. It consists of the rectangle whose diagonal is the line segment $[\boldsymbol{z}^{i},\boldsymbol{z}^{j}]$ such that all Pareto points $\boldsymbol{z}=(z_1,z_2)^{\top}\in \Psi^{*}$ in the interior of the search space satisfy $z^{i}_1
\leq z_1 \leq z^{j}_1$ and $z^{j}_2
\leq z_2 \leq z^{i}_2$.
\end{definition}

\begin{definition} \textit{Extreme points of a \subspace} -- Points $\boldsymbol{z}^{i}$ and $\boldsymbol{z}^{j}$ in \subspace $C^{l}(\boldsymbol{z}^{i},\boldsymbol{z}^{j})$ are referred to as extreme points of the \subspace. \end{definition}

\begin{figure}
	\centering
	\includegraphics[height=0.23\textheight]{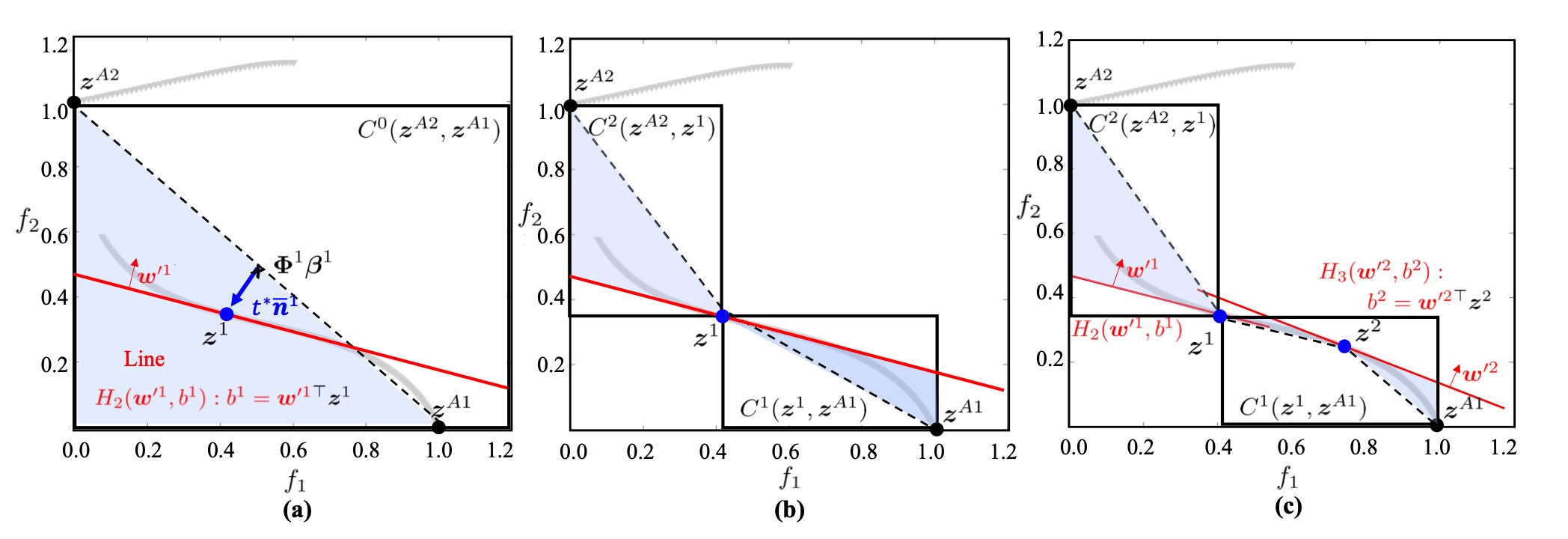}
	\vspace{-0.3cm}
	\captionsetup{font=small}\caption{Schematic illustrating the procedure for the decomposition of the bi-objective space in the case of a nonconvex Pareto front, with boundary points of the feasible region shown with grey symbols (\textcolor{gray}{$\bullet$}): (a) An initial objective space is characterized by \subspace $C^0(\boldsymbol{z}^{A2},\boldsymbol{z}^{A1}$) (black rectangular box), assumed to contain a convex Pareto front. A new nondominated point $\boldsymbol{z}^1$ is found by solving problem \eqref{eq:ch6:mNBI} for $\boldsymbol{\Phi}^{1}\boldsymbol{\beta}^{1}$; (b) the initial objective space is decomposed into \subspaces  $C^1(\boldsymbol{z}^{1},\boldsymbol{z}^{A1})$ and $C^2(\boldsymbol{z}^{A2},\boldsymbol{z}^{1})$ that are assumed to contain nonconvex and convex Pareto fronts, respectively. The inner approximation of the Pareto front in  $C^2$ is the facet $(\boldsymbol{z}^{A2},\boldsymbol{z}^{1})$ (black dashed line) and the outer approximation is given by the segment of the tangent at $\boldsymbol{z}^{1}$ (red solid line) connecting $\boldsymbol{z}^{1}$ to the $f_1=0$ line and by the $f_1=0$ line. For \subspace $C^1$, the inner approximation is given by the segment of the tangent at $\boldsymbol{z}^{1}$ (red solid line) connecting $\boldsymbol{z}^{1}$ to the $f_1=1$ line and by the $f_1=1$ line, while the outer approximation is given by the facet $(\boldsymbol{z}^{1},\boldsymbol{z}^{A1})$;  (c) the approximation of the Pareto front is improved by adding $\boldsymbol{z}^2$ and the \subspace $C^1$ is not decomposed since $\boldsymbol{w}'^{2
 \top}\boldsymbol{z} \leq b^{2}$ holds for all $\boldsymbol{z}\in \boldsymbol{Z}^{C^1}_{\boldsymbol{E}}$. For \subspace $C^1$, the inner approximation is given by the segment of the tangent at $\boldsymbol{z}^{2}$ (red solid line), while the outer approximation is given by the convexhull $(\boldsymbol{z}^{1},\boldsymbol{z}^{2},\boldsymbol{z}^{A1})$}\label{fig:ch6:SpaceDecompo}
\end{figure}

Let $C^{0}(\boldsymbol{z}^{A2},\boldsymbol{z}^{A1})$ denote an initial search space defined by the anchor points $\boldsymbol{z}^{A1}$ and $\boldsymbol{z}^{A2}$, as can be seen in Figure \ref{fig:ch6:SpaceDecompo}. The initial set of Pareto points is $Z_E=\{\boldsymbol{z}^{A1},\boldsymbol{z}^{A2}\}$.
It is initially assumed that the Pareto front in $C^{0}$ is convex. By solving subproblem \eqref{eq:ch6:mNBI}, a new nondominated point $\boldsymbol{z}^1$ is obtained that enables one to define the corresponding supporting line $H_2(\boldsymbol{w}'^{1},b^{1})$ and its half-plane $HS_2(\boldsymbol{w}'^{1},b^{1})$. Note that in Figure \ref{fig:ch6:SpaceDecompo}, $\boldsymbol{z}^1$ is located on the concave (or nonconvex) part of the Pareto front. The current objective space is then investigated to determine if the convexity assumption can be justified based on the current set of Pareto points. For this, we use the supporting hyperplane theorem, a proof of which can be found in \cite{luenberger1997opt}:

\begin{theorem}\label{Th:Supporting} Supporting hyperplane theorem -- 
Suppose $C \subseteq \mathbb{R}^{m}$ is a nonempty convex set and $\boldsymbol{z}_0$ is a point on its boundary $\boldsymbol{bd}(C)$, i.e., $\boldsymbol{z}_0\in \boldsymbol{bd}(C)$. Then, there exists a supporting hyperplane $\left\{\boldsymbol{z} \mid \boldsymbol{w}'^{\top}\boldsymbol{z}=\boldsymbol{w}'^{\top}\boldsymbol{z}_0\right\}$ such that $\boldsymbol{w}'^{\top}\boldsymbol{z} \geq \boldsymbol{w}'^{\top}\boldsymbol{z}_0$ for all $\boldsymbol{z}\in C$, $\boldsymbol{z}\neq 0$.
\end{theorem}
\begin{corollary} \label{Co:suppHyper}
Suppose convex set $V\subseteq C\subseteq \mathbb{R}^{m}$, where $m=2$, is a polyhedron defined by a finite number of half-planes and lines $V=\left\{\boldsymbol{z} \mid \boldsymbol{w}'^{k \top}\boldsymbol{z}\geq\boldsymbol{w}'^{k \top}\boldsymbol{z}^k, \ k=1,\dots,n\right\}$. If each point $\boldsymbol{z}^k$ is supported by line $H_2(\boldsymbol{w}'^{k},b^{k})$ where $b^{k}=\boldsymbol{w}'^{k \top}\boldsymbol{z}^{k}$, then every point $\boldsymbol{z}\in V$ must satisfy $\boldsymbol{w}'^{k \top}\boldsymbol{z} \geq {b^k}$ for all $k=1,\dots,n$.
\end{corollary}

It can be observed in Figure \ref{fig:ch6:SpaceDecompo}(a) that both inner approximation and outer approximation are not valid, cutting off the Pareto front. The nondominated point $\boldsymbol{z}^{A1}$ does not satisfy $\boldsymbol{w}^{1 \top}\boldsymbol{z}^{A1} \geq b^1$, implying that $\boldsymbol{z}^{A1}$ is located on a nonconvex part of the Pareto front in the \subspace $C^0$. On the other hand, the nondominated point $\boldsymbol{z}^{A2}$ satisfies $\boldsymbol{w}^{1 \top}\boldsymbol{z}^{A2} \geq \boldsymbol{b}^1$ so that the hyperplane can be assumed to be supporting in the \subspace $C^{2}(\boldsymbol{z}^{A2},\boldsymbol{z}^{1})$. Subsequently, the search region is decomposed, as shown in Figure \ref{fig:ch6:SpaceDecompo}(b), into $C^{2}(\boldsymbol{z}^{A2},\boldsymbol{z}^{1})$, with the assumption that the Pareto front is convex in this \subspace, and $C^{1}(\boldsymbol{z}^{1},\boldsymbol{z}^{A1})$, with the assumption that the Pareto front is nonconvex in this \subspace. For any \subspaces that fall under a convexity assumption, the inner approximation and outer approximation are obtained as described in Section \ref{section2-sandwich}, following the SD algorithm. By contrast, for \subspaces that fall under a nonconvexity assumption, the inner approximation is constructed using the supporting lines of the nondominated points in the \subspace,  $IPS=\{\boldsymbol{z} \ | \ \boldsymbol{W}'^{\top}\boldsymbol{z}\geq\boldsymbol{B}' \}$, where $\boldsymbol{w}'$ is a matrix with each column corresponding to a particular normal vector $\boldsymbol{w}'^k$ and $\boldsymbol{B}'$ is a column vector where the $k$-th element is obtained by the vector product between the transpose of the $k$-th column of $\boldsymbol{w}'$ and the corresponding nondominated solution $\boldsymbol{z}^k$. The outer approximation is constructed as the convex hull of the nondominated points in the \subspace. The addition of a further Pareto point $\boldsymbol{z}^2$ in \subspace $C^1$ is shown in Figure \ref{fig:ch6:SpaceDecompo}(c) and does not require further decomposition of this \subspace. 


The decomposition strategy can be generalized as follows:
\begin{enumerate}
    \item Consider a search \subspace $C^l(\boldsymbol{z}^{i},\boldsymbol{z}^{j})$ in which the Pareto front is assumed to be convex, an associated set of known nondominated points, $\boldsymbol{Z}^{C^l}_{\boldsymbol{E}}$, and a newly identified nondominated point $\boldsymbol{z}^*$ within the \subspace, with the corresponding line $H_2(\boldsymbol{w}'^*,b^*)$, where $b^*=\boldsymbol{w}'^{* \top}\boldsymbol{z}^*$. If $\boldsymbol{w}'^{*\top}\boldsymbol{z} \ge {b^*}$ holds for all $\boldsymbol{z}\in\boldsymbol{Z}^{C^l}_{\boldsymbol{E}}$, add  $\boldsymbol{z}^*$ to $\boldsymbol{Z}^{C^l}_{\boldsymbol{E}}$ and retain the convexity assumption on  \subspace $C^l$.
    
    \item Consider a search \subspace $C^l(\boldsymbol{z}^{i},\boldsymbol{z}^{j})$ in which the Pareto front is assumed to be nonconvex, an associated set of known nondominated points, $\boldsymbol{Z}^{C^l}_{\boldsymbol{E}}$, and a newly identified nondominated point $\boldsymbol{z}^*$ within the \subspace, with corresponding line $H_2(\boldsymbol{w}'^*,b^*)$, where $b^*=\boldsymbol{w}'^{* \top}\boldsymbol{z}^*$. If $\boldsymbol{w}'^{*\top}\boldsymbol{z} \le {b^*}$ holds for all $\boldsymbol{z}\in\boldsymbol{Z}^{C^l}_{\boldsymbol{E}} $, add  $\boldsymbol{z}^*$ to $\boldsymbol{Z}^{C^l}_{\boldsymbol{E}}$ and retain the nonconvexity assumption on \subspace $C^l$.
    
    
    \item If the conditions in 1 and 2 are violated, then partition  \subspace $C^{l}$ into two or more \subspaces by investigating all nondominated points, $\boldsymbol{z} \in\boldsymbol{Z}^{C^l}_{\boldsymbol{E}}\cup \{ \boldsymbol{z}^* \}$. $N^l$ new \subspaces are created based on the supporting hyperplane theorem, such that $\bigcup\limits_{dc=1}^{N^l} C^{l,dc} = \boldsymbol{Z}^{C^l}_{\boldsymbol{E}}\cup \{ \boldsymbol{z}^* \}$ and for each $dc=1,\dots,N^l$,  $C^{l,dc}$ must be such that: 
       \begin{equation} \label{eq:ch6:subconv}
        \hat{\boldsymbol{w}}^\top\boldsymbol{z}\geq \hat{\boldsymbol{w}}^\top \hat{\boldsymbol{z}} \; \mathrm{ for \; all } \; \boldsymbol{z}, \hat{\boldsymbol{z}}  \in C^{l,dc},
    \end{equation} or 
    \begin{equation} \label{eq:ch6:subnonconv}
       \hat{\boldsymbol{w}}^\top\boldsymbol{z}\leq \hat{\boldsymbol{w}}^\top \hat{\boldsymbol{z}} \; \mathrm{ for \; all } \; \boldsymbol{z}, \hat{\boldsymbol{z}}  \in C^{l,dc}.
    \end{equation} 
where $\hat{\boldsymbol{w}}$ is the vector associated with nondominated point $\hat{\boldsymbol{z}}$.

The number of \subspaces is chosen  to be as small as possible while meeting the supporting hyperplane criteria. A convexity (respectively nonconvexity) assumption is made on the Pareto front in a \subspace $C^{l,dc}$ if Equation \eqref{eq:ch6:subconv} (respectively Equation \eqref{eq:ch6:subnonconv}) holds. 
\end{enumerate}

\subsection{Identification of disconnected regions of a Pareto front} \label{sec-fathom}
\begin{figure}[t]
	\centering
	\includegraphics[height=0.33\textheight]{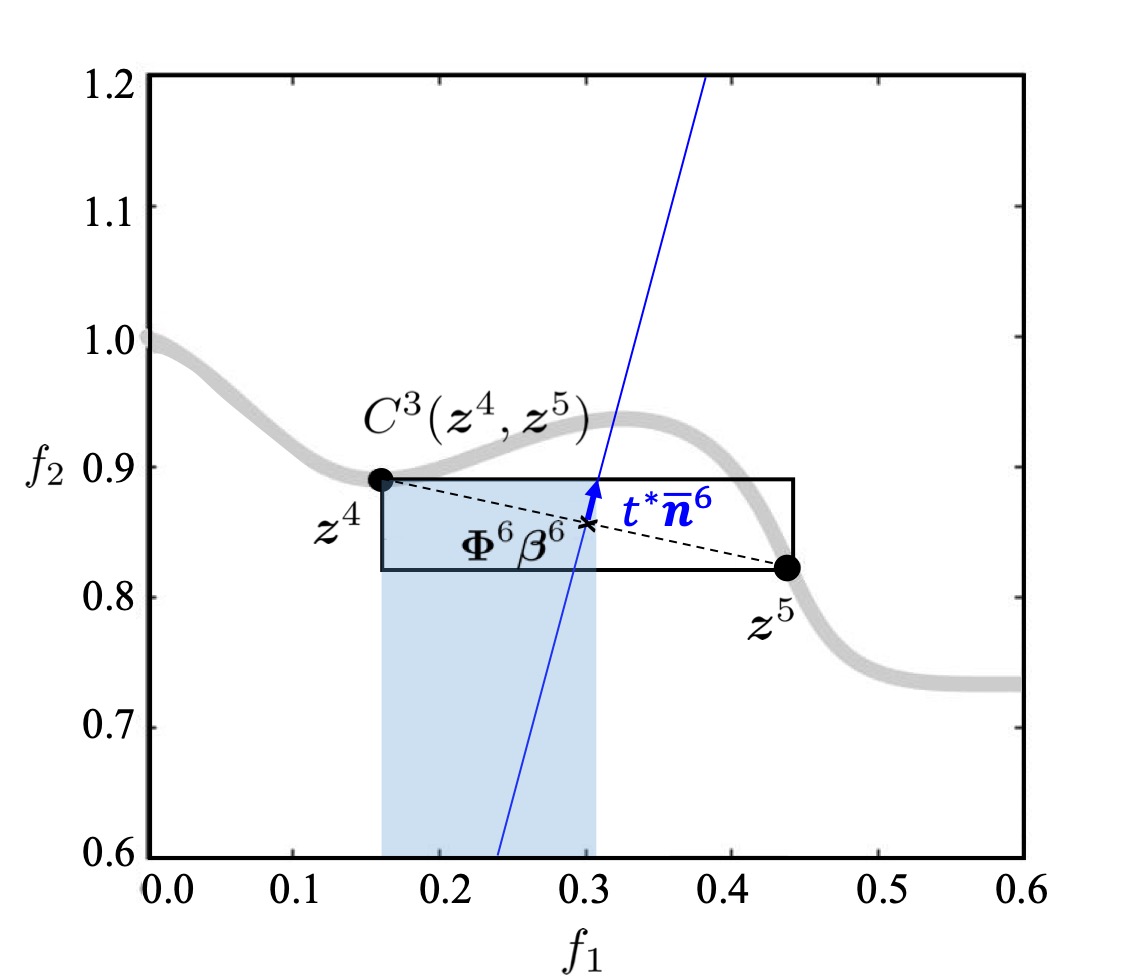}
	\vspace{-0.2cm}
	\captionsetup{font=small}\caption{A geometrical illustration of the solution of the \eqref{eq:ch6:mNBI} subproblem for a disconnected Pareto front. The boundary points of the feasible region are shown with grey symbols (\textcolor{gray}{$\bullet$}). We begin by considering \subspace $C^3(\boldsymbol{z}^4,\boldsymbol{z}^5)$ (black rectangular box). The solution of the \eqref{eq:ch6:mNBI} subproblem for facet $F_{S}(\boldsymbol{z}^4,\boldsymbol{z}^5)$ and $\boldsymbol{\beta}^6$ lies at $\boldsymbol{z}^4$. The blue shaded  region represents the area where $\boldsymbol{\Phi}^6\boldsymbol{\beta}^6 + t^*\boldsymbol{\bar{n}}^6 \ge \boldsymbol{z}^4-\boldsymbol{f}^{id}$.} \label{fig:ch6:ModifymNBI_1}
\end{figure}

\begin{figure}[t]
	\centering
	\includegraphics[height=0.3\textheight]{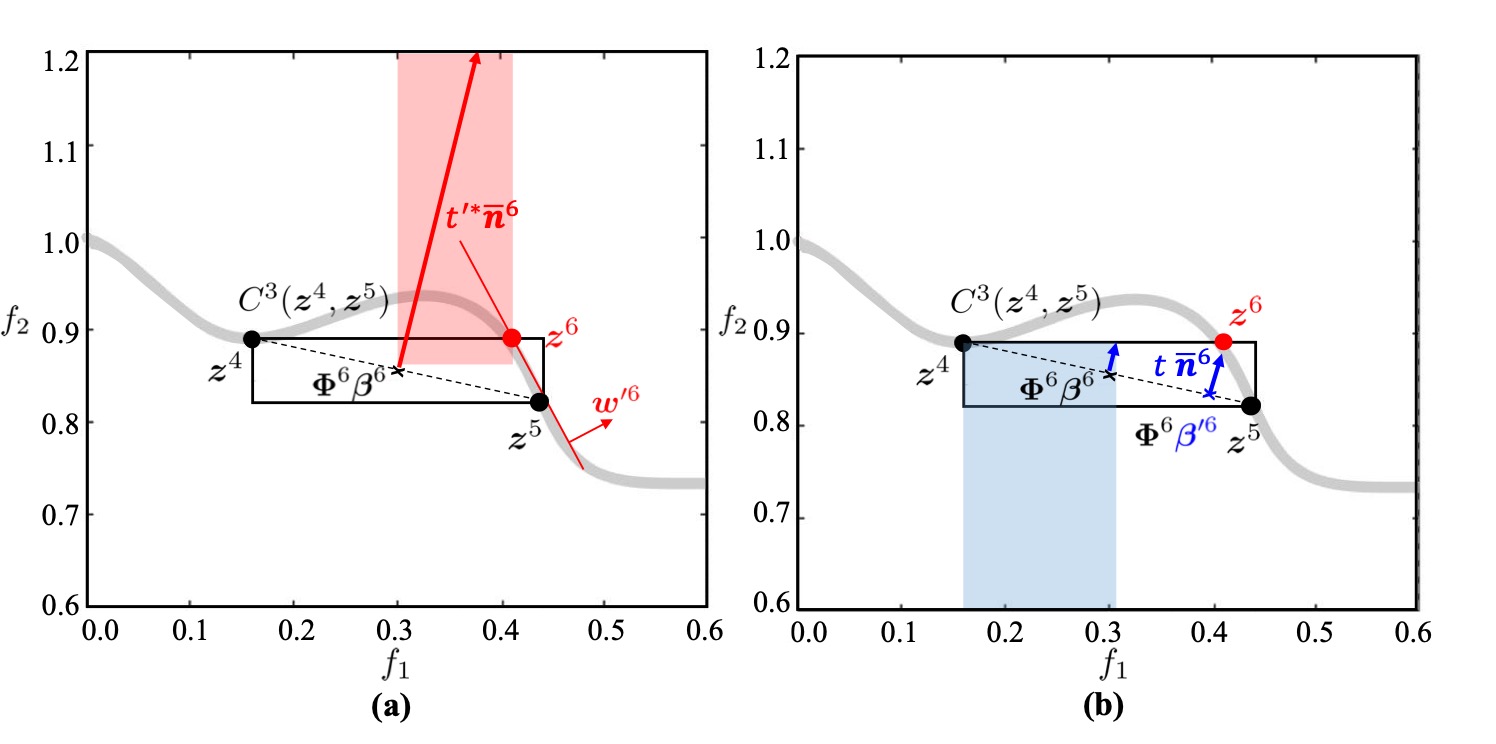}
	\vspace{-0.3cm}
	\captionsetup{font=small}\caption{A geometrical illustration of the SDNBI procedure for a disconnected Pareto front when using an alternative subproblem \eqref{eq:ch6:simpleNBIn}. The boundary points of the feasible region are shown with grey symbols (\textcolor{gray}{$\bullet$}). (a) A new Pareto point $\boldsymbol{z}^6$ is generated as a solution of \ref{eq:ch6:simpleNBIn}. The red shaded region represents the area where ($\boldsymbol{\Phi}^6\boldsymbol{\beta}^6 + t'^*\boldsymbol{\bar{n}}^6 \ge \boldsymbol{z}^6-\boldsymbol{f}^{id}) \cap (z^6_1 \geq z^4_1 + \epsilon_z$). The solid red line at $\boldsymbol{z}^6$ represents a tangent at the solution. The facet $F_{S}(\boldsymbol{z}^4,\boldsymbol{z}^6)$ generated by an inner approximation is discarded from the search space in the subsequent iterations. In (b), it can be seen that any choice of reference points on the line segment defined by $\theta\boldsymbol{z}^4 +(1-\theta) \boldsymbol{\Phi}^6\boldsymbol{\beta}^{'6}$, $\forall \ 0 
	< \theta \leq1$ (shown by the blue shaded region), produces the same solution $\boldsymbol{z}^4$. } \label{fig:ch6:ModifymNBI_2}
\end{figure}

One important aspect of the SDNBI  algorithm is its ability to increase the likelihood of finding a new Pareto point by avoiding convergence to a boundary point. However, the mNBI method may be inefficient when a Pareto front consists of many disconnected regions or in the extreme case when it consists of a finite set of points, i.e., for a purely integer problem. This is because the portion of the objective space that can be explored by the mNBI method is limited by the choice of parameters ($\boldsymbol{\bar{n}},\boldsymbol{\beta}$). If the search direction defined by these parameters does not lead to a region where as-yet-unknown nondominated points exist, the algorithm will converge to a previously identified nondominated solution, making it difficult to determine a new search direction. As a result, the algorithm may not identify Pareto points efficiently for some problems, leading to a higher computational cost. For example, a case involving the repeated computation of the same nondominated point is shown in Figure \ref{fig:ch6:ModifymNBI_1}. For a given \subspace defined by $C^{3}(\boldsymbol{z}^4, \boldsymbol{z}^5)$, the solution of subproblem \eqref{eq:ch6:mNBI} for facet normal $\boldsymbol{\bar{n}}^6$ with  reference point $\boldsymbol{\Phi}^6\boldsymbol{\beta}^6$ returns the already-known nondominated point $\boldsymbol{z}^4$. As a result, it is possible to discard the reference points belonging to the line segment defined by $\theta\boldsymbol{z}^4 + (1-\theta)\boldsymbol{\Phi}^6\boldsymbol{\beta}^6$, $\theta\in(0,1)$.
However, it remains difficult to ensure that further exploration will identify new nondominated points without a clear criterion for selecting appropriate mNBI parameters. In fact, as can be seen in Figure \ref{fig:ch6:ModifymNBI_2}(b), any choice of $\boldsymbol{\beta}$ such that the reference point lies on the line segment defined by  $\theta\boldsymbol{z}^4+ (1-\theta)\boldsymbol{\Phi}^6\boldsymbol{\beta}^{'6}$, $\forall \ 0 <\theta \leq 1$ returns the $\boldsymbol{z}^4$. Note that $\boldsymbol{z}^{6}$ is obtained by solving a modified mNBI subproblem \eqref{eq:ch6:simpleNBIn}, which we will discuss in the remaining part of this section and $\boldsymbol{\beta}'^{6}$ is calculated by using the objective constraints of the NBI subproblem (NBI$\beta$) for a given parameter $\boldsymbol{\bar{n}}^6$ and objective function value,  $\boldsymbol{f}=\boldsymbol{z}^6$.

To address this, two subproblems are introduced to fathom regions where no Pareto-optimal solutions exist. These are modifications of the original mNBI subproblem. For an arbitrary $p$-th facet denoted by $F^{p}_{S}(\boldsymbol{z}^{p1},\boldsymbol{z}^{p2})$ where  $z_1^{p1}<z_1^{p2}$, if the solution of subproblem \eqref{eq:ch6:mNBI} for given parameters ($\boldsymbol{\bar{n}}^k$, $\boldsymbol{\Phi}^k\boldsymbol{\beta}^k$) lies at $\boldsymbol{z}^{p1}$, then the following  subproblem is solved:

\begin{equation}\tag{mNBI$\bar{n}(a)$}\label{eq:ch6:simpleNBIn}
\begin{aligned}
&\begin{array}{c}
\max\limits_{\boldsymbol{x}\in X,t} \ t \\
\text { s.t. } \boldsymbol{\Phi}^k \boldsymbol{\beta}^k+t\boldsymbol{\bar{n}}^k\ge \boldsymbol{f}(\boldsymbol{x})-
\boldsymbol{f}^{id}\\
f_1(\boldsymbol{x}) \ge 
z^{p1}_1+\epsilon_z\\
\boldsymbol{x} \in \mathbb{R}^{n},\ t \in \mathbb{R}, \ \boldsymbol{\bar{n}}^k \in \mathbb{R}^{m}\\
\end{array}
\end{aligned}\end{equation}
where $\boldsymbol{z}^{p1}, \boldsymbol{z}^{p2}\in \boldsymbol{Z_E}$ are nondominated points obtained in previous iterations, $\boldsymbol{\bar{n}}^k$ is an unit normal vector to $F_{S}^p$  and $\epsilon_z$ is a small positive number.  

Similarly, for the case when the solution of subproblem \eqref{eq:ch6:mNBI} lies at extreme point $\boldsymbol{z}^{p2}$, the following subproblem is solved:
\begin{equation}\tag{mNBI$\bar{n}(b)$}\label{eq:ch6:simpleNBIn-b}
\begin{aligned}
&\begin{array}{c}
\max\limits_{\boldsymbol{x}\in X,t} \ t \\
\text { s.t. } \boldsymbol{\Phi}^k \boldsymbol{\beta}^k+t\boldsymbol{\bar{n}}^k\geq \boldsymbol{f}(\boldsymbol{x})-
\boldsymbol{f}^{id}\\
f_1(\boldsymbol{x}) \leq 
z^{p2}_1-\epsilon_z\\
\boldsymbol{x} \in \mathbb{R}^{n},\ t \in \mathbb{R}, \ \boldsymbol{\bar{n}}^k \in \mathbb{R}^{m}\\
\end{array}
\end{aligned}\end{equation}

\begin{theorem} \label{th:ch6:fathom}
Given two Pareto points $\boldsymbol{z}^{p1}$ and $\boldsymbol{z}^{p2}$, such that $z_{1}^{p1}<z_{1}^{p2}$, and the facet $F_{S}^{p}=(\boldsymbol{z}^{p1},\boldsymbol{z}^{p2})$ with normal vector $\boldsymbol{\bar{n}}^k$, if $(\boldsymbol{x}^{*}, \boldsymbol{z}^*,t^*)$ is a solution of subproblem \eqref{eq:ch6:simpleNBIn} and $\boldsymbol{z}^*$ is a nondominated point, then there is no nondominated point $\boldsymbol{z}=\left(z_1,z_2\right)^{\top}$ such that $z_1^{p1}+\epsilon_z< z_1 < z_1^{*}$. If $(\boldsymbol{x}^{*}, \boldsymbol{z}^*,t^*)$ is a solution of a subproblem \ref{eq:ch6:simpleNBIn-b} and $\boldsymbol{z}^*$ is a nondominated point, then there is no nondominated point $\boldsymbol{z}=\left(z_1,z_2\right)^{\top}$ such that $z_1^{*}< z_1 < z_1^{p2}-\epsilon_z$.
\end{theorem}

\begin{proof}
Let $\boldsymbol{z}^*$ be a nondominated point obtained at the solution $(\boldsymbol{x}^{*}, \boldsymbol{z}^*,t^*)$ of subproblem \eqref{eq:ch6:simpleNBIn} for facet $F_{S}^{p}=(\boldsymbol{z}^{p1},\boldsymbol{z}^{p2})$ as defined in Theorem \ref{th:ch6:fathom}, with parameters $\boldsymbol{\beta}^k = [\beta^k_1,1-\beta^k_1]$, $\boldsymbol{\Phi}^k$, and $\boldsymbol{\bar{n}}^k$.

Suppose there exists a nondominated point $\boldsymbol{z}'$ such that $z^{p_1}_{1}<z'_1 <z_1^*$. Because $\boldsymbol{z}'$ is nondominated, there must exist a point ($\boldsymbol{x}', \boldsymbol{z}', t'$) that is feasible for subproblem \ref{eq:ch6:simpleNBIn}. Then the first inequality of subproblem \ref{eq:ch6:simpleNBIn} must hold for both ($\boldsymbol{x}^{*}, \boldsymbol{z}^*, t^{*}$) and ($\boldsymbol{x}', \boldsymbol{z}', t'$): 
\begin{eqnarray}
 \boldsymbol{\Phi}^k \boldsymbol{\beta}^k+t^*\boldsymbol{\bar{n}}^k&\ge& \boldsymbol{z}^*-
\boldsymbol{f}^{id}\label{eq:ch6:profeq1}\\
 \boldsymbol{\Phi}^k \boldsymbol{\beta}^k+t'\boldsymbol{\bar{n}}^k&\ge& \boldsymbol{z}'-
\boldsymbol{f}^{id}\label{eq:ch6:profeq2}
\end{eqnarray}
By subtracting Equation \eqref{eq:ch6:profeq2} from Equation \eqref{eq:ch6:profeq1} and rearranging, we find:
\begin{equation}\label{eq:ch6:proof21}
\begin{array}{l}
(t^*-t^{'})\boldsymbol{\bar{n}}^k \geq \boldsymbol{z}^{*}-\boldsymbol{z}'
\end{array}
\end{equation}
where $\boldsymbol{\bar{n}}^k$ is the outer-facing normal vector of the facet $F^p_S$ and can be written as: \begin{equation}
\boldsymbol{\bar{n}}^k=-\boldsymbol{\Phi}^k \boldsymbol{e}=-\left[\begin{array}{cc}
0 & z^{p2}_{1}-z^{p1}_1 \\
z^{p1}_{2}-z^{p2}_2 & 0
\end{array}\right]\left[\begin{array}{l}
1 \\
1
\end{array}\right]=-\left[\begin{array}{c}
z^{p2}_{1}-z^{p1}_1  \\
z^{p1}_{2}-z^{p2}_2
\end{array}\right]
\end{equation}
Therefore, the components of $\boldsymbol{\bar{n}}^k$ are always negative and the non-principal elements of $\boldsymbol{\Phi}^k$ are positive. Note that the principal diagonal elements are always zero.
Equation \eqref{eq:ch6:proof21} is expressed component-wise as:
\begin{equation}\label{eq:ch6:proof21comp}
(t^{*}-t')\left[\begin{array}{c}
-(z^{p2}_{1}-z^{p1}_1)\\
-(z^{p1}_{2}-z^{p2}_2)
\end{array}\right]\\
\geq\left[\begin{array}{c}
z^{*}_{1}-z'_1 \\
z^{*}_{2}-z'_2
\end{array}\right]
\end{equation}
Since $t^{*}$ is a maximum value obtained as a solution of the problem \eqref{eq:ch6:simpleNBIn}, the strict inequality relation $t^{*} > t$ holds for all other feasible $t$. This means the components of $(t^*-t'){\boldsymbol{\bar{n}}^k}$ are always negative. Using Equation \eqref{eq:ch6:proof21comp}, where the components in the LHS are always negative, the elements of the RHS must be negative, meaning that $\boldsymbol{z}' > \boldsymbol{z}^*$ holds. Since we assumed that both $z^*$ and $z'$ are nondominated points, there exists no feasible Pareto-optimal solution $\boldsymbol{x}$ such that $f_j(\boldsymbol{x}) < z_j^*$ for all $j=1,2$, by Definition \ref{ParetoPoint}. However, this contradicts the assumption that both $\boldsymbol{z}^*$ and $\boldsymbol{z}'$ are nondominated points and $z^{p1}_1 <z'_1 < z^*_1$. Therefore $\boldsymbol{z}'$ cannot be a nondominated point and there exists no Pareto point $\boldsymbol{z}$  such that $z^{p1}_1<z_1<z^*_1$.
\end{proof}
 
Given a Pareto front, whenever the solution of subproblem \eqref{eq:ch6:mNBI} for the current facet $F^{p}_S=[\boldsymbol{z}^{p1}, \boldsymbol{z}^{p2}]$ generates a previously identified nondominated point, the facet is further explored by solving \eqref{eq:ch6:simpleNBIn} or \eqref{eq:ch6:simpleNBIn-b}, so that the empty part of the \subspace can be excluded from the search space. If the solution of \eqref{eq:ch6:simpleNBIn} or \eqref{eq:ch6:simpleNBIn-b} lies at one of the Pareto points obtained in a previous iteration, the entire facet represented by $F^{p}_S=[\boldsymbol{z}^{p1}, \boldsymbol{z}^{p2}]$ is discarded from the search space in a subsequent iteration. 

\subsection{Outline of the methodology}
In this section, we provide the pseudocode for the proposed SDNBI algorithm for the solution of BOO problems in Algorithm \ref{table:ch6:SDNBIalgorithm}.  The algorithm is built on the theory developed in previous sections.

In Step 1, the user specifies a convergence tolerance, $\epsilon$, as input to the algorithm. In Step 2, key quantities are initialized: the set of known Pareto points $\boldsymbol{Z_E}$ and the inner and outer approximations,  $IPS$ and $OPS$, respectively, are all empty. The approximation error $d_{max}$, the iteration counter $k$ and the \subspace counter $l$ are also initialized. The two anchor points, $\boldsymbol{z}^{A1}$ and $\boldsymbol{z}^{A2}$ with $z_1^{A2} < z_1^{A1}$, are also identified. This enables the definition of the first \subspace $C^0=(\boldsymbol{z}^{A2}, \boldsymbol{z}^{A1})$, which is assumed to be convex as indicated by setting the Boolean function $Is\_Convex(C^0)$ to $\mathtt{TRUE}$, and of the first facet $F_S^0=[\boldsymbol{z}^{A2}, \boldsymbol{z}^{A1}]$. There may be Pareto points in the \subspace defined by this facet (i.e., it $C^0$) so it is considered to be ``open" for further searching and the Boolean function $Is\_Open$ is set to $\mathtt{TRUE}$. The set of \subspaces ${\cal C}$ and the set of facets ${\cal F}$ are initialized. In Step 3, the facet is set to $F_S^p$, i.e., $p^1=0$.

The main loop begins at Step 4. The parameters of subproblem \eqref{eq:ch6:mNBI} at iteration $k$ are specified in Step 5. Subproblem \eqref{eq:ch6:mNBI} is then solved in Step 6, yielding $(\boldsymbol{z}^*,t^*)$. In Steps 7 to 13, the fathoming of facet $F_S^{p^k}$ is considered as described in Section~\ref{sec-fathom}. Specifically, if $\boldsymbol{z}^*$ is a member of $\boldsymbol{Z_E}$, i.e., a previously determined Pareto point, subproblem \eqref{eq:ch6:simpleNBIn} or \eqref{eq:ch6:simpleNBIn-b} is solved to generate a new $(\boldsymbol{z}^*,t^*)$. If the resulting vector $\boldsymbol{z}^*$ is once more a known Pareto point, the entire facet is discarded, i.e., $Is\_Open(F^{p^k}_S)$ is set to $\mathtt{FALSE}$.

If a new Pareto point has been found, the convexity property of the Pareto front in the current \subspace $C^{l^k}$ is examined in Steps 15-22 and new \subspaces are created as required, based on Section \ref{sec-decomposition}. Specifically, if all the Pareto points in the \subspace, $\boldsymbol{z} \in \boldsymbol{Z_E}^{C^{l^k}}$, are found to lie on the half-plane consistent with the convexity / nonconvexity assumption on $C^{l^k}$, then $\boldsymbol{z}^*$ can simply be added to the \subspace and the set of facets updated to take this new point into account. The convexity / nonconvexity of $C^{l^k}$ is unchanged. On the contrary, if there exists at least one point in the other half-plane, the set of Pareto points $\boldsymbol{Z_E}^{C^{l^k}}$ is updated by adding $\boldsymbol{z}^*$ and ${C^{l^k}}$ is partitioned to create \subspaces that only contain Pareto points that lie in the same half-space. The set of facets is updated as is the number of \subspaces.

The iteration counter is then increased.
The approximations $IPS$ and $OPS$ are updated in Step 25 and if there are no remaining facets to explore, the algorithm terminates regardless of whether the stopping tolerance has been achieved. Otherwise, the next facet $F_S^{p^k} = [\boldsymbol{z}^{p_2^{k}},\boldsymbol{z}^{p_1^k}]$ to be explored is identified in Steps 28 and 29 as the open facet at which the largest approximation error is observed. The corresponding \subspace is also identified, defined as the \subspace $C^{l^k} \in {\cal C}$ such that $\boldsymbol{z}^{p_2^{k}},\boldsymbol{z}^{p_1^k} \in C^{l^k}$.
The algorithm terminates when the maximum approximation error is less than $\epsilon$.

\begin{algorithm}[H]\setstretch{0.75}
\caption{SDNBI Algorithm for bi-objective optimization}\label{table:ch6:SDNBIalgorithm}
\begin{algorithmic}[1]
\NonumProcedure{SDNBI algorithm}{}
\Inputs{Set a quality threshold, $\epsilon$.}
\Initialize{\strut$\boldsymbol{Z_E} = IPS = OPS = \emptyset$; $d_{max} \gets \infty $, \ $k \gets 1$, $dc \gets 0$, $N_F \gets 1$. \\
Find anchor points $\boldsymbol{z}^{A1}$ and $\boldsymbol{z}^{A2}$ with $z_1^{A2} < z_1^{A1}$. \\
Define an initial search region, $C^0(\boldsymbol{z}^{A2},\boldsymbol{z}^{A1})$, set $Is\_Convex(C^0) = \mathtt{TRUE}$ and ${\cal C}=\{C^0\}$. \\
Define facet $F_S^0=[\boldsymbol{z}^{A2},\boldsymbol{z}^{A1}]$, set $Is\_Open(F_S^0)=\mathtt{TRUE}$ and  ${\cal F}=\{F_S^0\}$.}
\State Set $\boldsymbol{Z_E}=\{\boldsymbol{z}^{A1},\boldsymbol{z}^{A2}\}$, set $p^{k}=0$.
\While{$d_{max} \geq \epsilon$ }
\State \parbox[t]{\dimexpr\linewidth-\algorithmicindent}{Set the next normal $\boldsymbol{\bar{n}}^k$ to be the outer normal to $F^{p^k}_{S}$ and set the reference point $\boldsymbol{\Phi}^k\boldsymbol{\beta}^k$ by choosing $\boldsymbol{\beta}^k$ to be the midpoint of $F^{p^k}_{S}$.\strut}
\State Determine $(\boldsymbol{z}^*, t^{*})$ by solving \eqref{eq:ch6:mNBI} with $\boldsymbol{\Phi}^k\boldsymbol{\beta}^k$ and $\boldsymbol{\bar{n}}^k$.

 \If{($\boldsymbol{z}^k \in \boldsymbol{Z_E}$)}
 \State \parbox[t]{\dimexpr\linewidth-\algorithmicindent}{Obtain new $(\boldsymbol{z}^*,t^*)$ by solving \eqref{eq:ch6:simpleNBIn} or \eqref{eq:ch6:simpleNBIn-b}.\strut}
 \If{($\boldsymbol{z}^* \in Z_E$)}
  \State Fathom facet $F^{p^k}_{S}$ by setting $Is\_Open(F^{p^k}_S)=\mathtt{FALSE}$.
  \State Set  $k=k+1$ and return to step 6.
 \EndIf
 \Else
 \State \parbox[t]{\dimexpr\linewidth-\algorithmicindent}{Define the supporting line $H_2(\boldsymbol{w}'^k,b^k):b^k=\boldsymbol{w}'^{k \top}\boldsymbol{z}^*$ using Equation \eqref{eq:ch6:normal-mNBI}.\strut}
 \If {($Is\_Convex(\boldsymbol{Z}^{C^{l^k}}_{\boldsymbol{E}})==\mathtt{TRUE} \; \mathrm{and} \; \boldsymbol{w}'^{k \top}\boldsymbol{z} \ge {b^k}$ for all  $\boldsymbol{z}\in\boldsymbol{Z}^{C^{l^k}}_{\boldsymbol{E}})$ \\or ($Is\_Convex(\boldsymbol{Z}^{C^{l^k}}_{\boldsymbol{E}})==\mathtt{FALSE} \; \mathrm{and} \; \boldsymbol{w}'^{k \top}\boldsymbol{z} \le {b^k}$ for all $\boldsymbol{z}\in\boldsymbol{Z}^{C^{l^k}}_{\boldsymbol{E}}$)}
  \State \parbox[t]{\dimexpr\linewidth-\algorithmicindent}{Add $\boldsymbol{z}^*$ to $\boldsymbol{Z}^{C^{l^k}}_{\boldsymbol{E}}$ and replace $F_S^{p^k}\in {\cal F}$ by two new facets $F_S^{p_1^k}$ and $F_S^{p_2^k}$ \\ in which $\boldsymbol{z}^*$ is an extreme point. Set $Is\_Open(F_S^{p_1^k})=Is\_Open(F_S^{p_2^k})=\mathtt{TRUE}$.}
  \Else
  \State Update $\boldsymbol{Z}^{C^l}_{\boldsymbol{E}}=\boldsymbol{Z}^{C^{l^k}}_{\boldsymbol{E}}\cup \{ \boldsymbol{z}^* \}$
  \State \parbox[t]{\dimexpr\linewidth-\algorithmicindent}{Decompose current \subspace $C^{l^k}$ into two or more \subspaces by investigating \\ all nondominated points $\boldsymbol{z} \in\boldsymbol{Z}^{C^{l^k}}_{\boldsymbol{E}}$, and set corresponding $Is\_Convex$.\strut}
  \State Update facets for the \subspaces that contain $\boldsymbol{z}^*$.
  \State Set $dc=dc+N^l$, where $N^l$ is the number of new \subspaces created.
  \EndIf
  \EndIf
\State Set $k=k+1$.
\State \parbox[t]{\dimexpr\linewidth-\algorithmicindent}{For all \subspaces $C^{l}\in {\cal C}$, $l=1,...,dc$, construct $IPS$ and $OPS$.\strut}
\If{($Is\_Open{F^p_S}==\mathtt{FALSE}, \ \forall p$)} \textbf{break} \EndIf
 
\State \parbox[t]{\dimexpr\linewidth-\algorithmicindent}{ Compute the approximation error $d_{error,p}$ between $IPS$ and $OPS$ \\ for each  facet $F^{p}_{S}$ such that $Is\_Open(F_S^p)=\mathtt{TRUE}$, $p=1,\dots,N_F$. }
\State \parbox[t]{\dimexpr\linewidth-\algorithmicindent}{Select the facet  $F^{p^k}_{S}$ that presents the largest error, i.e., $p^k= \argmax\limits_{p=1,...,N_F}(d_{error,p})$,\\   and the corresponding \subspace $C^{l^k}$. Set $d_{max}=d_{error,p^k}$. \strut}
\EndWhile
\NonumEndProcedure{}{}
\end{algorithmic}
\end{algorithm}

\section{Performance of the algorithm}\label{sec:benchmark}
In this section, the numerical experiments to evaluate the performance of the proposed algorithm are discussed, and its behavior is compared with the SD algorithm and mNBI method for a set of test cases.

\subsection{Test functions}
To assess the performance of the three algorithms studied, five well-known problems are selected from the MOP1 \citep{van1999multiobjective}, ZDT \citep{Zitzler2000}, SCH \citep{schaffer1986some}, and TNK \citep{tanaka1995ga} test suites developed for testing evolutionary MOO algorithms. The test problems are chosen to vary in complexity in terms of problem size and numerical difficulty with convex, nonconvex and disconnected Pareto fronts. The test problems are summarized in Table \ref{table:ch6:TestProblem} where the number of decision variables $n$ (and of integer variables), the problem class, and the nature of the Pareto-optimal front are specified. Problem ZDT5 is a mixed-integer problem and the integer variables are denoted by the vector $\boldsymbol{y} \subseteq \mathbb{N}^{n_2}$, to distinguish them from the continuous variables, denoted by $\nu_i$, $i=1,\dots,11$, $g$  and $h$ in this formulation.

Before formulating the scalarized problems \eqref{eq:ch6:weightedsum} and \eqref{eq:ch6:mNBI}, each objective function is normalized with respect to the limits of objective space, i.e., the ideal and nadir points, $\boldsymbol{f}^{id}$ and $\boldsymbol{f}^{nd}$, to avoid biasing search direction towards a particular objective. Note that the mNBI method has been proven to be independent of the relative scales of the objective functions \citep{DASDENNIS1998}, but the normalization is applied to all methods in order to make a fair comparison. The $m$-dimensional objective vector of normalized functions is denoted by $\hat{\boldsymbol{f}}$ and each component is defined as follows:

\begin{equation}\label{eq:ch6:normalisedF}
\begin{aligned}
& \hat{f}_{j}(\boldsymbol{x}) = \frac{f_{j}(\boldsymbol{x})-f_{j}^{{nd}}(\boldsymbol{x})}{ f_{j}^{{id}}(\boldsymbol{x})-f_{j}^{nd}(\boldsymbol{x})},\ j\textrm{=1},\cdots,m.
\end{aligned}
\end{equation}

\begin{sidewaystable}[]
\caption{Two-objective test problems from the literature used in this work}\label{table:ch6:TestProblem}
\renewcommand{\arraystretch}{1}
\renewcommand{\tabcolsep}{8mm}
\resizebox{\columnwidth}{!}{%
\begin{tabular}{llcccc}
\hline
Name & Formulation & Number of variables, $n$ & Problem & Type of & Reference \\ 
&                &  (Number of integer variables)& class & Pareto front \\ \hline
MOP1 & $
\begin{array}{l}
\min\limits_{\boldsymbol{x}}\ f_{1}(\boldsymbol{x})=\frac{1}{x_1^2+x_2^2+1}\\
\min\limits_{\boldsymbol{x}}\  f_{2}(\boldsymbol{x})=x_1^2+3x_2^2+4\\
\text{s.t} \quad -3 \leq x_1 \leq 3, \ -3 \leq x_2 \leq 3 \\
\quad 
\end{array}$
 & 2 (0) & NLP &  \begin{tabular}[c]{@{}c@{}}Convex,\\ continuous \end{tabular} & \cite{van1999multiobjective} \\

SCH2 & $
\begin{array}{l}
\min\limits_{\boldsymbol{x}}\ f_{1}(\boldsymbol{x})=(x_1+x_2-7.5)^{2}+\frac{(x_2-x_1+3)^{2}}{4}\\
\min\limits_{\boldsymbol{x}}\  f_{2}(\boldsymbol{x})=\frac{(x_1-1)^{2}}{4}+\frac{(x_2-4)^{2}}{2}\\
\text{s.t} \quad \frac{(x_1-2)^{3}}{2}+x_2-2.5 \leq 0 \\
\quad\quad x_1+x_2-8(x_2-x_1+0.65)^{2}-3.85 \leq 0 \\
\quad\quad 0 \leq x_1 \leq 5, \ 0 \leq x_2 \leq 3 \\
\quad 
\end{array}
$ & 2  (0) &  NLP & \begin{tabular}[c]{@{}c@{}}Nonconvex,\\ disconnected \end{tabular}  & \cite{schaffer1986some}\\

TNK & $\begin{array}{ll}
\min\limits_{\boldsymbol{x}} \ f_{1}(\boldsymbol{x})=x_{1} \\ \min\limits_{\boldsymbol{x}} \ f_{2}(\boldsymbol{x})=x_{2} \\
\text { s.t. } \ x_{1}^{2}+x_{2}^{2}-1-0.1 \cos \left(16 \arctan \frac{x_{1}}{x_{2}}\right) \geq 0 \\
\quad\quad \left(x_{1}-0.5\right)^{2}+\left(x_{2}-0.5\right)^{2} \leq 0.5 \\
\quad\quad 0 \leq x_{1} \leq \pi, \ 0 \leq x_{2} \leq \pi\\
\quad
\end{array}$ & 2  (0) &  NLP & \begin{tabular}[c]{@{}c@{}}Nonconvex,\\ continuous \end{tabular}  & \cite{tanaka1995ga}\\

ZDT3 & $
\begin{array}{l}
\min\limits_{\boldsymbol{x},g,h} \ f_{1}(\boldsymbol{x})=x_{1}\\ \min\limits_{\boldsymbol{x},g,h}f_{2}\left(\boldsymbol{x}\right)=g\cdot h\\
\text { s.t. } g=1+\frac{9}{n-1} \sum_{i=2}^{11} x_{i} \\
\quad\quad h =1-\sqrt{x_{1} / g}-\left(x_{1} / g\right) \sin \left(10 \pi x_{1}\right)\\
\quad\quad 0 \leq x_i \leq 1, \   i=1,\dots,30 \\
\quad \quad g \in \mathbb{R}, h \in \mathbb{R}
\\\quad 
\end{array}
$ & 32  (0) &  NLP & \begin{tabular}[c]{@{}c@{}}Nonconvex,\\ disconnected\end{tabular} & \cite{Zitzler2000}  \\

ZDT5 & $
\begin{array}{l}
\min\limits_{\boldsymbol{y}, \boldsymbol{\nu}, g, h}\ f_{1}(\boldsymbol{y}) =1+ y_{1} \\ 
\min\limits_{\boldsymbol{y}, \boldsymbol{\nu}, g, h} f_{2}(\boldsymbol{y})=g\cdot h\\
\text { s.t. } g =\sum_{i=2}^{11} \nu_i, h=\frac{1}{1+y_1} \\
\quad\quad \nu_i =\left\{\begin{array}{ll}
2+y_i & \text { if } y_i<5 \\
1 & \text { if } y_i=5 \\
\end{array}\right\}, \   i=2,\dots, 11 \\
\quad\quad y_1 \in \{0,\ldots,30\}, \ y_{i}\in \{0,\ldots,5\}, \ i=2,\dots,11 \\
\quad\quad g , h, \nu_i \in \mathbb{R} \\
\end{array}
$ & 23 (11)&  MINLP & \begin{tabular}[c]{@{}c@{}}Integer, \\ convex when \\ integrality is relaxed \end{tabular} & \cite{Zitzler2000} \\ \hline
\end{tabular}
} 
\end{sidewaystable}

\subsection{Performance metrics}
The main goal when solving MOO problems is to generate, in the minimum time possible, a diverse set of nondominated points that are evenly distributed along the Pareto front. For nonconvex problems, there is the added challenge of ensuring that the points found are as close as possible to the true Pareto front.   

Based on these considerations, we use the following performance metrics for evaluating the quality of the solution set obtained by each algorithm.

\begin{enumerate}
\item The number of unique non-dominated solutions ($N_{unq}$):  $N_{unq}$ denotes the cardinality of the set of unique Pareto points obtained in a given run, serving as a measure of the diversity of the solution set.

\item Hypervolume (HV): The hypervolume \citep{Zitzler2003} of a set of non-dominated solutions is the volume of the $m$-dimensional region in objective space enclosed by the  nondominated solutions obtained and a reference point, $f_{ref},j=1,...,m$.  The larger HV, the better the solution set in terms of convergence to the true Pareto front and/or in terms of diversity of the solutions. In the evaluation of the hypervolume in this work, the reference point is chosen by selecting the worst (largest) value of each objective across all nondominated solutions generated by the three algorithms considered.

\item Distribution metric (DM):  The distribution metric was suggested by \cite{zheng2017new} to capture the spread of the solution set over the Pareto front approximation and the extent of the true Pareto front covered by the nondominated points. This metric addresses some of the deficiencies of the metrics of \citet{schott1995fault} and \citet{wu2001metrics}. Smaller values of DM indicate better distributed solution sets. It is given by the following set of equations:
\begin{equation}\label{eq:ch6:DM}
\begin{aligned}
\begin{array}{c}
    
DM(\boldsymbol{Z_E})=\dfrac{1}{|\boldsymbol{Z_E}|} \sum\limits_{j=1}^{m}\left(\dfrac{\sigma_{j}}{\tau_{j}}\right)\left(\dfrac{\left|f^{id}_{j}-f^{nd}_{j}\right|}{R_{j}}\right),\\ [10pt]
\sigma_{j}=\dfrac{1}{|\boldsymbol{Z_E}|-2} \sum_{e=1}^{|\boldsymbol{Z_E}|-1}\left(d_{e}^{j}-\tau_j\right)^{2}, \  \tau_{j}=\dfrac{1}{|\boldsymbol{Z_E}|-1} \sum_{e=1}^{|\boldsymbol{Z_E}|-1} d_{e}^{j}, \; j = 1,\dots,m,\\ [10pt]
R_{j}=\max\limits_{\boldsymbol{x} \in \boldsymbol{X_E}} f_{j}(\boldsymbol{x})-\min\limits_{\boldsymbol{x} \in \boldsymbol{X_E}} f_{j}(\boldsymbol{x}), \; j=1,\dots,m,
\end{array}
\end{aligned}
\end{equation}
where $|\boldsymbol{Z_E}|$ is the number of nondominated points, $d^j_{e}, \ e=1, \dots, |\boldsymbol{Z_E}|-1$, is the absolute difference between the values of the $j^{\mathrm{th}}$ objective at two adjacent nondominated points in objective space, $\sigma_j$ and $\tau_j$ are the standard deviation and mean of $d^j_{e}$, respectively, over all $e=1,\dots,|\boldsymbol{Z_E}|$, for the $j^{\mathrm{th}}$ objective, $\sigma_j/\tau_j$ is the coefficient of variance for the $j^{\rm th}$ objective, $R_j$ is the difference between the maximum and minimum values of the $j^{\mathrm{th}}$ objective over all points in $\boldsymbol{Z_E}$, and $\boldsymbol{X_E}$ is the set of all Pareto solutions that have been identified.

\item CPU time: The computational time to approximate the Pareto frontier is a critical aspect for computationally expensive MOO problems. Both the total CPU time to generate all solutions ($t_{cpu,t}$) and the average CPU time taken to generate one nondominated point ($t_{cpu,a}$), defined as $t_{cpu,t}/N_{unq}$ are reported.

\end{enumerate}

\subsection{Implementation overview}
All algorithms are implemented in Matlab 2018a using the same common subfunctions. For the mathematical benchmark problems, each subproblem derived from the scalarization method is solved through the GAMS modeling environment interfaced with CONOPT \citep{drud1994conopt} and SNOPT \citep{murray1997sequential,gill2005snopt} for nonlinear programs and DICOPT  \citep{grossmann2002gams} for the mixed-integer programs. Both are local solvers. The codes are available at \url{https://github.com/MolChemML/MOO-SDNBI}.

To increase the likelihood of identifying globally optimal Pareto points, a multi-start approach \citep{pal2013comparison} is adopted in the solution of the subproblems. The multiple starting points ($N_{sobol}$) are generated via a Sobol' sequence \citep{Sobol2011} to ensure diversity in the design space. For the test functions, a deterministic global optimization algorithm such as BARON \citep{ryoo1995BARON} could be used instead. 

For the SD algorithm and mNBI method, the stopping criterion is defined as a fixed number of iterations $N_{iter,fix}$, selected based on the iterations taken for SDNBI to reach a pre-defined error tolerance. This is because, although the SD algorithm benefits from a stopping criterion based on a convergence tolerance, it tends to terminate early as it cannot explore nonconvex regions of the Pareto front. Furthermore, the only stopping criterion for the mNBI algorithm is related to the computational effort with no means to estimate the quality of the approximation of the Pareto front generated. 

The initial set, $C^{ref}$, of $N_{\beta}-1$ reference points for the mNBI method is chosen such that consecutive reference points $\boldsymbol{\Phi\beta}^{k}$ and $\boldsymbol{\Phi\beta}^{k+1}$ are equally spaced, with a spacing $\delta$. This is  expressed as:
\begin{equation}\label{eq:ch6:beta_initial_set}
\beta^{k}_2= k \times \delta, \ \beta^{k}_1=1-\beta^{k}_2, \ \text{ for } k=1, ..., N_{\beta}-1, \delta=1/(N_{\beta}-1).
\end{equation}
Once the parameters in the initial set have been enumerated by the algorithm, the next $\boldsymbol{\beta}$ to be chosen, i.e., $\boldsymbol{\beta}^{N_{\beta}}$ is determined by taking the midpoint between two adjacent $\boldsymbol{\beta}$ vectors used in previous iterations of the algorithm. As previously mentioned, the values of the parameter $\boldsymbol{\beta}$ for the SDNBI algorithm are chosen as facet midpoints, i.e., $\boldsymbol{\beta}^k = [0.5, 0.5]^{\top}$. All other parameters used in the numerical tests are specified in Table \ref{Table:ch6:parameter}.

\begin{table}[]\centering
\caption{Algorithmic parameters used in the test problems}
\label{Table:ch6:parameter}
\resizebox{0.75\textwidth}{!}{%
\begin{tabular}{llllll}
\hline
Parameter             & MOP1     & SCH2       & TNK        & ZDT3       & ZDT5       \\ \hline
$N_{sobol}$             & 20       & 20         & 20         & 50         & 30         \\
$\epsilon$ in SDNBI   & 0.001    & 0.001      &  0.002    & 0.005       & 0.005       \\
$N_{iter,fix}$        & 33       & 27         & 59         & 36         & 40           \\
$N_{\beta}$            & 10       & 10         & 15         & 10         & 10           \\
$\boldsymbol{f}^{id}$ & {[}0.053, 1{]} & {[}7.251, 0.5{]} & {[}0.0416, 0.0416{]} & {[}0, -0.7733{]} & {[}0, 0.3226{]} \\
$\boldsymbol{f}^{nd}$ & {[}1, 37{]} & {[}18.5, 3.045{]} & {[}1.0384, 1.0384{]} & {[}0.8518, 1{]} & {[}31, 10{]} \\ \hline
\end{tabular}%
}
\end{table}

\section{Results and discussion}\label{sec:results}
\begin{figure}
	\centering
	\includegraphics[height=0.4\textheight]{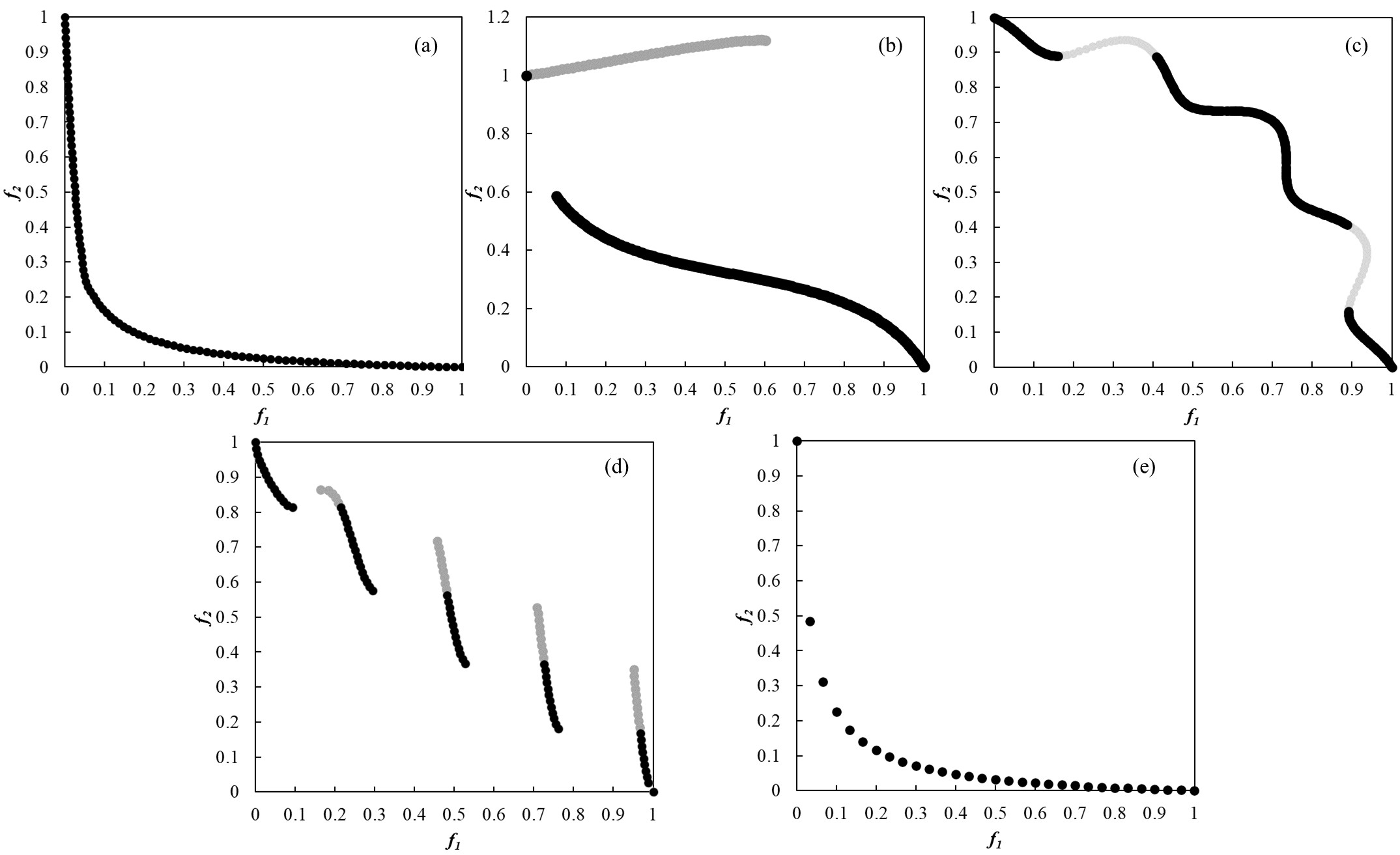}
	\vspace{-0.3cm}
	\captionsetup{font=small}\caption{Best-known Pareto front (black symbols $\bullet$)  and boundary points of the feasible region (grey symbols \textcolor{gray}{$\bullet$}) in bi-objective space for the test functions: (a) MOP1, (b) SCH2, (c) TNK, (d) ZDT3, and (e) ZDT5 obtained by solving problem (NBI$\beta$) for many reference points $\boldsymbol{\Phi\beta}$.}\label{fig:ch6:TrueParetofront}
\end{figure}

In this section, we compare the performance of the SD, mNBI, and SDNBI algorithms. The ``true" Pareto front and the boundary of the feasible region are generated by solving (NBI$\beta$) for MOP1, SCH2, TNK, and ZDT3 and \eqref{eq:ch6:mNBI} for ZDT5 using exhaustive enumeration of a large set ($N_{finite}$) of reference points $\boldsymbol{\Phi\beta}$ in order to provide a benchmark for the quality of the solutions. The results of these calculations are referred to as the best-known approximations of the true Pareto front and used for graphical comparison. The values of relevant metrics for the best-known Pareto front obtained for each test problem are shown in Table \ref{Table:ch6:benchmark-metrics}. Graphical representations of the Pareto fronts and boundary points in bi-objective space are given in Figure \ref{fig:ch6:TrueParetofront}.

\begin{table}[]\centering
\caption{Values of the relevant metrics for the best-known Pareto fronts of the test problems}
\label{Table:ch6:benchmark-metrics}
\begin{tabular}{llllll}
\hline
Quality Metrics & MOP1   & SCH2   & TNK    & ZDT3   & ZDT5   \\ \hline
$N_{finite}$    & 100    & 200    & 301    & 100     & 100     \\
$N_{unq}$       & 100    & 153    & 200    & 67     & 31     \\
HV ($10^{-2}$)  & 92.96  & 63.78  & 30.82  & 51.46     & 89.57     \\
DM              & 0.0157 & 0.0389 & 0.0238 & 0.0455 & 0.0944 \\ \hline
\end{tabular}%
\end{table}

\subsection{Test problem MOP1}
\begin{figure}
	\centering
	\includegraphics[height=0.2\textheight]{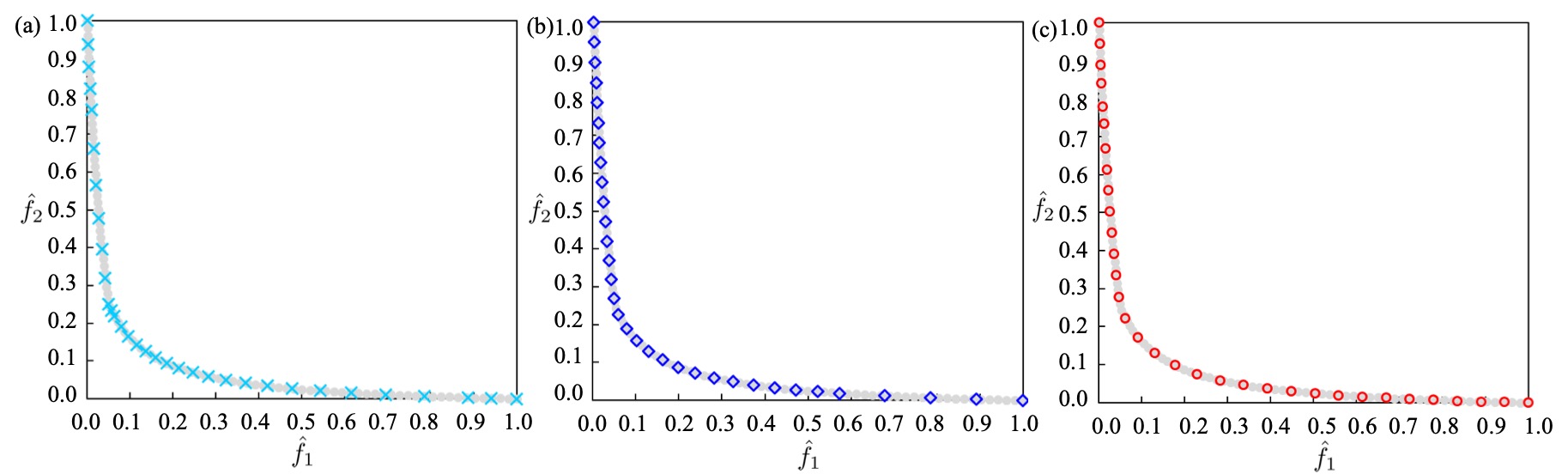}
	\vspace{-0.2cm}
	\captionsetup{font=small}\caption{Pareto points generated by (a) SD (\textcolor{skyblue}{$\boldsymbol{\times}$}), (b) mNBI (\textcolor{blue}{$\mathbf{\triangle}$}), and (c) SDNBI (\textbf{\textcolor{red}{\textbf{$\circ$}}}) for test problem MOP1. The best-known Pareto front is shown in grey.}\label{fig:ch6:MOP1-PF}
\end{figure}

The results generated by applying the SD, mNBI, and SDNBI algorithms to problem MOP1 and the relative performance of the algorithms are summarized in Table \ref{Table:ch6:MOP1SCH2}. The algorithms are run for 33 iterations, corresponding to the convergence of the SDNBI algorithm with the chosen value of $\epsilon$ (see Table \ref{Table:ch6:parameter}).
For test problem MOP1, all three algorithms show similar performance in terms of HV and $N_{unq}$, as can be seen in Table \ref{Table:ch6:MOP1SCH2}. This is attributed to the convexity of the Pareto front that allows the three algorithms to generate a new nondominated point at every iteration, resulting in an equal number of Pareto points at termination. While the SD algorithm is marginally better in terms of HV, the SDNBI algorithm yields the lowest DM, indicating better-distributed solutions along the Pareto front. This can be observed in Figure \ref{fig:ch6:MOP1-PF}, in which the approximate Pareto front at the 33$^{\rm rd}$ iteration is shown. This better performance in DM is mainly due to the fact that in the SDNBI algorithm, the next search direction is chosen by taking into account both spacing and an approximation error: the spacing between two adjacent Pareto points is used to set the parameter $\boldsymbol{\beta}$ as the midpoint of a facet and the facet itself is chosen based on the distance between the inner and outer approximation, which in turn determines the value of the parameters $\boldsymbol{\bar{n}$ and $\Phi}$, in order to minimize the error. By contrast, in the mNBI and SD algorithms, either spacing or error is considered, but not both.

\begin{figure}
	\centering
	\includegraphics[height=0.36\textheight]{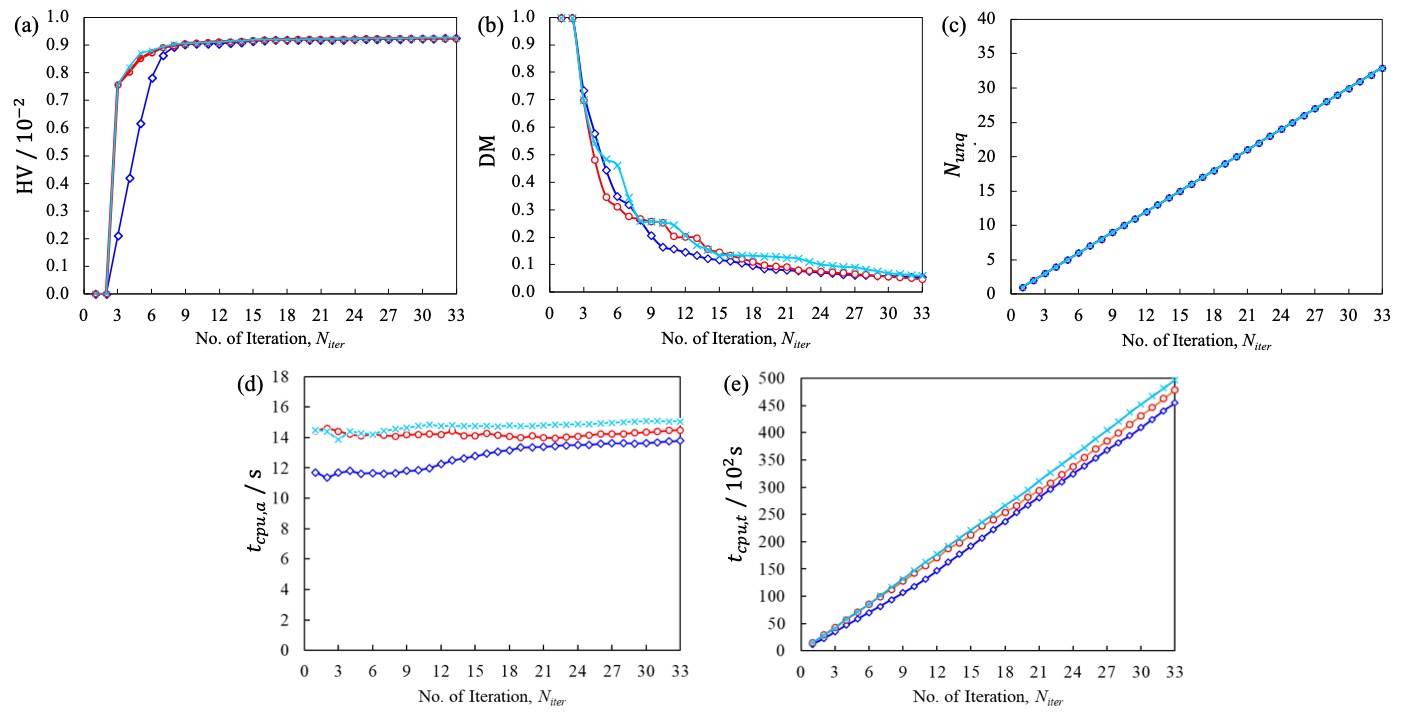}
	\vspace{-0.2cm}
	\captionsetup{font=small}\caption{Performance metrics as a function of iteration number for the application of the SD algorithm (\textcolor{skyblue}{--$\times$--}), mNBI algorithm (\textcolor{blue}{--$\diamond$--}), and SDNBI algorithm  (\textbf{\textcolor{red}{--$\circ$--}}) to MOP1 for 33 iterations: (a) Scaled hypervolume, HV, (b) Scaled distribution metric, DM, (c) Number of unique Pareto points, $N_{unq}$, (d) CPU time spent per iteration in seconds, $t_{cpu,a}$, (e) overall CPU time in seconds, $t_{cpu,t}$.} 
	\label{fig:ch6:MOP1-all}
\end{figure}

Although the mNBI algorithm outperforms the SDNBI in terms of HV, it is important to note that the accuracy of the mNBI-generated approximation depends highly on the initial set of  $\boldsymbol{\beta}$ values chosen. Yet, the number of $\boldsymbol{\beta}$ values required in the initial set and the number of additional values required to achieve a particular quality of solution set is unknown and there is no systematic method for determining an appropriate number. This is reflected in Figure \ref{fig:ch6:MOP1-all}(a) and \ref{fig:ch6:MOP1-all}(b) where the HV and DM values for the mNBI method exhibit the lowest initial performance. The SD and SDNBI algorithms, on the other hand, achieve high HV and low DM values that are close to the reference values (92.69 for HV and 0.0157 for DM) in few iterations, most likely because the scalarization parameters are chosen to maximize the accuracy of the approximation. 

The computational cost of generating Pareto points comes at a relatively higher price for SD and SDNBI (9.2\% and 5.2\% higher $t_{cpu,t}$ value, cf. Table \ref{Table:ch6:MOP1SCH2} and Figure \ref{fig:ch6:MOP1-all}(d) and \ref{fig:ch6:MOP1-all}(e)). The higher computational cost is mainly associated with the time taken for convex hull generation and for linear optimization runs for the calculation of the error between the inner and outer approximations of each facet.

\subsection{Test problem SCH2}
\begin{figure}[]
	\centering
	\includegraphics[height=0.23\textheight]{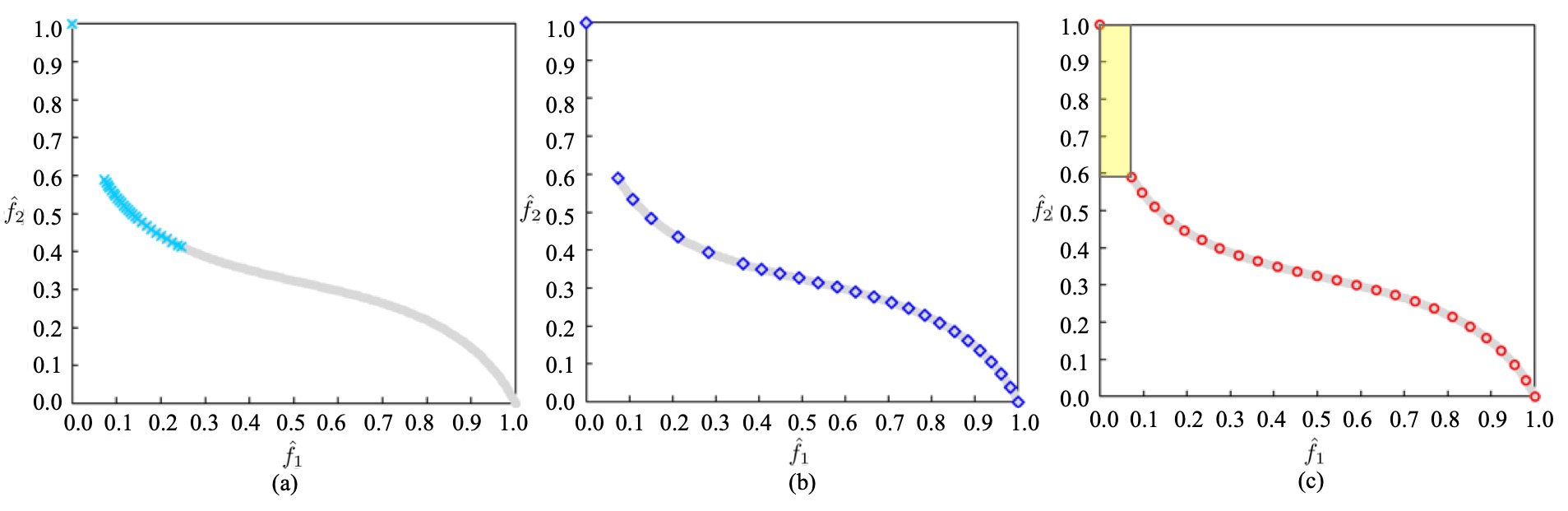}
	\vspace{-0.4cm}
	\captionsetup{font=small}\caption{Pareto points generated by (a) the SD algorithm (\textcolor{skyblue}{$\boldsymbol{\times}$}), (b) the mNBI algorithm (\textcolor{blue}{$\mathbf{\diamond}$}), and (c) the SDNBI algorithm (\textbf{\textcolor{red}{\textbf{$\circ$}}}) for test problem SCH2, including the region discarded as being devoid of Pareto points (yellow shading). The best-known Pareto front is shown in grey.}\label{fig:ch6:SHC2-PF}
\end{figure}

Test problem SCH2 is much more challenging than MOP1 because the Pareto-optimal set consists of two disconnected regions, as shown in Figure \ref{fig:ch6:SHC2-PF}. The first consists of a single point and the second is a large continuous, nonconvex region.  There is a large gap between these two regions.

The desired approximation quality is achieved with the SDNBI algorithm in 27 iterations. As shown in Table \ref{Table:ch6:MOP1SCH2}, the highest value of HV and the lowest value of DM are achieved when using SDNBI, indicating that the nondominated points generated by SDNBI are the closest approximation to the best-known Pareto front and offer a better distribution. This can also be observed by comparing the three panels in Figure \ref{fig:ch6:SHC2-PF}. It can further be seen that the points identified by the SD algorithm lie only on the convex part of the true Pareto front, as is expected, resulting in the lowest value of $N_{unq}$ and HV, and the highest value of DM. In view of this, subsequent discussion is focused on comparing the mNBI and SDNBI algorithms. 

\begin{table}[]
    \centering
    \caption{Performance metrics for the test problems MOP1 with $N_{iter,fix}=33$ and SCH2 with $N_{iter,fix}=27$, comparing the SD, mNBI, and SDNBI algorithms.}\label{Table:ch6:MOP1SCH2}
    \resizebox{0.7\textwidth}{!}{%
    \begin{tabular}{l|lll|lll}
    \hline
    Problem & \multicolumn{3}{c|}{MOP1} & \multicolumn{3}{c}{SCH2} \\ \hline
    Algorithm & SD & mNBI   & SDNBI              & SD & mNBI    &  SDNBI\\ 
    $N_{unq}$       & 33       & 33     & 33      & 26       & 25      & 26 \\
    HV ($10^{-2}$)  & 92.57    & 92.54  & 92.45  & 53.09    & 62.70   & 62.94\\
    DM              & 0.0586   & 0.0534 & 0.0478 & 0.2512   & 0.0998  & 0.0910\\
    $t_{cpu,a}$ (s) & 15.05    & 13.78  & 14.49   & 9.35     & 9.46    & 9.41  \\
    $t_{cpu,t}$ (s)  & 4.97$\times$10$^{2}$    & 4.55$\times$10$^{2}$  & 4.78$\times$10$^{2}$ & 2.43$\times$10$^{2}$   & 2.36$\times$10$^{2}$  & 2.44$\times$10$^{2}$  \\ \hline
    \end{tabular}%
    }
\end{table}

\begin{figure}[h]
	\centering
	\includegraphics[height=0.35\textheight]{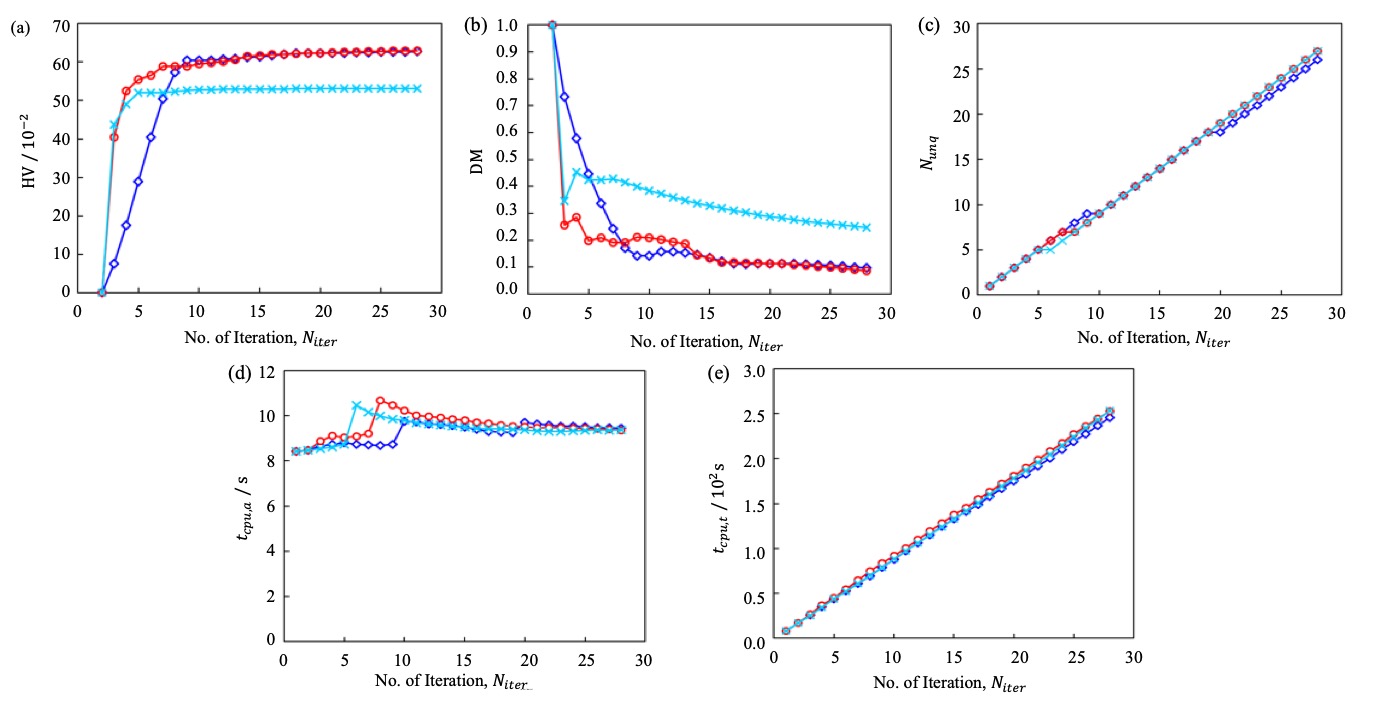}
	\vspace{-0.2cm}
	\captionsetup{font=small}\caption{Performance metrics as a function of iteration number for the application of the  SD (\textcolor{skyblue}{--$\times$--}), mNBI (\textcolor{blue}{--$\diamond$--}), and SDNBI (\textbf{\textcolor{red}{--$\circ$--}}) algorithms to SCH2 for 27 iterations:	(a) Scaled hypervolume, HV, (b) Scaled distribution metric, DM, (c) Number of unique Pareto points, $N_{unq}$, (d) CPU time spent per iteration in seconds, $t_{cpu,a}$, (e) overall CPU time in seconds, $t_{cpu,t}$.}
	\label{fig:ch6:SCH2-all}
\end{figure}

The overall performance of the mNBI algorithm is similar to that of SDNBI in terms of HV and DM, although the mNBI algorithm trails SDNBI in early iterations (see Figure \ref{fig:ch6:SCH2-all}(a) and Figure \ref{fig:ch6:SCH2-all}(b)). 

The better overall quality of the solution set obtained with SDNBI is partly as a result of recognizing discontinuities in the Pareto front. Once the corresponding facet has been fathomed, a new Pareto point can be generated at every subsequent iteration. 
For example, the disconnected region of the Pareto front between $0\leq \hat{f_1} \leq 0.07$ and also $0.59 \leq \hat{f_2} \leq 1$ (see yellow shaded box in Figure \ref{fig:ch6:SHC2-PF}) is removed from the search space at the 8$^{\rm th}$ iteration by solving an alternative SDNBI subproblem \eqref{eq:ch6:simpleNBIn-b}. This event is also observed as a kink in the plot of $N_{unq}$ for the SDNBI algorithm (Figure \ref{fig:ch6:SCH2-all}(c)). Fewer points are identified by the mNBI algorithm within the fixed number of iterations. This is mainly because, at several iterations, the solution of the mNBI subproblem fails to result in a new Pareto point as searching takes place along the large disconnected region. In fact, this region is repeatedly explored during the mNBI search procedure as iterations proceed. This algorithmic behavior is illustrated in Figure \ref{fig:ch6:SCH2-all}(c), where a lower value of $N_{unq}$ is observed with the mNBI algorithm than with the SDNBI algorithm after the 20$^{\rm th}$ iteration. 

Finally, all three algorithms exhibit similar performance in terms of CPU time per Pareto point and overall, as seen in Table \ref{Table:ch6:MOP1SCH2} and in Figure \ref{fig:ch6:SCH2-all}(d) and \ref{fig:ch6:SCH2-all}(e). This example illustrates that with the SDNBI, it is possible to obtain a better distribution of Pareto points  and to gain additional confidence that the Pareto front is described well, given the knowledge that no Pareto point exists in the relatively large discarded region.

\begin{table}[]
    \centering
    \caption{Performance metrics for the test problems:  TNK with $N_{iter,fix}=59$,  ZDT3 with $N_{iter,fix}=36$, and  ZDT5 with $N_{iter,fix}=40$, comparing the SD, mNBI, and SDNBI algorithms. SD is omitted for TNK as it does not identify any Pareto point.}\label{Table:ch6:TNKZDT35}
     \centering
    \resizebox{0.9\textwidth}{!}{%
    \begin{tabular}{l|ll|lll|lll}
    \hline
    Problem  & \multicolumn{2}{c|}{TNK} & \multicolumn{3}{c|}{ZDT3} & \multicolumn{3}{c}{ZDT5}  \\ \hline
    Algorithm  & mNBI   & SDNBI     & SD  & mNBI    & SDNBI         & SD  & mNBI    & SDNBI \\ 
    $N_{unq}$        & 56     & 59   & 36        & 35      & 36    & 23        & 23      & 31  \\
    HV ($10^{-2}$)   & 30.34  & 30.44 & 48.85     & 51.12   & 51.21   & 89.46     & 89.48   & 89.57\\
    DM               & 0.0674 & 0.0629 & 0.1356    & 0.0776  & 0.0667  & 0.1377   & 0.1230  & 	0.0944  \\ 
    $t_{cpu,a}$ (s)  & 15.88  & 16.24  & 45.06     & 47.38   & 46.10 & 73.40     & 52.86   & 51.92 \\
    $t_{cpu,t}$ (s)   & 8.89$\times$10$^{2}$ & 9.58$\times$10$^{2}$  & 16.22$\times$10$^{2}$   & 16.58$\times$10$^{2}$ & 16.59$\times$10$^{2}$ & 1688.20   & 1215.88 & 1609.64  \\\hline
    \end{tabular}%
    }
    \hspace{0.01\textwidth} 
%
    

%
  
%
\end{table}

\subsection{Test problems TNK and ZDT3}
Similar trends are observed for problems TNK and ZDT3 in which the nondominated sets comprise two disconnected regions for TNK and four for ZDT3, as shown in Figure \ref{fig:ch6:TrueParetofront}(c) and \ref{fig:ch6:TrueParetofront}(d). 

A key challenge in solving these two problems is to determine a set of parameter values so that the Pareto points identified are evenly distributed and, more importantly, so that the solution of the optimization subproblems that result in convergence to an already known Pareto point is avoided.  

\begin{figure}
	\centering
	\includegraphics[height=0.3\textheight]{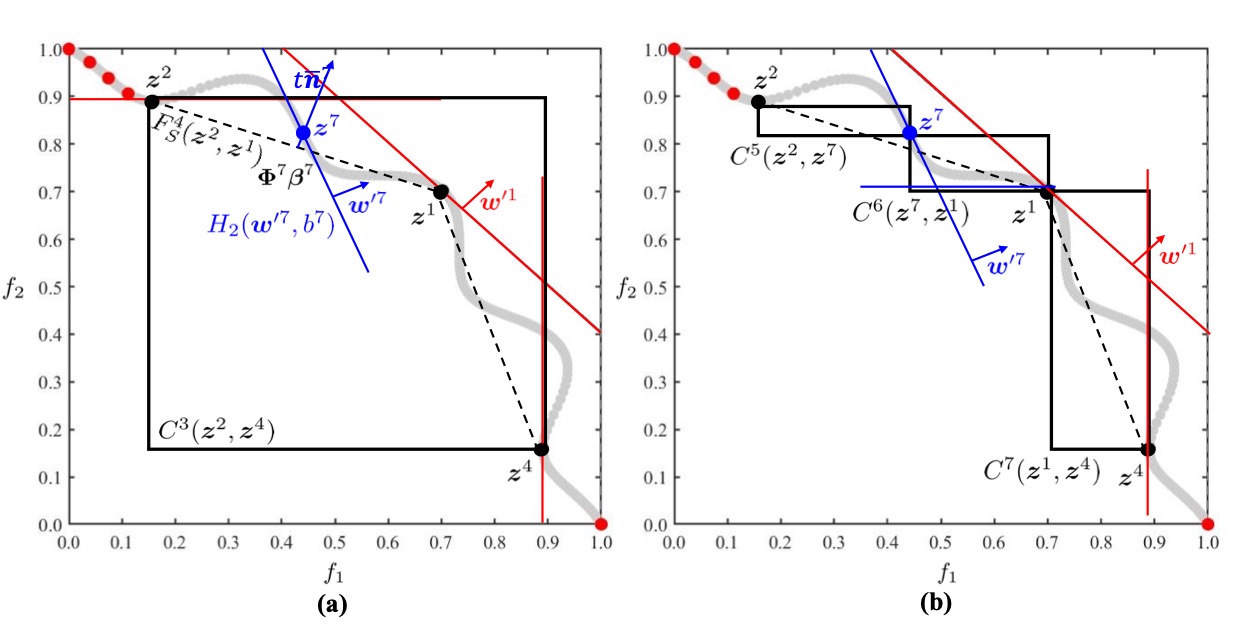}
	\vspace{-0.2cm}
	\captionsetup{font=small}\caption{A geometrical illustration of the SDNBI procedure for problem TNK. (a) A new Pareto point $\boldsymbol{z}^{7}$ is generated as a solution of \eqref{eq:ch6:mNBI} for facet $F^4_S$ of \subspace $C^3(\boldsymbol{z}^{2},\boldsymbol{z}^{4})$ ($Is\_Convex(C^3)=\mathtt{FALSE}$). The black dashed line represents a facet obtained by convexhull($\boldsymbol{Z}^{C^3}_{\boldsymbol{E}}$). The solid red and blue lines represent supporting lines at the solutions. (b) The \subspace is decomposed into three \subspaces, such that all points in each \subspace satisfy Equation \eqref{eq:ch6:subconv} or \eqref{eq:ch6:subnonconv}. The \subspaces are $C^{5}(\boldsymbol{z}^{3},\boldsymbol{z}^{4})$, with $Is\_Convex(C^5)=\mathtt{FALSE}$,  $C^{6}(\boldsymbol{z}^{3},\boldsymbol{z}^{4})$,  with $Is\_Convex(C^5)=\mathtt{TRUE}$, and $C^{7}(\boldsymbol{z}^{3},\boldsymbol{z}^{4})$, with $Is\_Convex(C^5)=\mathtt{FALSE}$, respectively.}
	
	\label{fig:ch6:TNK-decompo}
\end{figure}

In Figure \ref{fig:ch6:TNK-decompo}, we display the 7$^{\rm th}$ iteration of the SDNBI algorithm when applied to problem TNK, to exemplify how a \subspace is decomposed and information on the shape or convexity of the Pareto front is obtained. At the 7$^{\rm th}$  iteration, the \subspace $C^{3}(\boldsymbol{z}^{2},\boldsymbol{z}^{4})$ is selected as containing the facet with the largest approximation error (facet $F^{4}_{S}=[\boldsymbol{z}^2,\boldsymbol{z}^1]$). The \subspace comprises three nondominated points,  $\boldsymbol{Z}^{C^3}_{\boldsymbol{E}} =\{\boldsymbol{z}^{2},\boldsymbol{z}^{1},\boldsymbol{z}^{4}\}$ and is assumed to be nonconvex, i.e., $Is\_Convex(C^3)=\mathtt{FALSE}$. We begin by constructing inner and outer approximations as shown in Figure \ref{fig:ch6:TNK-decompo}. The red lines are tangents (inner approximations) at each nondominated point and the black dashed lines represent the outer approximations generated by convexhull($ \boldsymbol{Z}^{C^3}_{\boldsymbol{E}}$). 
The mNBI parameters ($\boldsymbol{\Phi}^7\boldsymbol{\beta}^7,\boldsymbol{\beta}^7$) are selected for facet $F_S^4$ and the corresponding  \eqref{eq:ch6:mNBI} subproblem is solved. From visual inspection of Figure \ref{fig:ch6:TNK-decompo}(a), it is clear, given the newly generated point $\boldsymbol{z}^{7}$,  $\boldsymbol{w}'^{7 \top}\boldsymbol{z}\leq b^{7}$ does not hold for $\boldsymbol{z}=\boldsymbol{z}^1$ and $\boldsymbol{z}=\boldsymbol{z}^4$, indicating that the Pareto front in \subspace $C^3$ 
is partly convex. Therefore, \subspace $C^{3}(\boldsymbol{z}^{3},\boldsymbol{z}^{4})$ is decomposed into three \subspaces:  $C^{5}(\boldsymbol{z}^{3},\boldsymbol{z}^{4})$, with $Is\_Convex(C^5)=\mathtt{FALSE}$,  $C^{6}(\boldsymbol{z}^{3},\boldsymbol{z}^{4})$,  with $Is\_Convex(C^5)=\mathtt{TRUE}$, and $C^{7}(\boldsymbol{z}^{3},\boldsymbol{z}^{4})$, with $Is\_Convex(C^5)=\mathtt{FALSE}$, such that all points in each \subspace satisfy Equation \eqref{eq:ch6:subconv} or \eqref{eq:ch6:subnonconv}. 

\begin{figure}[b]
	\centering
	\includegraphics[height=0.21\textheight]{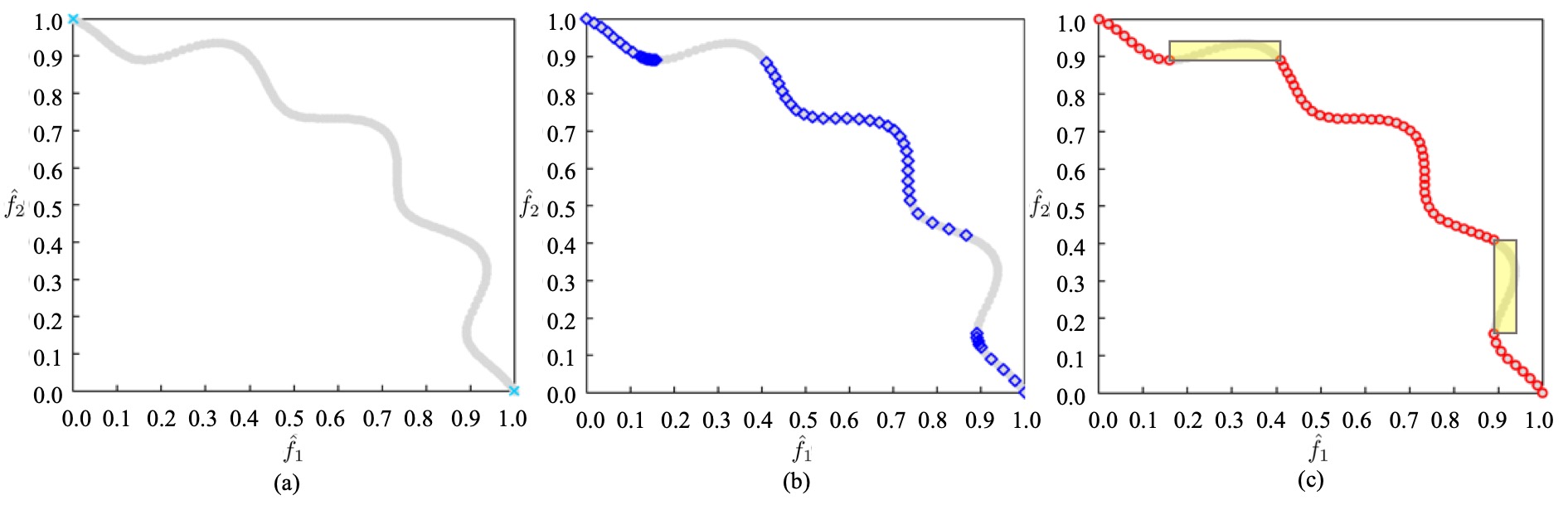}
	\vspace{-0.2cm}
	\captionsetup{font=small}\caption{Pareto points generated by (a) the SD algorithm (\textcolor{skyblue}{$\boldsymbol{\times}$}), (b) the mNBI method (\textcolor{blue}{$\mathbf{\diamond}$}), and (c) the SDNBI algorithm (\textbf{\textcolor{red}{\textbf{$\circ$}}}) for test problem TNK, including the regions discarded as being devoid of Pareto points (yellow shading). The best-known Pareto front is shown in grey.}\label{fig:ch6:TNK-PF}
\end{figure}

\begin{figure}[t]
	\centering
	\includegraphics[height=0.34\textheight]{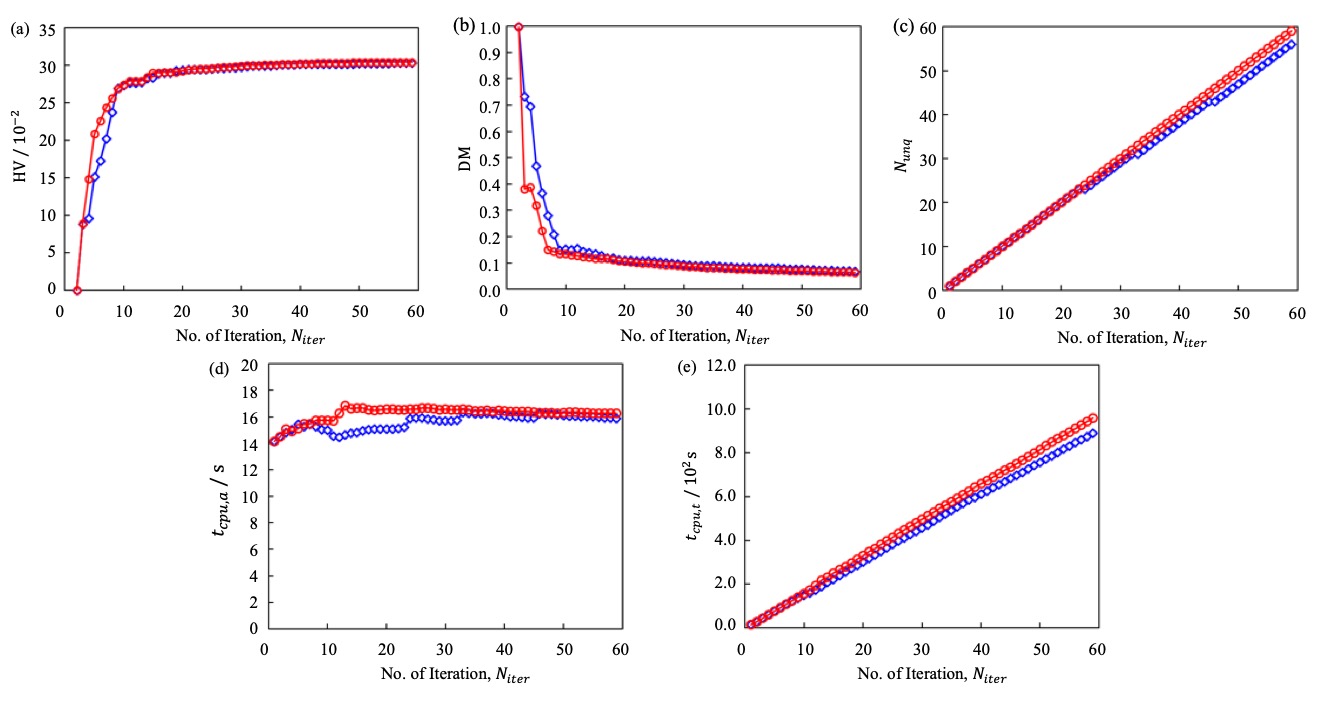}
	\vspace{-0.2cm}
	\captionsetup{font=small}\caption{Performance metrics as a function of iteration number for the application of the  mNBI (\textcolor{blue}{--$\diamond$--}) and SDNBI (\textbf{\textcolor{red}{--$\circ$--}}) algorithms to problem TNK for 59 iterations:	(a) Scaled hypervolume, HV, (b) Scaled distribution metric, DM, (c) Number of unique Pareto points, $N_{unq}$, (d) CPU time spent per iteration in seconds, $t_{cpu,a}$, (e) overall CPU time in seconds, $t_{cpu,t}$.}
	\label{fig:ch6:TNK-all}
\end{figure}

\begin{figure}[]
	\centering
	\includegraphics[height=0.21\textheight]{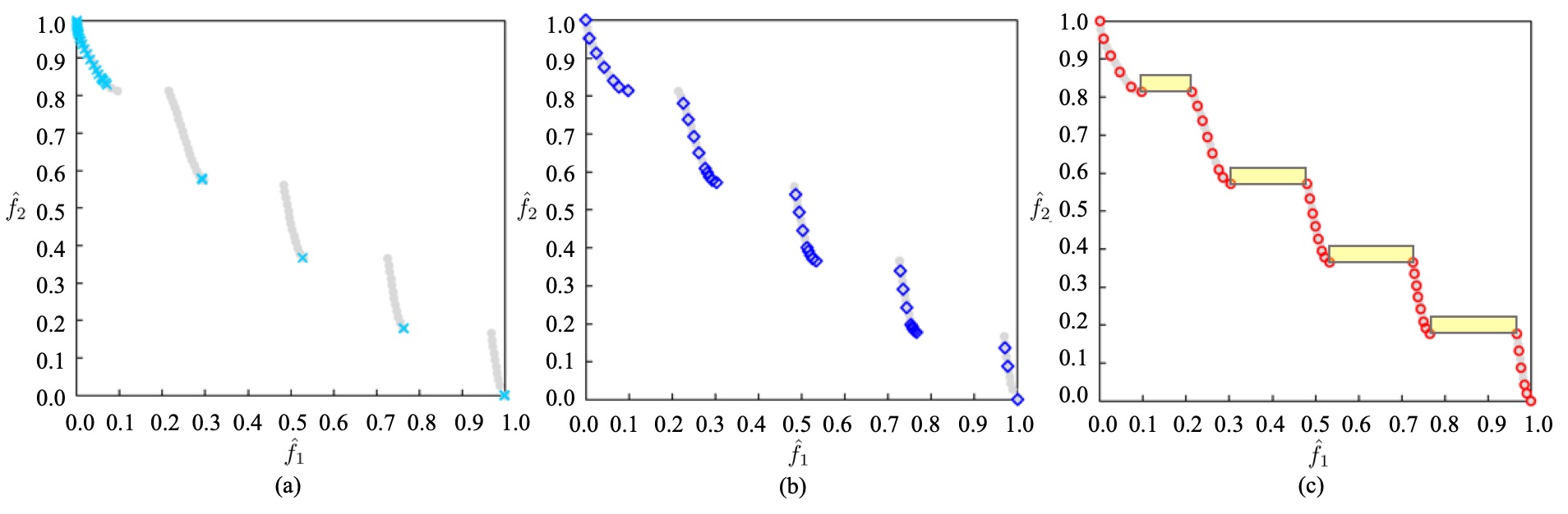}
	\vspace{-0.3cm}
	\captionsetup{font=small}\caption{Pareto points generated by (a) SD algorithm (\textcolor{skyblue}{$\boldsymbol{\times}$}), (b) mNBI (\textcolor{blue}{$\mathbf{\diamond}$}), and (c) SDNBI (\textbf{\textcolor{red}{\textbf{$\circ$}}}) for the test problem ZDT3, including the regions discarded as being devoid of Pareto points (yellow shading). The best-known Pareto front is shown in grey.}\label{fig:ch6:ZDT3-PF}
\end{figure}

\begin{figure}
	\centering
	\includegraphics[height=0.35\textheight]{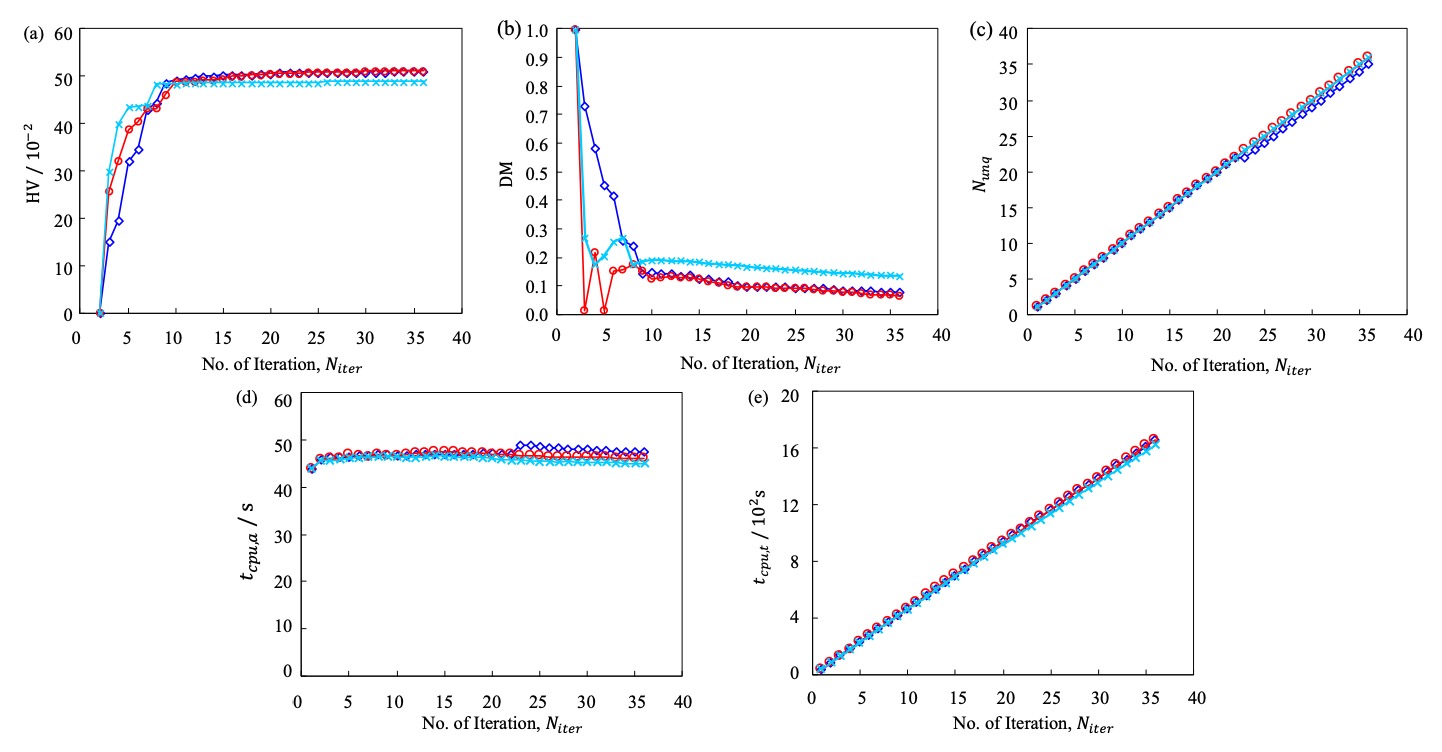}
	\vspace{-0.3cm}
	\captionsetup{font=small}\caption{Performance metrics as a function of iteration number for the application of the  SD algorithm (\textcolor{skyblue}{--$\times$--}), mNBI method (\textcolor{blue}{--$\diamond$--}), and SDNBI algorithm (\textbf{\textcolor{red}{--$\circ$--}}) to ZDT3 for 36 iterations: (a) Scaled hypervolume, HV, (b) Scaled distribution metric, DM, (c) Number of unique Pareto points, $N_{unq}$, (d) CPU time spent per iteration in seconds, $t_{cpu,a}$, (e) overall CPU time in seconds, $t_{cpu,t}$.}
	\label{fig:ch6:ZDT3-all}
\end{figure}

The overall results are shown in Tables \ref{Table:ch6:TNKZDT35}, and Figures \ref{fig:ch6:TNK-PF}-\ref{fig:ch6:ZDT3-all}. Note that results with the SD algorithm are not presented in Table \ref{Table:ch6:TNKZDT35} nor in Figure \ref{fig:ch6:TNK-all} for problem TNK, since there are no convex parts on the Pareto front and solutions other than the anchor points cannot be identified, as can be seen in Figure \ref{fig:ch6:TNK-PF}(a). 

From the tabulated and graphical results, we can conclude that the SDNBI algorithm outperforms the SD algorithm and the mNBI method in all metrics except for $t_{cpu,t}$ and $t_{cpu,a}$. The higher total computational cost of SDNBI may be explained by the additional calculation procedures for: the generation of inner and outer approximations of the Pareto front; the removal of \subspaces that do not contain Pareto points, and the addition of objective constraints inherent in the mNBI subproblem. At the 20$^{\rm th}$ iteration, the SDNBI algorithm identifies all disconnected regions of the Pareto front, with confidence that no Pareto points exist in these areas. By solving the \eqref{eq:ch6:simpleNBIn} subproblems, not only two empty regions are discarded, but new Pareto points are identified. As a result, the SDNBI algorithm succeeds in improving the Pareto approximation efficiently at each iteration, whilst the performance of the mNBI is limited when it encounters a large disconnected region. These aspects are illustrated in Figures \ref{fig:ch6:TNK-all} and \ref{fig:ch6:ZDT3-all}, where the progress of each quality measure can be compared. As can be seen in Figure \ref{fig:ch6:TNK-all}(d) the differences in $t_{cpu,a}$ between the algorithms decreases despite the rapid increase in $t_{cpu,t}$ with SDNBI (see Figure \ref{fig:ch6:TNK-all}(e)), since a larger number of Pareto points ($N_{unq}$) are identified with SDNBI as iterations proceed. In Figure \ref{fig:ch6:ZDT3-all}(d), lower values of $t_{cpu,a}$ are observed with SDNBI after the 25$^{\rm th}$ iteration compared to mNBI, as a result of identifying more Pareto points. This suggests that the SDNBI algorithm may be more computationally efficient when the Pareto front consists of many disconnected regions.

\begin{figure}[]
	\centering
	\includegraphics[height=0.21\textheight]{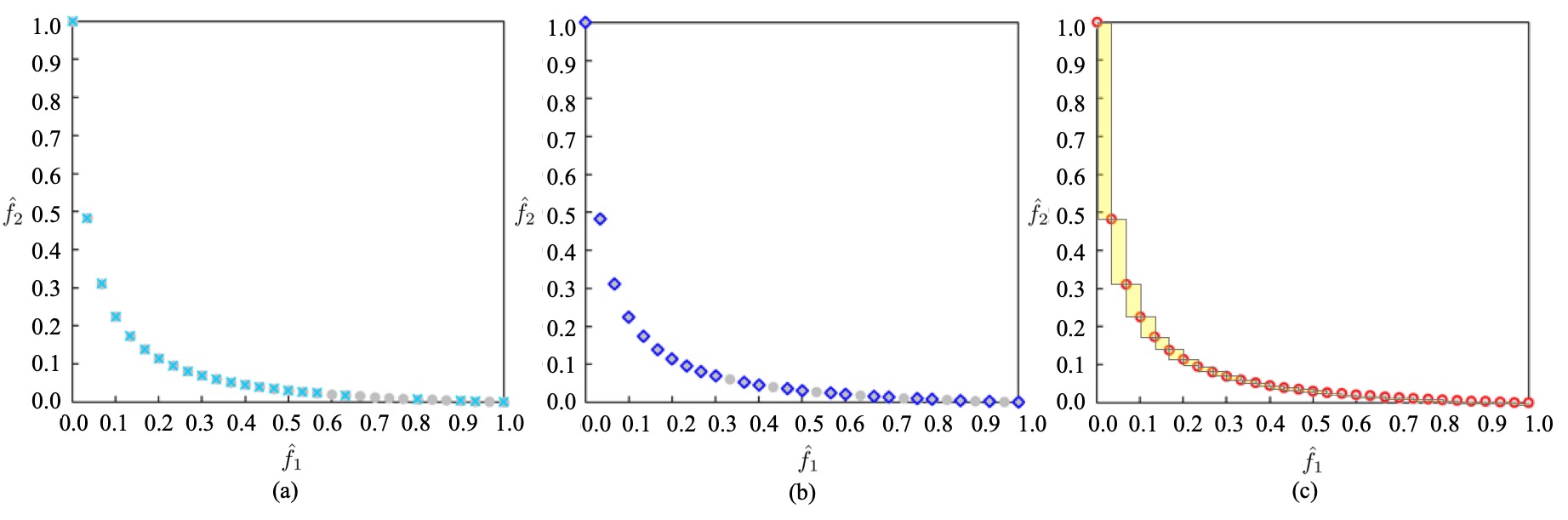}
	\vspace{-0.3cm}
	\captionsetup{font=small}\caption{Pareto points ({$\boldsymbol{z}$}) generated by the (a) SD algorithm (\textcolor{skyblue}{$\boldsymbol{\times}$}), (b) the mNBI method (\textcolor{blue}{$\mathbf{\triangle}$}), and (c) the SDNBI algorithm (\textbf{\textcolor{red}{\textbf{$\circ$}}}) for test problem ZDT5 with  $N_{iter,fix}=40$, including the regions discarded as being devoid of Pareto points (yellow shading). The grey points represent the set all Pareto points.}\label{fig:ch6:ZDT5-PF2}
\end{figure}

\begin{figure}[]
	\centering
	\includegraphics[height=0.35\textheight]{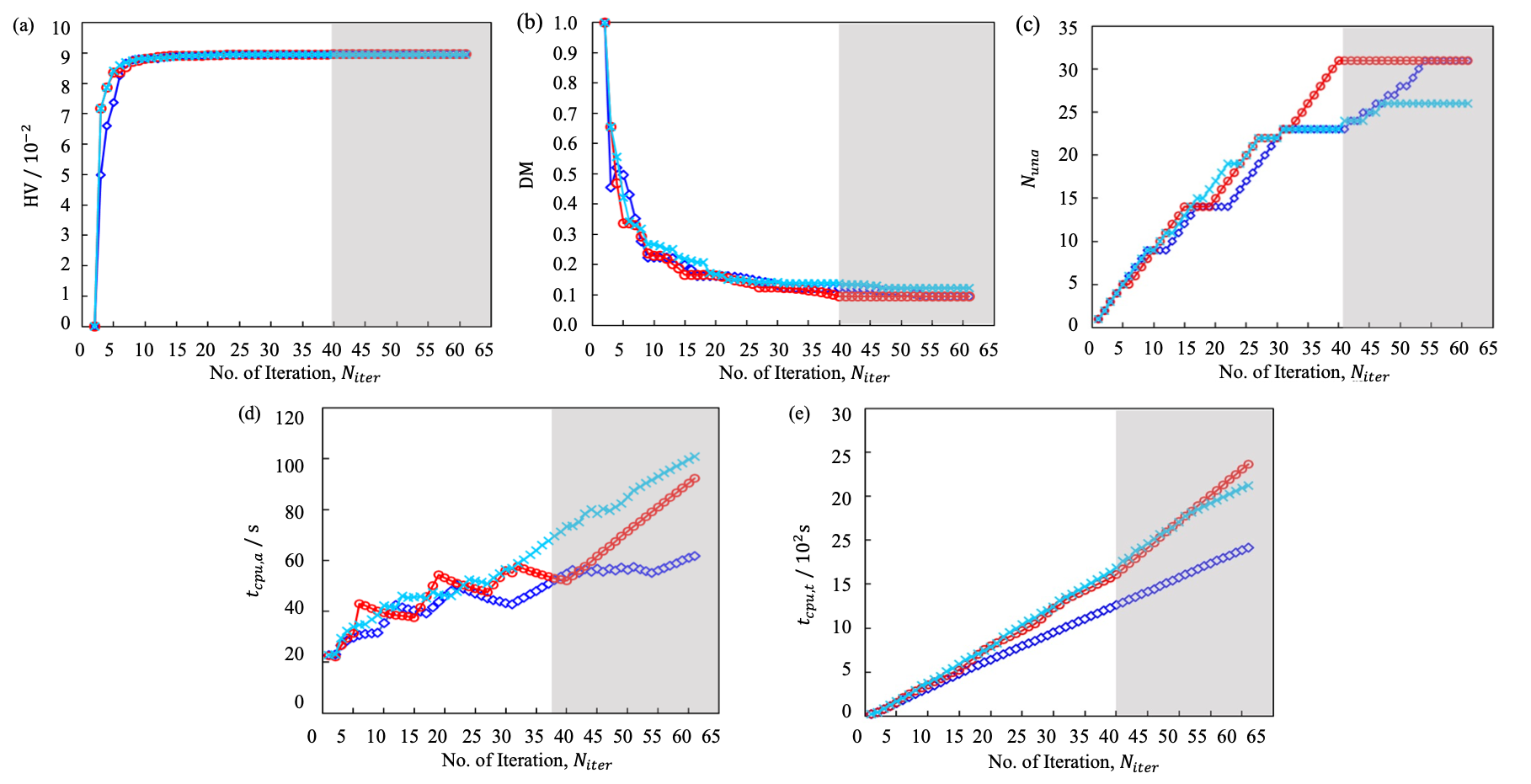}
	\vspace{-0.3cm}
	\captionsetup{font=small}\caption{Performance metrics as a function of iteration number for the application of the SD  (\textcolor{skyblue}{--$\times$--}), mNBI (\textcolor{blue}{--$\diamond$--}), and SDNBI algorithms (\textbf{\textcolor{red}{--$\circ$--}}) to ZDT5 for up to 64 iterations. The white area corresponds to the 40 iterations consistent with the SDNBI stopping criterion, while the grey shaded area corresponds to additional iterations needed for eliminating all empty area and so terminating the algorithm before reaching the stopping tolerance.}
	\label{fig:ch6:ZDT5-all}
\end{figure}
\subsection{Test problem ZDT5}
The results for test problem ZDT5 obtained with SD, mNBI and SDNBI  for 40 iterations are summarized in Table \ref{Table:ch6:TNKZDT35} and Figures \ref{fig:ch6:ZDT5-PF2} and \ref{fig:ch6:ZDT5-all}. For this test problem, the true Pareto front is represented as a set of 31 integer Pareto points and the SDNBI identifies the complete set of solutions for a given tolerance, while the mNBI method and the SD algorithm appear to be relatively ineffective in achieving high performance in DM and $N_{unq}$ with the given number of iterations. These results demonstrate the effectiveness of the SDNBI algorithm at finding a diverse and accurate set of Pareto points over the discrete bi-objective domain. As can be observed in Figure \ref{fig:ch6:ZDT5-all}(d), the performance of the SDNBI algorithm in terms of $t_{cpu,a}$ improves as iterations proceed with the lowest computational cost is obtained at iteration 40, by which convergence is reached. 

In order to investigate how the SDNBI algorithm works when the complete Pareto set has been identified, the convergence tolerance $\epsilon$ is reduced to 10$^{-3}$ and all algorithms are run for the number of iterations required to reach this tolerance with SDNBI. As can be seen in Figure \ref{fig:ch6:ZDT5-all}(d), the $t_{cpu,a}$ values increase rapidly with the SDNBI algorithm as no new Pareto points can be found following the identification of the complete set of Pareto points (at 40$^{\rm th}$ iteration). This additional cost is associated with the solution of the subproblems \eqref{eq:ch6:simpleNBIn} and \eqref{eq:ch6:simpleNBIn-b} through which facets are removed, but no additional Pareto points are identified. It is notable that all regions where no Pareto point exists are removed from the search space after 61 iterations (see Figure \ref{fig:ch6:ZDT5-all}). This allows the SDNBI algorithm to terminate although the pre-defined error tolerance, $\epsilon$, is not satisfied. Although not shown here, without setting a maximum number of iterations, the SD algorithm would iterate until there are no remaining facets to be investigated (at the 66$^{\rm th}$ iteration).

As can be observed from Figure \ref{fig:ch6:ZDT5-all}(c), by iteration 56, all 31 Pareto points have been identified with the mNBI method. However, $\boldsymbol{\beta}$ values could be enumerated until the  87$^{\rm th}$ iteration, when the set of possible choices of reference points, $\boldsymbol{\Phi\beta}$, becomes almost empty. This suggests that the mNBI method can fail to recognize that all Pareto points have been found when a highly disconnected or discrete Pareto front is involved, due to the use of an arbitrary termination criterion.

\section{Application to molecular and process design problems}
In previous sections, the SDNBI algorithm was evaluated through the solution of five standard mathematical benchmark problems. To further assess the practical relevance and robustness of the algorithm, we now consider two chemically motivated case studies (CSs) drawn from computer-aided molecular (and process) design (CAM(P)D). These problems are substantially more complex than the synthetic benchmarks, formulated as MINLP with discrete variables, nonlinear constraints, and nonconvex trade-offs arising from thermodynamic models, molecular structure–property relationships, and process design constraints. They serve as a testbed for evaluating BOO algorithms.

The CSs considered here are based on previously published CAM(P)D formulations from \cite{Lee2020}. Case Study 1 (CS1) is derived by reformulating a four-objective problem (Equation (6) in that work) into a BOO problem. In Case Study 2 (CS2),   the bi-objective formulation of Case Study 3 (Equation (8)) in the same work is adopted, with minor modifications for compatibility with the scalarization procedures of the SDNBI, mNBI and SD algorithms. Readers are referred to \cite{Lee2020} for the full problem descriptions, including molecular design space, property models, and constraint bounds. Details of the mathematical formulations of the two case studies and complementary information on solution strategies are provided in Section S4 of Supplementary Materials.

\subsection{Case study 1: Solvent design for chemical absorption of CO$_2$}
A CAMD application for the design of optimal amine-based solvents for chemical absorption of CO$_2$ is considered in CS1. Two performance metrics are selected as objective functions: the liquid-phase heat capacity, $C_P$, which correlates with the sensible heat required for solvent regeneration, and the relative energy difference, $RED$, which correlates inversely with CO$_2$ loading capacity, i.e., the ability to dissolve gaseous CO$_2$. These objectives are chosen based on prior analysis of the original multi-objective formulation. Analysis of the reduced problem shows that the corresponding Pareto front is highly nonconvex, making it a suitable test for evaluating the performance of the SDNBI algorithm.
 The resulting optimization problem is as follows:
\begin{equation}\label{eq:ch7:test1-CO2} \tag{PCS1}
\begin{aligned}
\min_{\boldsymbol{n}} \quad & {C_{P}}, RED \\ 
\textrm{s.t.} \quad & \boldsymbol{g}(\boldsymbol{n})\leq \textbf{0} \\
& \boldsymbol{h}(\boldsymbol{n}) = \textbf{0} \\
& \boldsymbol{Cn}\leq \boldsymbol{d}
\end{aligned}
\end{equation}
where $\boldsymbol{n}\in \boldsymbol{N}\subseteq \mathbb{N}^{n_1}$ is a $n_1$-dimensional vector of non-negative integer variables that make up the molecular design space. The molecular design space is defined over 13 functional groups: NH$_2$CH$_2$, NH$_2$CH, NH$_2$C, NHCH$_3$, NHCH$_2$, NHCH, NCH$_3$, NCH$_2$, CH$_3$, CH$_2$, CH, C, and OH. $\boldsymbol{g}(\boldsymbol{n})$ is a vector of inequality constraints on the density ($\rho$), the normal melting point ($T_m$), the normal boiling point ($T_b$), the viscosity ($\mu$), and the surface tension ($\sigma$), $\boldsymbol{h}(\boldsymbol{n})$ is a vector of equality constraints that include structure-property models that link the structure of the solvent to its thermophysical properties, and $\boldsymbol{Cn}\leq \boldsymbol{d}$ is a set of linear constraints that represent molecular feasibility constraints so that only chemically feasible molecules are considered in the optimization. The structure-property models are group contribution methods expressed as a set of linear and nonlinear algebraic equations. 

\subsection{Case study 2: Integrated working fluid and ORC process design}
CS2 involves the simultaneous design of a pure working fluid and an Organic Rankine Cycle (ORC). The net power output ($P_{net}$) and the total cost of investment ($TCI$) are considered as the two objective functions in order to investigate trade-offs between the efficiency of the working fluid in extracting power from the ORC for a specified heat source and heat sink, and the overall economic performance of the combined working fluid-ORC process. The problem formulation thus includes structure-property models as well as a process model of the ORC. The optimization problem is expressed as:
\begin{equation}\label{eq:ch7:test2-ORC} \tag{PCS2}
\begin{aligned}
\min_{\boldsymbol{x,n}} \quad &  -{P_{net}, TCI}\\
\textrm{s.t.} \quad & \boldsymbol{g}(\boldsymbol{x,n}) \leq \textbf{0}   \\
& \boldsymbol{h}(\boldsymbol{x,n}) = \textbf{0} \\
& \boldsymbol{Cn} \leq \boldsymbol{d} \\
\end{aligned}
\end{equation}
where  $\boldsymbol{x}$ is a $n_2$-dimensional vector of continuous variables that represent process conditions, including the mass flow of working fluid ($m_{WF}$), the pump inlet pressure ($P_{in,pump}$), the pump outlet pressure ($P_{out,pump}$), and the extent of superheating in an evaporator ($\Delta T_{superheat}$). The molecular design space is defined over nine functional groups: CH${_3}$, CH${_2}$, CH${_2}$=, CH=, eO (end group oxygen), cO (central oxygen), COOH, COO and OH. eO represents an oxygen atom connected to one CH${_3}$ and one CH${_2}$ group, and cO describes an oxygen atom bonded to two CH${_2}$ groups. The constraints $\boldsymbol{g}(\boldsymbol{x,n})$ and $\boldsymbol{h}(\boldsymbol{x,n})$ define process feasibility, thermodynamic relations, and other structural-property relationships. The model is significantly more challenging to evaluate than that is CS1 due to the need to include nonideal vapour-liquid equilibrium calculations, with key thermodynamic properties calculated with the SAFT-$\gamma$ Mie group-contribution method \cite{Papaioannou2014, haslam2020expanding}.   

For Case Studies 1 and 2, while the overall algorithmic framework is implemented in MATLAB, the scalarized (single-objective) subproblems are solved using a standard outer approximation algorithm with equality relaxation and augmented penalty (OA/ER/AP) \citep{viswanathan1990combined}. In this approach, each MINLP subproblem is decomposed into a nonlinear programming (NLP) subproblem, referred to as the primal problem, and a mixed-integer linear programming (MILP) subproblem, referred to as the master problem. The primal problems are evaluated via a C++ interface to gPROMS ModelBuilder 7.0.7 \citep{gproms}, and the master problems are solved using Gurobi 9.0 \citep{gurobi2020}. All runs are performed on a single Intel(R) Xeon(R) Gold 5122 CPU @ 3.60GHz processor with 384 GB of RAM. 

The objective functions are normalized using the ideal and nadir values listed in Table \ref{table:ch7:nadir_utopia}. While the global search strategy and stopping criteria are consistent with those used for the benchmark tests, several implementation aspects---including constraint handling and infeasibility resolution---are tailored to the CAM(P)D applications. Full solver specifications and implementation details are provided in Section S4 of the Supplementary Materials.

\begin{table}[htb]\centering
\caption{Nadir and utopia point used for the normalization of the objective functions in CS1 and CS2}
\label{table:ch7:nadir_utopia}
\renewcommand{\arraystretch}{0.95}
\renewcommand{\tabcolsep}{1mm}
\begin{tabular}{l|cclcc}
\hline
\multicolumn{1}{c|}{$i$} & \begin{tabular}[c]{@{}c@{}}$C_p$\\ (kJ mol-K$^{-1}$)\end{tabular} & \multicolumn{2}{c}{\begin{tabular}[c]{@{}c@{}}$RED$\\ (-)\end{tabular}} & \begin{tabular}[c]{@{}c@{}}$P_{net}$\\ (MW)\end{tabular} & \begin{tabular}[c]{@{}c@{}}$TCI$\\ (10$^6$ US\$)\end{tabular} \\ \hline
$f_j^{nd}$  & 0.434          & \multicolumn{2}{c}{4.181}                                             & 5.01$\times$10$^{-8}$                                             & 7.50                                                     \\
$f_j^{id}$               & 0.139                                                     & \multicolumn{2}{c}{0.226}                                             & 31.62                                                & 6.00$\times$10$^{-2}$                                                 \\ \hline
\end{tabular}%
\end{table}

\subsection{True/Best-known Pareto fronts}
\begin{figure}
	\centering
	\includegraphics[height=0.26\textheight]{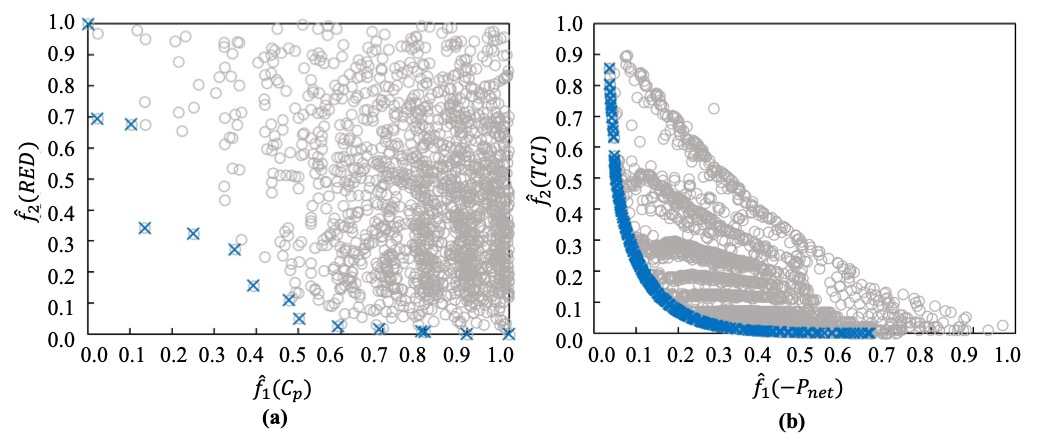}
	\vspace{-0.2cm}
	\captionsetup{font=small}\caption{True/best-Known Pareto front (blue marker,\textcolor{frenchblue}{$\times$}) and feasible points (grey marker,\textcolor{gray}{ $\circ$}) of (a) the normalized heat capacity, $C_p$ ($\hat{f}_{1}$) versus the normalized relative energy difference, $RED$ ($\hat{f}_{2}$) of CS1; and (b) the normalized net power output, $-P_{net}$ ($\hat{f}_{1}$) versus the normalized total capital investment, $TCI$ ($\hat{f}_{2}$) of CS2, resulting from the enumeration of all possible combination of molecules.}
	\label{fig:ch7:trueall}
\end{figure}

The performance of the SD, mNBI, and SDNBI algorithms on the two CAM(P)D CSs is evaluated relative to the reference Pareto fronts shown in Figure \ref{fig:ch7:trueall} and Table \ref{table:ch7:truePareto_quality}. For CS1, the true Pareto front was obtained by enumerating all possible combinations of molecular structures. For C2, an exhaustive solution is infeasible due to the presence of continuous process variables. Instead, the best-known Pareto front is approximated by combining three sets of non-dominated points: (i) results from 123 independent mNBI runs, each with a distinct $\boldsymbol{\beta}$ vector, (ii) solutions generated by SD and SDNBI runs across the entire search space, and (iii) the non-dominated solutions reported in \cite{Lee2020}, obtained via exhaustive scalarization of 267 feasible working fluids. The merged set yields 326 unique non-dominated points and is used as a reference for all comparative evaluations.
A graphical representation of the true/best-known Pareto front for each case study is shown in Figure \ref{fig:ch7:trueall}, in the space of objective functions. Quantitative characteristics of the true/best-known Pareto front obtained for each test is provided in Table \ref{table:ch7:truePareto_quality}. As can be seen in Figure \ref{fig:ch7:trueall}(a), each point in the feasible region corresponds to a distinct solvent structure, meaning that the Pareto front for CS1 is entirely disjoint. In CS2, one can observe two disconnected parts, corresponding to points where one of the objective constraints of the mNBI method become inactive, implying that there be two disconnected regions.

The list of Pareto-optimal molecules for each case study is given in Tables S2 and S3 in the Supplementary Material. It should be noted that each Pareto-optimal working fluid in CS2 typically appears in more than two Pareto points, characterized by different optimal process conditions. Therefore, we only present the distinct molecular structures in Table S2 and omit objective function values.

\begin{table}[]\centering
\caption{Enumeration of all solutions for all possible combinations of functional groups in each test. $N_{structure}$: the number of molecules that satisfy structural feasibility, $N_{feasible}$: the number of molecules that satisfy all constraints. The value in (.) of $N_{unq}$ of CS2 is the number of different molecular structures contained within this set.}
\label{table:ch7:truePareto_quality}
\begin{tabular}{lcc}
\hline
                & CS1              & CS2      
                \\ \hline
$N_{structure}$ & 67,664              & 3,175 \\
$N_{feasible}$  & 3,160               & 267   \\
$N_{unq}$        & 15                   & 326 (6)    \\
HV ($10^{-2}$)  & 79.12               & 91.62    \\
DM              & 0.1366              & 0.0144 \\ \hline
\end{tabular}
\end{table}

\subsection{Results for CS1}
\begin{figure}[h]
	\centering
	\vspace{4cm}
	\includegraphics[height=0.21\textheight]{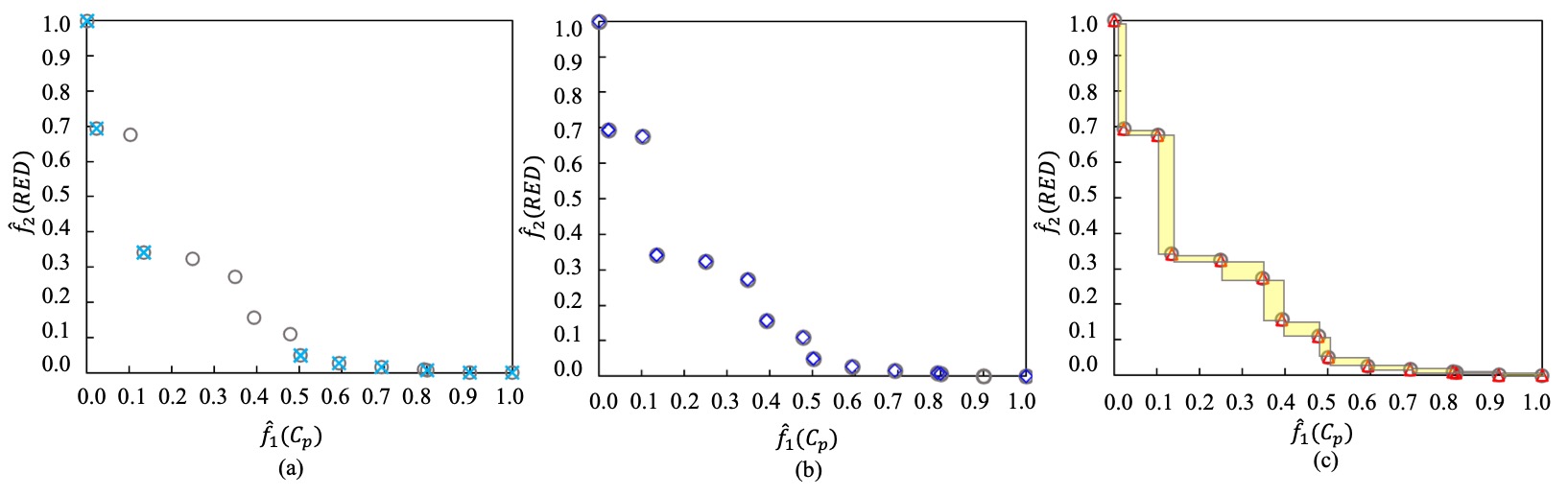}
	\vspace{-0.3cm}
	\caption{Pareto points generated by (a) the SD algorithm (\textcolor{cyan1}{$\times$}), (b) the mNBI method (\textcolor{blue}{$\boldsymbol{\diamond}$}), and (c) the SDNBI algorithm (\textcolor{red}{\textbf{$\boldsymbol{\triangle}$}}) for CS1, including the regions discarded as being devoid of Pareto points (yellow shading). The best-known Pareto front is shown in grey.}
	\label{fig:ch7:test1_result}
\end{figure}

The results obtained for CS1 after a fixed number of iterations for each algorithm are summarized in Table \ref{table:ch7:test1_result} and Figure \ref{fig:ch7:test1_result}. It can be seen from Figure \ref{fig:ch7:test1_result}  that both the SDNBI and mNBI methods generate well-distributed and globally Pareto-optimal solutions, whereas the SD algorithm fails to produce an accurate representation of the Pareto front. This is because the weighted-sum subproblem over the large nonconvex part of the Pareto front, converges to the previously identified optimal solution, forcing the algorithm to premature termination. The corresponding HV and DM values show a large discrepancy relative to the reference quantities.

None of the algorithms achieves convergence under the maximum error criterion as the Pareto front is entirely disjoint, and the points are separated by large gaps. This highlights the importance of identifying subspaces that contain no Pareto points in order to reduce unnecessary computational effort in attempting to improve the approximation quality. The SDNBI method achieves this by objective-space fathoming and terminates at iteration 27 with the highest HV and lowest DM. The higher computational cost of SDNBI is partly attributed to the solution of 14 additional subproblems \eqref{eq:ch6:simpleNBIn} and \eqref{eq:ch6:simpleNBIn-b}. This overhead also reflects the added constraints in the SDNBI formulation, which complicate the identification of feasible regions. As a result, subproblem evaluations may require more effort to locate feasible points, and additional optimization runs are often triggered during the global search phase when local runs fail due to infeasible starting points.

\begin{figure}[h]
	\centering
	\includegraphics[height=0.22\textheight]{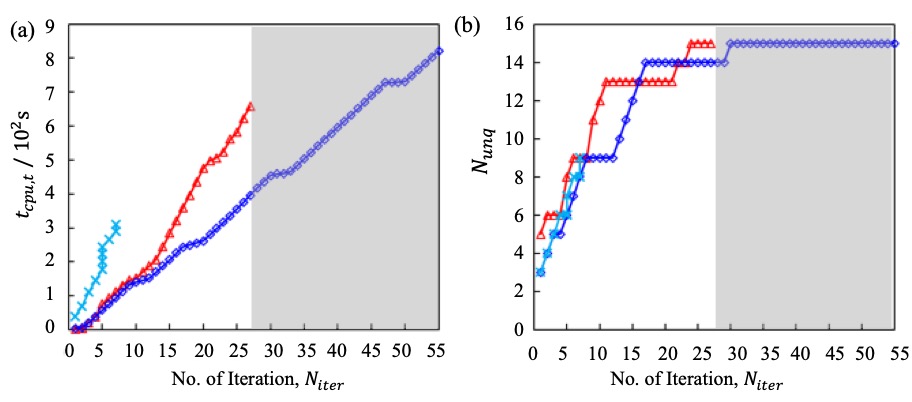}
	\vspace{-0.3cm}
	\caption{(a) Overall CPU time in seconds, $t_{\mathrm{cpu},t}$, and (b) Number of unique Pareto points, $N_{\mathrm{unq}}$, as a function of iteration number for the application of the SD algorithm (\textcolor{cyan1}{–$\times$–}), mNBI method (\textcolor{blue}{–$\boldsymbol{\diamond}$–}), and SDNBI algorithm (\textcolor{red}{–$\boldsymbol{\triangle}$–}) to CS1 up to 55 iterations. The white region corresponds to the 27 iterations completed by the SDNBI algorithm, which terminated after identifying all empty objective-space regions. The grey-shaded area shows additional iterations taken by the mNBI method, which continued despite no new Pareto points being found.}
	\label{fig:ch7:test1_result_cpu}
\end{figure}

Despite the higher computational cost, a key advantage of the SDNBI algorithm over mNBI and SD is its ability to discard regions of the objective space that do not yield non-dominated solutions. This leads to a more efficient search and a more informative representation. In the case of the mNBI method, the algorithm continues until the set of $\boldsymbol{\beta}$ is exhausted, unless additional stopping criteria, such as the progression metric \citep{van1999multiobjective}, are employed. To investigate such algorithmic behavior, we consider an extended set of $\boldsymbol{\beta}$. The cumulative CPU time and number of Pareto points obtained as a function of iteration are shown in Figure \ref{fig:ch7:test1_result_cpu}. At the 30$^{\rm th}$ iteration, all 15 Pareto-optimal solvents are identified by the mNBI method. However, a higher total CPU time is observed with mNBI as the iteration proceeds because the objective search regions that do not include any new Pareto points are only partially eliminated by mNBI and thus the algorithm iterates over the remaining search domain in an attempt to achieve the pre-defined approximation quality tolerance. These results indicate that when the Pareto front is highly disconnected, mNBI may be inefficient, and further iterations do not necessarily yield improvements in approximation quality.

\begin{table}[t]\centering
\caption{Performance metrics for molecular design problems: CS1 with $N_{iter,fix}=27$, and CS2 with $N_{iter,fix}=61$, comparing the SD, mNBI, and SDNBI algorithms}
 \resizebox{0.7\textwidth}{!}{%
\begin{tabular}{l|ccc|ccc}  
\hline
Problem  & \multicolumn{3}{c|}{CS1} & \multicolumn{3}{c}{CS2}  \\ \hline
Algorithm & SD       & mNBI      & SDNBI     & SD       & mNBI      & SDNBI  \\ 
$N_{unq}$   & 9        & 14        & 15       & 59(3)      & 58(4)        & 60(4)          \\
HV(10$^{-2}$)& 76.43 & 79.05    & 79.12   & 91.23   & 91.47    &91.44      \\
DM         & 0.2323   & 0.1521    & 0.1366   & 0.1130    & 0.0472   & 0.0388            \\
$t_{cpu,a}$ (s)  & 34.55   & 28.31   & 43.77  & 30.97    & 144.99     &115.78  \\
$t_{cpu,t}$ (s) & 3.11$\times$10$^2$ & 3.96$\times$10$^2$   & 6.58$\times$10$^2$  & 1.68$\times$10$^3$    & 8.41$\times$10$^3$ & 6.95$\times$10$^3$     \\ \hline
\end{tabular}%
}
\label{table:ch7:test1_result}
\end{table}

\subsection{Results for CS2}
\begin{figure}[h]
	\centering
	\includegraphics[height=0.21\textheight]{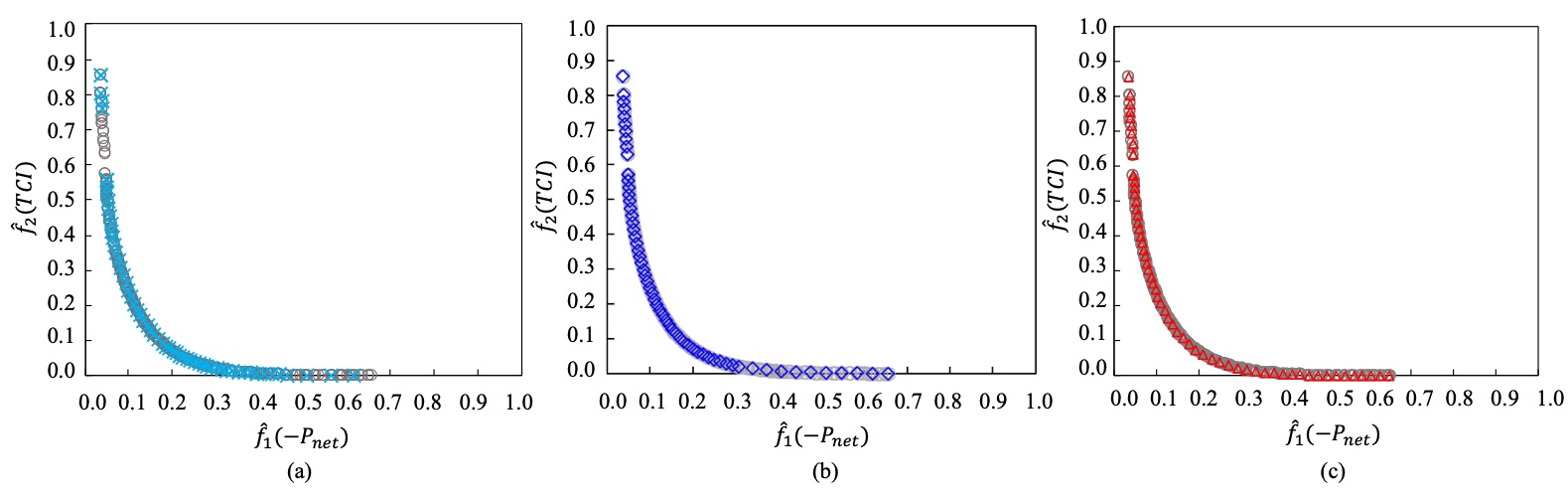}
	\vspace{-0.7cm}
	\caption{Pareto points generated by (a) the SD algorithm (\textcolor{cyan1}{$\times$}), (b) the mNBI algorithm (\textcolor{blue}{$\boldsymbol{\diamond}$}), and (c) the SDNBI algorithm (\textcolor{red}{\textbf{$\boldsymbol{\triangle}$}}) for CS2. The best-known Pareto front is shown in grey.}
	\label{fig:ch7:test2_result}
\end{figure}

The Pareto points obtained with each algorithm are shown in Figure \ref{fig:ch7:test2_result} and compared with the best-known Pareto front. The quality measures, calculated based on the fixed number of iterations, are summarized in Table \ref{table:ch7:test1_result}. It is clear from Figure \ref{fig:ch7:trueall} that most parts of the Pareto front are characterized as convex, allowing the SD algorithm to access a diverse set of Pareto-optimal solutions along the Pareto surface. The differences in HVs (91.23$\times$10$^{-2}$ for SD, 91.47$\times$10$^{-2}$ for mNBI, and 91.44$\times$10$^{-2}$ for SD) are indeed very small, and all values are close to the reference value (91.62$\times$10$^{-2}$). Similar values of $N_{unq}$ are also achieved. This is also explained in the same manner, although the SDNBI algorithm is marginally better at identifying more diverse combinations of working fluids and thermal-economic performance. 

As shown in Table \ref{table:ch7:test1_result}, the mNBI and SDNBI algorithms exhibit higher computational cost than the SD algorithm. As in CS1, this is mainly due to the increased numerical complexity of the mNBI/SDNBI subproblems, where two additional objective constraints are introduced into the original MINLP formulation. The increase in computational time is more pronounced in CS2, mainly due to the evaluation cost of the ORC process models, which involve complex mass- and heat-balance equations. In such cases, the performance of the OA-based framework becomes highly sensitive to the quality of starting points. When many molecular candidates fail to satisfy the added constraints, the algorithm requires more iterations to locate feasible regions or may fail entirely. Furthermore, since starting points are sampled in the space of molecular structures, continuous variables are determined by solving the primal problem with fixed discrete variables. This further increases computational burden when starting from infeasible points, requiring either more samples to be evaluated or additional time for the primal problem to converge.

Among the three methods, the SDNBI algorithm produces the lowest DM value, indicating the most uniform and accurate coverage of the Pareto front. This is visually evident in Figure \ref{fig:ch7:test2_result}, where the solutions from SDNBI exhibit broader and more uniform coverage. In contrast, the approximation by the SD algorithm reveals sparsely populated regions within $0.48 \leq f_1 \leq 0.6$ and $0.59 \leq f_2 \leq 0.7$, suggesting the presence of nonconvex or disconnected segments of the Pareto front. A notable advantage of SDNBI lies in its ability to identify and characterize such disconnected regions, thereby providing insight into the shape and structure of the true Pareto front. For example, the steep surface around both anchor points suggests that considerable economic improvement can be achieved with only modest reductions in power output. The front obtained with the SD algorithm, by contrast, suggests limited trade-off potential in that region.

\begin{figure}[h]
	\centering
	\includegraphics[height=0.25\textheight]{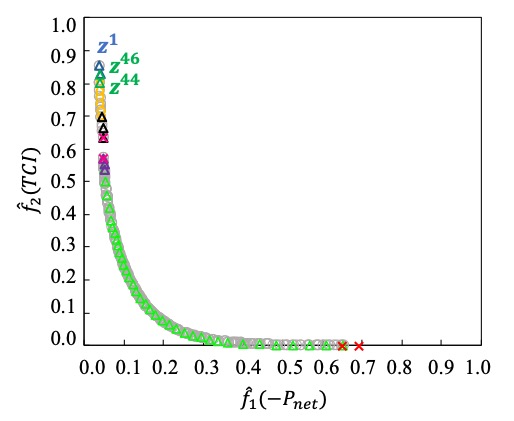}
	\vspace{-0.2cm}
	\caption{Schematic of the solutions (\textbf{$\boldsymbol{\triangle}$}) generated by SDNBI in a bi-objective space for CS2. Each color denotes each of eight different subspaces. The empty regions between the cross marker(\textcolor{red}{$\boldsymbol{\times}$}) are identified by solving SDNBI subproblems \eqref{eq:ch6:simpleNBIn} and \eqref{eq:ch6:simpleNBIn-b}. Solution $\boldsymbol{z^{46}}$ is identified as a dominated solution after investigating if it is a member of the best-known Pareto front.}
	\label{fig:ch7:test2_subspace}
\end{figure}

At the termination of the SDNBI algorithm, the objective space is decomposed into a total of eight subspaces, and they are shown in Figure \ref{fig:ch7:test2_subspace}, where the Pareto points included in each subspace are depicted with different colors. As denoted with cross markers, two subspaces are successfully eliminated from consideration by solving subproblems \eqref{eq:ch6:simpleNBIn} and \eqref{eq:ch6:simpleNBIn-b}, while such regions are only partially screened out when using mNBI.  The discontinuity of the Pareto front observed at $0.8022<f_2 <0.8547$ is misinterpreted and a sub-optimal, i.e., dominated solution, $\boldsymbol{z^{44}}$, is produced. The facets represented by $F^1_S(\boldsymbol{z^1, z^{46}})$ and $F^2_S(\boldsymbol{z^{46}, z^{44}})$ in Figure \ref{fig:ch7:test2_subspace}, which are constructed by the convex hull generation are discarded from the search space in a subsequent iteration, as they satisfy the quality tolerances rather than because an empty part of the Pareto front has been recognized. Similarly, mNBI and SD each produce one dominated solution. This suboptimality is identified by comparing the set of Pareto points generated by each algorithm with the elements of the best-known Pareto front, not within the approximated set. Overall, the comparative analysis confirms that the SDNBI algorithm shows excellent performance in generating a well-spaced and accurate set of Pareto-optimal working fluids and process conditions, within a reasonable computational budget.

\section{Conclusions}\label{sec:conclusion}
A novel algorithm, SDNBI, that combines features of the SD and NBI methods with the aim of generating accurate approximations of nonconvex and/or discontinuous Pareto fronts in bi-objective space has been presented. The main benefits of the proposed algorithm are focused on the exploration of nonconvex regions of the Pareto front and on the identification of regions where no further nondominated solution exists, indicating that parts of the Pareto front are disconnected. Key aspects of the SDNBI algorithm include: 1) the characterization of inner and outer approximations such that the accuracy of the incumbent set of Pareto points can be assessed during the course of the algorithm, 2) a systematic way of setting the scalarization parameters of the mNBI subproblem, 3) the decomposition of an objective space based on the convexity and nonconvexity of parts of the Pareto front, and 4) the refinement of the scalarized subproblem to avoid unnecessary iterations over disconnected or empty regions of the objective space.

To assess the performance of the proposed algorithm,  numerical tests were conducted for five bi-objective benchmark problems (MOP1, SCH2, TNK, ZDT3, and ZDT5) and two case studies on molecular and process design. The performance of the algorithm in terms of the accuracy of the approximation of the Pareto front constructed in disconnected nonconvex objective domains was compared with two MOO approaches: the SD algorithm and the mNBI method. This provided clear evidence of the effectiveness of the SDNBI algorithm in generating a more diverse and better-distributed set of Pareto points over disconnected and/or nonconvex Pareto fronts compared to the approximations generated with the mNBI method and the SD algorithm. The solution of test problem ZDT5, in which the true Pareto front comprises 31 discrete points, highlighted the robustness of the SDNBI algorithm:  the complete set of  Pareto points was generated with relatively few iterations, leading to high computational efficiency. The SDNBI not only provides information on the quality of the Pareto front approximation (as does the SD algorithm), but it also reveals areas of bi-objective space that are devoid of Pareto points.

The computational studies on the highly constrained molecular and process design problems further highlighted the effectiveness of the SDNBI algorithm on engineering applications. The ability to identify and eliminate regions of objective space that do not contain any Pareto points, a unique feature of the proposed algorithm, was found to be especially beneficial in reaching convergence.  

Further work will focus on improving algorithmic performance and testing the SDNBI algorithm on additional literature and engineering problems. The extension of the approach to handle problems with three or more objective functions will also be considered in future work as such problems are often relevant to practical applications of multi-objective optimization.\\

\textbf{Acknowledgements}\\
Ye Seol Lee gratefully acknowledges financial support from the Fund For Women Graduates (FfWG) of the British Federation of Women Graduates and from  Imperial College London via a Roger Sargent scholarship.\\

\textbf{Data Availability}\\
Code for the SD algorithm, mNBI method, and SDNBI is available at \url{https://github.com/MolChemML/MOO-SDNBI}, and all data pertaining to the results are available on the Zenodo repository (\url{https://doi.org/10.5281/zenodo.11223677}), under CC-BY licence.

\begin{spacing}{1.42} 
\bibliography{SDNBI_ref}
\end{spacing}
\end{document}